\documentclass[12pt]{amsart}

\usepackage{amssymb,amsfonts}

\usepackage{amsthm}
\usepackage[dvipsnames]{xcolor}
\usepackage{graphicx}
\usepackage{xfrac}
\usepackage[all,arc]{xy}
\usepackage[colorlinks=true, allcolors=blue]{hyperref}
\usepackage[only,llbracket,rrbracket]{stmaryrd}
\usepackage{enumerate}
\usepackage{mathrsfs}
\usepackage{tikz-cd}
\usepackage{import}
\usepackage[utf8]{inputenc}
\usepackage{hyperref}
\usepackage{cancel}
\usepackage{bbm}
\usepackage{comment}
\usepackage[left=1in,top=1in,right=1in,bottom=1in]{geometry}
\usepackage{mathtools}

\usepackage{amsthm}
\theoremstyle{plain}
\newtheorem{theorem}{Theorem}[section]
\newtheorem{lemma}[theorem]{Lemma}
\newtheorem{prop}[theorem]{Proposition}
\newtheorem{cor}[theorem]{Corollary}

\theoremstyle{definition}
\newtheorem{definition}[theorem]{Definition}
\newtheorem{definition-lemma}[theorem]{Definition/Lemma}
\newtheorem{example}[theorem]{Example}

\newtheorem{observation}[theorem]{Observation}
\newtheorem{problem}[theorem]{Problem}
\newtheorem{hypothesis}[theorem]{Hypothesis}

\theoremstyle{remark}
\newtheorem{remark}[theorem]{Remark}

\newcommand{\F}{\mathbb{F}}
\newcommand{\Z}{\mathbb{Z}}
\newcommand{\C}{\mathbb{C}}
\newcommand{\N}{\mathbb{N}}
\newcommand{\R}{\mathbb{R}}
\newcommand{\Q}{\mathbb{Q}}

\newcommand{\ol}[1]{\overline{#1}}

\newcommand{\longto}{\longrightarrow}

\newcommand{\Hom}{\text{Hom}}

\newcommand{\Aut}{\mathrm{Aut}}

\newcommand{\coker}{\mathrm{coker }}

\newcommand{\std}{\operatorname{std}}

\newcommand{\into}{\hookrightarrow}
\newcommand{\onto}{\twoheadrightarrow}

\renewcommand{\H}{\mathrm{H}}

\newcommand{\desus}{\mathbf{s}^{-1}}

\newcommand{\gr}[1][]{\mathrm{gr}_{#1}~}
\newcommand{\wt}{\operatorname{wt}}
\newcommand{\tot}{\operatorname{tot}}
\newcommand{\defect}{\operatorname{def}}

\newcommand{\maps}{\colon}

\newcommand{\other}{\mathrm{otherwise}}

\DeclareMathOperator{\Id}{Id}

\renewcommand{\defect}{\operatorname{def}}

\newcommand{\lead}{\operatorname{lead}}

\begin{document}
%Header-Make sure you update this information!!!!

\noindent
\title[Letter-braiding: bridging group theory and topology]{
  Letter-braiding: bridging combinatorial group theory and topology
  }

\author{Nir Gadish}
\address{Department of Mathematics, University of Pennsylvania, Philadelphia, PA}
\email{\url{ngadish@math.upenn.edu}}



\begin{abstract}
We define invariants of words in arbitrary groups, measuring how letters in a word are interleaving, perfectly detecting the dimension series of a group.
These are the letter-braiding invariants. On free groups, braiding invariants coincide with coefficients in the Magnus expansion.  In contrast with Magnus' coefficients, our invariants are defined on all groups and over any PID. They respect products in the group and are a complete invariant of the
dimension series, so they are the coefficients of a universal multiplicative finite-type invariant, depending functorially on the group.

Letter-braiding invariants arise from the bar construction on a cochain model of a space with a prescribed fundamental group. This approach specializes to simplicial presentations of a group as well as to more geometric contexts, which we illustrate in examples.

As an application, we define variants of the Johnson filtration and the Johnson homomorphism on the automorphisms of arbitrary groups, and use them to constrain automorphisms of finite $p$-groups.
\end{abstract}
\maketitle

\section{Introduction}
How can one detect whether a particular element in a group belongs to the $k$-th term in the lower central series? 
%And, given two such elements, how can one distinguish them?
In this work we describe a collection of topologically motivated invariants of words in arbitrary groups, called \emph{letter-braiding invariants}.  They are explicitly computable for every word, given a presentation of the group. The invariants distinguish between linear combinations of group elements with coefficients over any PID up to powers of the augmentation ideal in the group ring, and thus give obstructions to being a $k$-fold commutator. This theory is, in a sense, a common generalization of Chen's iterated integrals \cite{chen1971pi1} and Fox' free differential calculus \cite{Fox1} -- both deeply influential works.
 
Our invariants grew out of Monroe--Sinha \cite{monroe-sinha}, who discovered that %, motivated by geometric constructions in rational homotopy theory,
ideas from classical linking theory generalize to letters in words $w\in F_n$ -- the free group on $n$ generators. 
For example, in
the commutator $(xyx^{-1}y^{-1})$ the $x$'s link with the $y$'s once.
%-- see Example \ref{ex:intro-braiding} for more. 
They introduced combinatorial linking of letters, and show that, when defined, their linking invariants coincide with coefficients in Magnus' classical power series expansion of $F_n$ \cite{Magnus_1935}.

Below, we give an ``open-string" generalization of these letter-linking invariants, which we therefore call \emph{letter-braiding}\footnote{Recall that links are closures of braids -- the former consisting of closed strings while the latter having open strands. Similarly, Monroe--Sinha's letter-linking is defined when all `strings of letters' are closed, while we consider strands with ends.}. These vastly extend the scope of  Mornoe--Sinha's letter-linking, as well as Magnus' expansion. They also give a simplicial analog to Chen's celebrated iterated integrals \cite{chen1971pi1} -- one that is combinatorial and can detect torsion phenomena.

\subsection{Main Results}

%We use constructions that are standard in algebraic topology.  
Let $\Gamma$ be a 
group and  
let $X$ be any connected simplicial set with fundamental group $\Gamma$. 
In applications, choices of $X$ can range 
 from small combinatorial models
 to functorial models and the singular simplicial set of 
a topological space. Let $A$ be a PID.

The cohomology ring of $X$ with coefficients in $A$ is computed by the associative algebra of simplicial cochains $C^*_{\cup}(X;A)$. The $A$-valued letter-braiding invariants of $\Gamma$ are parametrized by the $0$-th (associative) bar cohomology of $C^*_{\cup}(X;A)$. This is denoted by
\[
\H^0_{Bar}(X;A) := \H^0\left( Bar(C^*_{\cup}(X;A)) \right)
\]
and is a subquotient of the tensor (co)algebra $\bigoplus_k C^1_{\cup}(X;A)^{\otimes k}$, see \S\ref{sec:bar construction} for definitions.

Every $T\in \H^0_{Bar}(X;A)$ defines an $A$-valued invariant $\ell_T$ of elements $w\in \Gamma$, as follows. Representing $w$ by a map $\tilde{w}\maps S^1\to X$, it defines a pullback on bar constructions
\[
\tilde{w}^* \maps \H^0_{Bar}(X;A) \to \H^0_{Bar}(S^1;A).
\]
We exhibit a canonical invariant $\ell\maps \H^0_{Bar}(S^1;A) \to A$, through which one obtains the \emph{letter braiding number} by composition $\ell_T(w):= \ell(\tilde{w}^*(T)) \in A$. A combinatorial formula for $\ell_T(w)$ on free groups is this: for homomorphisms $H_1,\ldots,H_k\maps\Gamma\to A$ and generators $s_1,\ldots s_n \in \Gamma$,
\[
T=[H_1|\ldots|H_k] \text{ and } w=s_1^{\epsilon_1} s_2^{\epsilon_2} \ldots s_n^{\epsilon_n} \quad\rightsquigarrow \quad\ell_T(w) = \sum_{i_1<_w \ldots <_w i_k} H_1(s_{i_1}^{\epsilon_{i_1}})H_2(s_{i_2}^{\epsilon_{i_2}})\ldots H_k(s_{i_k}^{\epsilon_{i_k}}),
\]
where order relation $i<_w j$ between natural numbers is the standard $<$ relation when $\epsilon_i=1$ and $\leq$ when $\epsilon_i = -1$.

%We give an explicit weight-reduction algorithm 
%in Proposition \ref{prop:weight-reduction}, inspired by Monroe--Sinha's linking of letters \cite{monroe-sinha},  which defines an isomorphism $L\maps \H^0_{Bar}(S^1;A) \xrightarrow{\sim} A[t].$
%Our letter-braiding invariants are the composite $L \circ \tilde{w}^*\in A[t]$.  

%Thus $T\in \H^0_{Bar}(X;A)$ determines a polynomial invariant $L_{T}(w)= L(\tilde{w}^*(T))\in A[t]$. The linear coefficient of the polynomial $L_T(w)$ is the \emph{letter-braiding number}  $\ell_T(w)$.

This collection of functions $\ell_T(\bullet)$ for $T\in \H^0_{Bar}(X;A)$ turns out to depend only on the group $\Gamma = \pi_1(X,*)$, as explained in \S\ref{sec:general groups}, and is thus denoted $\H^0_{Bar}(\Gamma;A)$. Extending $A$-linearly,

\begin{definition}
    The letter-braiding numbers define an $A$-bilinear pairing with the group ring $A[\Gamma]$, called the \emph{letter-braiding pairing},
\begin{equation}\label{eq:intro-pairing}
    \langle - , \bullet \rangle \maps \H^0_{Bar}(\Gamma;A) \otimes A[\Gamma] \to A \quad\text{given on group elements by} \quad \langle T, w\rangle := \ell_T(w).
\end{equation}
\end{definition}
When $\Gamma$ is presented by generators and relations, computing both $\H^0_{Bar}(\Gamma;A)$ and its pairing on words is given by explicit, simple algorithms. When $\Gamma$ is the fundamental group of some manifold or variety, braiding invariants have compelling geometric interpretations.

The letter-braiding pairing is compatible with the natural structures on $A[\Gamma]$ and $\H^0_{Bar}(\Gamma;A)$. Recall that $A[\Gamma]$ is a ring under the convolution product, and it is filtered by powers of the augmentation ideal $I= \ker(\epsilon\maps A[\Gamma]\to A)$. This filtered algebra structure is functorial with respect to group homomorphisms. 

Our main result, proved in \S\ref{sec:group rings} and \S\ref{sec:completeness}, is the following.

\begin{theorem}[Fundamental theorem of letter-braiding]\label{thm:intro-pairing}
    Let $A$ be a principal ideal domain. For any group $\Gamma$,  its
    bar cohomology $\H^0_{Bar}(\Gamma;A)$ is a filtered coalgebra, functorial in $\Gamma$.  Under the letter-braiding pairing from \eqref{eq:intro-pairing}, this structure is dual to the filtered algebra structure on $A[\Gamma]$, in the following sense.
    \begin{enumerate}
        \item Functoriality: for any group homomorphism $h\maps \Gamma_1\to \Gamma_2$,
        \[
        \langle h^*(T), w\rangle = \langle T, h_*(w)\rangle.
        \]
        \item Product-coproduct duality\footnote{Thinking of letter-braiding numbers as coefficients in a power series expansion, product--coproduct duality is equivalent to the statement that this expansion is \emph{multiplicative}. E.g. the Kontsevich integral on pure braids is a multiplicative expansion, while Mostovoy--Willerton's $\Z$-valued expansion is not. The $\Z$-valued braid invariant produced by our procedure would be multiplicative, but it would take some work to find the ring over which it is defined -- see Example \ref{ex:braids} for coefficients of degree $\leq 3$.}: setting $\langle T_1\otimes T_2, x\otimes y\rangle := \langle T_1,x\rangle\langle T_2,y\rangle$, the deconcatenation coproduct $\Delta$ on $\H^0_{Bar}(\Gamma;A)$ is dual to the convolution product on $A[\Gamma]$,
        \[
        \langle \Delta T, x\otimes y\rangle = \langle T, x\cdot y\rangle.
        \]
        \item \label{case:intro-filtrations} Filtrations: let $I^n\leq A[\Gamma]$ denote the $n$-th power of the augmentation ideal, then
        \[
        (I^{n})^{\perp} = \H^0_{Bar}(\Gamma;A)_{<n},
        \]
        the sub-coalgebra of weight $< n$ tensors. In particular, this coalgebra vanishes on the $n$-th term of the lower central series of $\Gamma$.

        \item Extending Magnus: when $\Gamma$ is a free group, the letter-braiding pairing coincides with coefficients in Magnus' power series expansion: $F_n \to \Z\langle X_1,\ldots,X_n\rangle^{\wedge}$ given by
        $$s_i\mapsto 1+X_i  \quad \text{ and } \quad s_i^{-1} \mapsto 1-X_i+X_i^2-X_i^3+\ldots $$
        and extended to words as a group homomorphism, see Example \ref{ex:intro-free group}.
    \end{enumerate}
\end{theorem}
Note that functoriality and the evaluation formula on free groups together determine the letter braiding pairing uniquely. Moreover, letter braiding invariants are maximal given the constraints enforced by the last theorem, as explained next.

\begin{theorem}[Completeness/Universal finite-type invariant]
    \label{thm:intro-complete} When $\Gamma$ is finitely generated, the letter-braiding pairing defines an isomorphism of filtered coalgebras
        \[
        \H^0_{Bar}(\Gamma;A) \xrightarrow{\sim} \varinjlim \Hom_{A}( A[\Gamma]/I^n , A ).
        \]
        This map is defined and injective even without assuming finite generation (see Remark \ref{rmk:infinite generation}).
\end{theorem}
For context on Theorem \ref{thm:intro-complete}, recall that in braid theory, an invariant of the braid group $B_n$ has \emph{finite type} (also called a \emph{Vassiliev invariant}) whenever it vanishes on some power of $I\leq \Z[B_n]$ (see \cite{Mostovoy-Willerton:magnus_for_braids} for this perspective). Extending this to general groups, Property \eqref{case:intro-filtrations} of the pairing says that every $T\in \H^0_{Bar}(\Gamma;A)$ defines a finite type invariant $\ell_T(\bullet)\maps \Gamma \to A$.
Furthermore, a finite type invariant $\Psi\maps \Gamma\to R$ is \emph{universal} if every finite type invariant of $\Gamma$ factors uniquely through $\Psi$.  Theorem \ref{thm:intro-complete} says that every $A$-valued finite type invariant is of the form $\ell_T( \bullet)$ for a unique $T\in \H^0_{Bar}(\Gamma;A)$ for $\Gamma$ finitely generated, exhibiting universality. 

For infinitely generated groups, the following generalization holds. First, recall that a counital coalgebra $(C,\Delta,\eta)$ with an element $1\in C$ such that $\eta(1)=1$ is said to be \emph{conilpotent} every fixed element $x\in \ker(\eta)\subseteq C$ is annihilated by some iterate of the reduced coproduct
\[
\overline\Delta(x) := \Delta(x)-(x\otimes 1+1\otimes x) \text{ where } \overline\Delta^{k+1} := (\overline\Delta\otimes \Id^{\otimes k})\circ \overline\Delta^k.
\]
The data of $1\in C$ is called a \emph{coaugmentation} of $C$, and a homomorphism of coaugmented coalgebras $(C,1_C)\to (C',1_{C'})$ is a homomorphism that furthermore takes $1_C$ to $1_{C'}$.
\begin{definition}\label{def:universal pairing}
    Given a coaugmented coalgebra $(C,1)$ over $A$, say that a pairing
    \[
    \langle -,\bullet\rangle: C\otimes A[\Gamma]\to A
    \]
    is \emph{multiplicative} if $\langle 1,x\rangle \equiv 1$ and $\langle \Delta T,x\otimes y\rangle = \langle T,x\cdot y\rangle$ for all $x,y\in \Gamma$, where we define $\langle T_1\otimes T_2,x\otimes y\rangle := \langle T_1,x\rangle \cdot \langle T_2,y\rangle$.
    
    Say that a multiplicative pairing $\langle-,\bullet\rangle_C$ of $(C,1_C)$ with $\Gamma$ is \emph{universally conilpotent} if $(C,1_C)$ is conilpotent and for every multiplicative pairing $Z\otimes A[\Gamma]\to A$ with $(Z,1_Z)$ conilpotent there exists a unique homomorphisms of augmented coalgebras $(Z,1_Z)\to (C,1_C)$ pulling back the pairing on $C$ to the one on $Z$.
\end{definition}
By standard categorical arguments, if $(C,1_C)$ and $(C',1_{C'})$ both support universally conilpotent pairings with $\Gamma$ then they are canonically isomorphic as augmented coalgebras. When this happens we call $(C,1_C)$ the \emph{conilpotent dual} of $\Gamma$. For infinitely generated groups $\Gamma$, the universal finite-type invariant takes the following form.
\begin{theorem}[Universal finite-type invariant of arbitrary groups]\label{intro-thm:infinite type universal pairing}
    Let $\Gamma$ be an arbitrary possibly infinitely-generated group, and let $A$ be any principal ideal domain. Then the letter-braiding pairing
    \[
    \langle -,\bullet\rangle:\H^0_{Bar}(\Gamma;A)\otimes A[\Gamma] \to A
    \]
    is a universally conilpotent multiplicative pairing. In particular, it induces an injection
    \[
    \H^0_{Bar}(\Gamma;A)\into \Hom_{A}( A[\Gamma] , A )
    \]
    whose image is the maximal $A$-submodule containing the functional $1:x\mapsto 1$ for all $x\in \Gamma$ and on which the operation \[
    m^*:\Hom_{A}( A[\Gamma] , A )\to \Hom_{A}( A[\Gamma]\otimes A[\Gamma] , A ),
    \]
    sending $f:A[\Gamma]\to A$ to the functional $ (f\circ m):(x\otimes y\mapsto f(x\cdot y))$ factors through a conilpotent coproduct structure.
\end{theorem}

Nearly a century ago, Magnus discovered a universal finite-type invariant for free groups \cite{Magnus_1935}, valued in a noncommutative power series ring. Later, Fox developed his `free differential calculus' to compute the coefficients of this power series \cite{Fox1}. These two works lead to an explosion of activity and countless applications spanning topology and group theory.

A generalization of Magnus' expansion to other groups has until now been elusive, and was only known in a small number of special classes of `sufficiently free' groups: right-angled Artin groups (RAAGs) are partially free and partially commutative, and a Magnus expansion was constructed by Duchamp--Krob \cite{Duchamps-Krob-RAAGS}; while pure braid groups are semi-direct products of free groups, and Magnus expansions over $\Z$ were found by Mostovoy--Willerton \cite{Mostovoy-Willerton:magnus_for_braids} and Papadima \cite{papadima-braidinvariants-overZ}, but their expansions are not \emph{multiplicative} in that they do not respect the product structure on braids. The most general construction of invariants for general presented groups might still be that of Chen--Fox--Lyndon \cite{chen-fox-lyndon}, who give a word-rewriting algorithm that computes the lower central series of groups from a presentation. %They build on Fox calculus as a dual basis for the lower central series Lie algebra of free groups.  
Their algorithm produces functionals on the lower central series subquotients of a group, but
it is inductive with equations that must be solved at each stage.

The construction in this paper explains why expansions for arbitrary groups have been hard to find, by exhibiting the obstruction: as discussed in \S\ref{sec:invariants from Massey}, expansion coefficients of degree $k$ require knowledge of all Massey products $\mu_n\maps \H^1(\Gamma;A)^{\otimes n}\to \H^2(\Gamma;A)$ of size $n\leq k$. These products are often unknown and hard to compute. For RAAGs and braid groups all higher Massey products vanish as a consequence of formality over $\Q$, so their expansions are easier to determine completely. 

We discuss and compare with older work throughout, but let us highlight some advantages of our approach over previous ones.
\begin{itemize}
    \item General groups: every group $\Gamma$ can be probed by a collection of letter-braiding invariants, which are well-defined even on non-free or infinitely generated groups.
    \item Globally defined: all letter-braiding invariants are well-defined on the entire group, not just high-order commutators or elements of dimension
    subgroups.
    \item Cohomological interpretation: letter-braiding invariants are built from combinations of classes in $\H^1(X;A)$, for any space with $\pi_1(X,*)\cong \Gamma$, and thus carry topological meaning alongside their combinatorial manifestation.
    \item Coefficients in PIDs: invariants can be made to take values in any PID, and thus detect torsion elements.
\end{itemize}

\tableofcontents

\subsection{First examples, motivation and background}

Letter-braiding invariants enjoy an explicit and effective combinatorial formulation, but they come from intersection theory in the geometry of loop spaces. If $M$ is a manifold, classes $\alpha\in \H^1(M;\Q)$ are represented by counting transverse intersections of loops with real hypersurfaces in $M$. Given a sequence of such hypersurfaces, one may track the order in which intersections
occur and count series of intersections in prescribed orders -- this is morally what letter-braiding invariants are measuring (mirroring Hain's construction using iterated integrals \cite{hain1985linking}). Such ordered intersection counts are generally not homotopy invariant, e.g. when the hypersurfaces intersect, but one can find linear combinations of counts that are invariant through analysis of such intersections, as reflected by cup products and Massey products. 

The algebraic framework for this analysis is governed by the classical bar construction, which has been used to model cochains on (a completion of) loop spaces for over sixty years \cite{adams1956cobar, EilenbergMoore, dwyer1975exotic}.
From this perspective, one is studying the fundamental group $\pi_1(X,x)$ through the zeroth cohomology of the loop space $\Omega_x X$. 

In this introduction, we will recall several classical approaches to studying fundamental groups using the bar construction, but few of those give rise to invariants of words in groups that are explicit enough to, say, be calculated on a computer. Letter-braiding invariants achieve this: they give concrete and fast algorithms for obtaining $\H^0_{Bar}(\Gamma;A)$ and its letter-braiding pairing with $\Gamma$, as \href{https://github.com/benwalt/coLie/}{implemented in Sage} by Aydin Ozbek, Ben Walter  and others\footnote{Since the first draft of this paper, Ozbek has tragically passed away. His original code \href{https://github.com/AydinOzbek/Hopf_Invariants}{is still available}.}.

\begin{example}\label{ex:intro-free group}
    The free group $F_n = \langle s_1,\ldots,s_n\rangle$ is the fundamental group of a wedge
\[
W_n = \bigvee_{i=1}^n S^1
\]
whose $1$-cochains are $C^1(W_n;A) = A^n$ with trivial cup product. Its letter-braiding invariants are parametrized by the tensor (co)algebra
\[
T(A^n) = A\langle S_1,\ldots,S_n\rangle
\]
of noncommutative polynomials in the standard basis elements $S_1,\ldots,S_n\in A^n$. 

Fox calculus \cite{Fox1} supplies operators on the group algebra $A[F_n]$ called `derivatives' with which every group element is expanded as a Taylor series. Theorem \ref{thm:magnus} shows that for every word $w\in F_n$, the letter-braiding number associated with the monomial $(S_{i_1}S_{i_2}\ldots S_{i_k})$ is \[
\ell_{(S_{i_1}S_{i_2}\cdots S_{i_k})}(w) = \epsilon\left( \frac{\partial}{\partial s_{i_1}}\cdots \frac{\partial}{\partial s_{i_k}}w \right) ,
\]
where $\frac{\partial}{\partial s_{i}}$ is the $i$-th Fox derivative from \cite{Fox1} and $\epsilon\maps A[F_n] \to A$ is the standard augmentation map. Equivalently, this letter-braiding number coincides with the coefficient of $t_{i_1}\ldots t_{i_k}$ in the classical Magnus expansion for $F_n$ (see \cite{Magnus_1935}), which is the same as Fox's Taylor coefficient.
In this way, letter-braiding invariants are a posteriori seen to be functorial extensions of the classical Magnus expansion to arbitrary groups and coefficients in any PID. At the same time, they are substantially easier to calculate than Fox' or Magnus' invariants.

In \S\ref{sec:weight reduction} we give a computational approache to letter-braiding on free groups, inspired by Monroe--Sinha's weight reduction procedure \cite{monroe-sinha}. This leads to fast and flexible calculations, with one implementation resulting in the following closed form expression.

\begin{prop}\label{intro-prop:interated sum}
    On the free group $F_n=\langle s_1,\ldots,s_n\rangle$, the letter braiding numbers $\ell_{(S_{i_1}S_{i_2}\cdots S_{i_k})}(\bullet)$ count the number of $<_w$-increasing subsequences of letters $(\ldots s_{i_1}^{\pm} \ldots s_{i_2}^{\pm} \ldots (\cdots)\ldots s_k^{\pm}\ldots )$ in a word, counted with signs:
    \[
    \ell_{(S_{i_1}S_{i_2}\cdots S_{i_k})}(s_{j_1}^{\epsilon_1}s_{j_2}^{\epsilon_2} \cdots s_{j_n}^{\epsilon_n}) = \sum_{1\leq r_1{<}_w r_2{<}_w \ldots {<}_w r_k\leq n} \epsilon_{r_1}\epsilon_{r_2}\cdots \epsilon_{r_k} S_{i_1}(s_{j_{r_1}})S_{i_2}(s_{j_{r_2}})\cdots S_{i_k}(s_{j_{r_k}}),
    \]
    where $S_{i}(s_j)=\delta_{ij}$ is the Kronecker delta, and $r<_w r'$ is the order relation $r<r'$ when $\epsilon_i=1$ and it is $r\leq r'$ when $\epsilon_i=-1$. 
\end{prop}
For example, consider $F_2 = \langle x,y\rangle$ and let $X$ and $Y$ be the indicator functions of the respective basis element. Then last the summation formula identifies the letter braiding number of $(XXYX)$ with a signed count of ``$x$'s before $x$'s before $y$'s before $x$'s". On the word $w = [xy,x^{-2}] = xyx^{-1}x^{-1}y^{-1}x$ this calculation may be represented figuratively as
\[
\begin{aligned}
\begin{tikzpicture}

\node at (0,0) {$^{\phantom{-1}}x^{\phantom{-1}}$};
\node at (1,0) {$^{\phantom{-1}}y^{\phantom{-1}}$};
\node at (2,0) {$x^{{-1}}$};
\node at (3,0) {$x^{{-1}}$};
\node at (4,0) {$y^{{-1}}$};
\node at (5,0) {$^{\phantom{-1}}x^{\phantom{-1}}$};

\draw[->] (-.4,.25) -- (.4,.25);
\draw[->] (.6,.25) -- (1.4,.25);
\draw[->] (2.4,.25) -- (1.6,.25);
\draw[<-] (2.6,.25) -- (3.4,.25);
\draw[<-] (3.6,.25) -- (4.4,.25);
\draw[->] (4.6,.25) -- (5.4,.25);

\node (v1) at (0,0.5) {$+1$};
\node (v2) at (2,0.5) {$-1$};
\node (v3) at (3,0.5) {$-1$};
\node (v4) at (5,0.5) {$+1$};

\draw[->] (.5,.25) -- (.5,1.05) -- (6,1.05);
\draw[->, densely dotted] (1.5,.25) -- (1.5,.95) -- (6,.95);
\draw[->, densely dotted] (2.5,.25) -- (2.5,.85) -- (6,.85);
\draw[->] (5.5,.25) -- (5.5,.75) -- (6,.75);

\node at (0,.25) {$\bullet$};
\node at (2,.25) {$\circ$};
\node at (3,.25) {$\circ$};
\node at (5,.25) {$\bullet$};

\node  at (2,1.05) {$\circ$};
\node  at (2,.95) {$\bullet$};

\node  at (3,1.05) {$\circ$};
\node  at (3,.95) {$\bullet$};
\node  at (3,.85) {$\bullet$};

\node  at (5,1.05) {$\bullet$};
\node  at (5,.95) {$\circ$};
\node  at (5,.85) {$\circ$};

\node (v12) at (3,1.35) {$+1$};
\node (v22) at (5,1.35) {$-1$};

\draw[->] (2.5,1.15) -- (2.5,1.75) -- (6,1.75) ;
\draw[->, densely dotted] (5.5,1.15) -- (5.5,1.65) -- (6,1.65) ;

\node at (4,1.75) {$\circ$};
\node at (4,2.05) {$-1$};

\draw[->, densely dotted] (3.5,1.8) -- (3.5,2.35) -- (6,2.35) ;

\node at (5,2.35) {$\circ$};
\node at (5,2.65) {$-1$};
\draw  (4.7,2.9) rectangle (5.3,2.45);
%\draw[->, densely dotted] (5.5,2.4) -- (5.5,2.9) -- (6,2.9) ;

\node at (-1,0.5) {X};
\node at (-1,1.4) {X};
\node at (-1,2.1) {Y};
\node at (-1,2.7) {X};

\end{tikzpicture}
\end{aligned}
\implies \ell_{XXYX}([xy,x^{-2}]) = -1
\]
The bottom row represents the word itself, with segments labeled by letters, oriented according to the sign. Every dot records a letter that contributes to the sum, and it has an arrow coming up from the tip of the corresponding segment to indicate the positions that satisfy order relation $<_w$ relative to the dot's position. Solid dots and arrows count as $(+1)$, while the non-solid ones count as $(-1)$. The figure shows that the Magnus expnasion of $w$ includes the term $(-1)XXYX$, detecting e.g. that $w$ is not a $4$-fold nested commutator.
\end{example}

\begin{example}[Torsion commutators]
Let $\Gamma=\mathscr{H}_3(\F_2)$ be the Heisenberg group of matrices
\[
\begin{pmatrix}
    1&X&Z\\
    0&1&Y\\
    0&0&1
\end{pmatrix}\in \operatorname{Gl}_3(\F_2).
\]
The functions $X,Y,Z\maps \mathscr{H}_3(\F_2)\to \F_2$ determine cochains in $C^1(\mathscr{H}_3(\F_2);\F_2)$, and Example \ref{ex:heisenberg} shows $XY+Z\in \H^0_{Bar}(\mathscr{H}_3(\F_2);\F_2)$ defines an invariant on words in $\mathscr{H}_3(\F_2)$, which pairs nontrivially with the commutator $\begin{psmallmatrix}
1&0&1\\0&1&0\\0&0&1    
\end{psmallmatrix}
=\left[\begin{psmallmatrix}
1&1&0\\0&1&0\\0&0&1    
\end{psmallmatrix},\begin{psmallmatrix}
1&0&0\\0&1&1\\0&0&1    
\end{psmallmatrix}\right]$.

With this, letter-braiding has a similar flavor to Chen's theory of iterated integrals \cite{chen1971pi1}, but can be defined over $\mathbb{F}_p$ and detect torsion elements.
\end{example}

\subsection{Hopf invariants} In a direct predecessor of this work, Sinha--Walter \cite{sinha-walter} use the bar construction, and a variant for commutative algebras called the \emph{Harrison complex}, to define invariants of maps $\alpha\maps S^n \to X$ with $n\geq 2$ using higher linking numbers in spheres $S^n$. They call these \emph{Hopf invariants} and show that they distinguish rational homotopy classes of maps into simply connected spaces. The present paper grew out of a joint effort to define Hopf invariants of maps $S^1\to X$ using `linking' of $0$-manifolds in $S^1$, together with with A. Ozbek, D. Sinha, and B. Walter \cite{joint}. The theory we develop there is complementary to letter-braiding, starting with a commutative model for cochains in characteristic zero, and some of its benefits are the following.
\begin{itemize}
    \item Hopf invariants of maps $S^1\to X$ detect and obstruct the Malcev Lie algebra of the group $\pi_1(X,*)$. The Harrison complex is much smaller than the bar complex, and it is naturally a Lie coalgebra. As such, it removes redundancy from functionals on the group ring, and computes a maximal Lie coalgebraic dual of the Malcev Lie algebra of any $\Gamma$.
    \item The theory is naturally geometric: de Rham forms supported near submanifolds (a.k.a. Thom forms) realize integration and wedge products as intersection of submanifolds. As such, letter-braiding theory is restated in terms of the meeting times of a smooth loop with with various hypersurfaces, as described above.
    \item Hopf invariants specialize to Chen's \cite{chen1971pi1} theory of iterated integrals, on the one hand, while on the other hand having a combinatorial formulation relative to group presentations.
\end{itemize}
Hopf invariants differ subtly from letter-braiding invariants studied here, which are based on simplicial cochains with their noncommutative cup product. In \cite{joint} we prove that the former are obtained from the latter as the invariants of $\log(w) = \sum_{n\geq 1}\frac{(-1)^{n+1}}{n}(w-1)^n$, as expected from functionals on the Malcev Lie algebra.

Topological applications of this theory include the following (see \cite{joint} for terminology).

%\begin{prop}
%    The complexified Malcev Lie algebra of a complex hyperplane arrangement complement is combinatorially determined by the associated matroid. In fact, incidences of triples of hyperplanes are sufficient.
%\end{prop}
%The is surprising in light of Rybnikov's example \cite{rybnikov2011fundamental} of two hyperplane arrangements with the same matroid but non-isomorphic nipotent quotients of fundamental groups $\pi/\gamma_4 \pi$.
\begin{prop}
    Milnor invariants of links are computed from the intersections of a link component with any collection of Seifert surfaces for the link, along with further secondary surfaces: they are the letter-braiding numbers of the intersection pattern.
    
    Further invariants are well-defined beyond the weights which Milnor initially considered, and for any element in the link group.
\end{prop}
Here, Seifert surfaces are allowed to intersect, at the cost of adding secondary surfaces. In ongoing work in progress, we use this additional freedom to give algorithms for calculating Milnor invariants from a planar projection of the link.

\begin{comment}
    
\begin{example}
    Let $PB_n$ denote the pure braid group on $n$ strands, and denote $\omega_{ij}\maps PB_n \to \Z$
    \[
    \omega_{ij}(\sigma) = (\text{winding number of strand $i$ with strand $j$ in braid $\sigma$}).
    \]
    Then tensors of the form
    $T_{ijk}=(\omega_{ij}|\omega_{jk}+\omega_{jk}|\omega_{ik}+\omega_{ik}|\omega_{ij}) $  for $1\leq i<j<k\leq n$ define additive invariants on commutator braids $\ell_\bullet(T_{ijk})\maps [PB_n,PB_n]\to A$, when on basic commutators $[\sigma_1,\sigma_2]$ one takes
        \[
        \ell_{[\sigma_1,\sigma_2]}(\omega_1|\omega_2) = \omega_1(\sigma_1)\omega_2(\sigma_2) - \omega_1(\sigma_2)\omega_2(\sigma_1).
        \]
        See Example \ref{ex:braids} for details and for further invariants on triple commutators.
\end{example}

\end{comment}

\subsection{First Applications}

The invariants $T\in\H^0_{Bar}(\Gamma;A)$ are determined by any space or simplicial set $X$ with $\pi_1(X,*)\cong \Gamma$, thus different choices of model give invariants of distinct flavors, blending Algebra, Combinatorics, Geometry and Topology.

\subsubsection{Combinatorial group theory}
A presentation $\Gamma = \langle S|R \rangle$ determines a $2$-complex $X_{S|R}$ whose $1$-cells are labeled by $S$ and $2$-cells come from relations $R$. Its letter-braiding invariants are therefore presented in terms of the generators and relation,
as follows.

For any set $S$, let $F_S$ denote the free group generated by $S$. Its letter-braiding invariants are computed from the wedge $X_S = \vee_{s\in S} S^1$ and are simply the tensor (co)algebra
\[
T(A^S) = \bigoplus_{k\geq 0} (A^S)^{\otimes k}.
\]
As mentioned, letter-braiding numbers of words $w\in F_S$ coincide with the output of Fox' classical derivatives on free groups (see Theorem \ref{thm:magnus}). We thus think of them as differential operators: $f\in A^S \longleftrightarrow \sum_{s\in S} f(s)\frac{\partial}{\partial s}$.

Now, the inclusion $X_S \into X_{S|R}$ identifies $\H^0_{Bar}(\Gamma;A)$ as a sub-coalgebra of $T(A^S)$, consisting of possibly nonhomogeneous differential operators, so that the letter-braiding numbers of $w\in \langle S|R\rangle$ are given by the Magnus expansion coefficients of any lift $\tilde{w}\in F_S$.

But to get well-defined invariants of words in $\langle S|R\rangle$, we may only consider differential operators that happen to kill all relations $\langle R \rangle \leq F_S$. Examples \ref{ex:relation1},\ref{ex:relation2} discuss basic cases of relations constraining the collection of invariants. In fact, we show that vanishing on relations, suitably understood,
is the only constraint.

\begin{theorem}
    For $\Gamma = \langle S|R\rangle$, the quotient map $p\maps F_S \onto \Gamma$ induces an inclusion
    \[
    p^*\maps \H^0_{Bar}(\Gamma;A)\into T(A^S)
    \]
    whose image is the maximal sub-coalgebra of differential operators that vanish on all $r\in R$.
\end{theorem}

In other words, a differential operator on $\Gamma$ must annihilate
all relations $r\in R$, but this by itself is insufficient -- all of its iterated coproducts must annihilate the relations as well. Over the real numbers and restricting to manifolds, this was proved by Chen using iterated integrals on the bar construction of the de Rham complex \cite{Chen_iterated_integrals77}; the result over torsion PIDs appears to be new.

Theorem \ref{thm:lifting criterion} gives a precise and useful generalization of this theorem. With it, one can compute the letter-braiding invariants of $\Gamma$ as follows: start with the linear differential operators that vanish on $r\in R$; then find among their degree $2$ compositions those linear combinations that continue to vanish on $r\in R$; and so on.

\subsubsection{Geometric topology}
Suppose $\Gamma$ arises as the fundamental group of a manifold $M$ rather than from a presentation. Then its letter-braiding invariants are expressible in geometric terms. Examples include surface groups, Artin groups, knot and link groups and more general 3-manifold groups, and (finite-index subgroups of) mapping-class groups.

 Lemma \ref{lem:Bar spectral sequence} recounts that $\H^0_{Bar}(M;A)$ is determined by sums of tensors of cohomology classes $\H^1(M;A)^{\otimes k}$, all of whose cup and higher Massey products vanish -- see Examples \ref{ex:relation4}, \ref{ex:relation5}. By Poincar\'{e} duality, classes $\alpha\in \H^1(M;A)$ are represented by intersection with a dual hypersurface $H_\alpha \subset M$. A word $w\in \Gamma$ is represented by some loop in $M$, which intersects those hypersurfaces at various times. The letter-braiding numbers of $w$ are then determined by the interleaving of these intersections -- see \S\ref{sec:weight reduction}. With $A=\R$ coefficients, this description was given in \cite{hain1985linking}.

However, different choices of representative loops for $w$ might meet the hypersurfaces $H_\alpha$ in a different order. The condition of vanishing cup and Massey products turns out to be equivalent to the resulting letter-braiding number being a homotopy invariant. Since cup products are dual to intersection, relations in $\pi_1(M,*)$ constrain the collections of hypersurfaces that may define invariants, and therefore force nontrivial intersection patterns of hypersurfaces. 

For example, the fundamental group of the Borromean ring complement is 
\[
\Gamma_3 = \langle x,y,z\mid[x,[y^{-1},z]],[y,[z^{-1},x]],[z,[x^{-1},y]]\rangle.
\]
whose double-commutator relations are detected by the invariant $XYZ$ of $F_3=\langle x,y,z\rangle$, reflecting a triple Massey product $\langle X,Y,Z\rangle\neq 0$. Geometrically, $X,Y$ and $Z$ are dual to Seifert surfaces of the three link components. Thus, the Massey product implies that no choice of Seifert surfaces can be made disjoint, as originally observed by Massey.

\subsubsection{Positive characteristic}
Taking coefficients rings $A$ of positive charactertistic, letter-braiding can detect torsion elements and commutators in finite groups. Recall that the lower exponent $p$-central series of a group $\Gamma$ is the filtration induced by $w\in F_k\Gamma\iff (w-1)\in I^k\leq \F_p[\Gamma]$. It follows that $\F_p$-valued letter-braiding invariants obstruct this $p$-central series. They are particularly informative when $\Gamma$ is a finite $p$-group, since in this case the lower $p$-central series terminates and, therefore, any two elements $w\neq w' \in \Gamma$ can be distinguished by letter-braiding invariants.

Conversely, the structure of $\F_p[\Gamma]$ provides information about Massey products in $\F_p$-cohomology.
\begin{example}
    For $\Gamma = C_p$, the cyclic group of order $p$, the group ring is $p$-dimensional, while $\H^1(C_p;\F_p) = \F_p\cdot T$ generated by the $\F_p$-group homomorphism $T(i) = i$. The braiding invariants for $\Gamma$ are thus spanned by $T^k+(\text{lower weight correction terms})$ for $0\leq k\leq p-1$, see Example \ref{ex:finite cyclic} for explicit cochain representatives for which
    \[
    \langle T^k + \text{ (lower weights) }, \, i \rangle = \frac{i^k}{k!} \mod{p}.
    \]
    That these already span the dual of $\F_p[C_p]$ means that $T^p$ does not extend to an invariant.
    \begin{cor}
        For any space $X$ with $\pi_1(X,*)\cong C_p$, the generator $T\in \H^1(X;\F_p)$ has vanishing cup and Massey powers up to order $p-1$, but its $p$-th Massey power
        \[
        \langle\underbrace{ T,\ldots,T }_{p \text{ times}}\rangle \neq 0.
        \]
    \end{cor}
    For $C_{p^2}$ the same $p$-fold Massey power is $0$, and $T^p$ detects the nontrivial element $p\in C_{p^2}$.
\end{example}
See \S\ref{sec:Massey in finite groups} for a generalization to all finite groups.

\subsubsection{Relative Massey products vs. relations}
The previous two subsections exhibit a dictionary between Magnus expansions of relations in $\Gamma$ and Massey products on $\H^1(\Gamma;A)$. That such a connection exists has been well-known, going back to Dwyer \cite{Dwyer} and Fenn--Sjerve \cite{other-Massey-paper}.

We complete this dictionary in \S\ref{sec:relative Massey} by defining \emph{relative Massey products}.  Given a pair of spaces $(Y,X)$ there are higher product operations, that in suggestive simple cases become
\[
\mu_n\maps \H^1(X;A)^{\otimes n} \to \H^2(Y,X;A).
\]
When $Y$ is obtained from $X$ by attaching `relation' $2$-discs along a set $R\subseteq \pi_1(X,*)$, the relative cohomology $\H^2(Y,X;A)\cong A^R$, recording functions on $R$, and Proposition \ref{prop:massey=braiding} shows that the relative Massey product $\mu_n$ computes the letter-braiding numbers of all $r\in R$. Our setup thus vastly extends the picture painted by Fenn-Sjerve \cite{other-Massey-paper} in a number of ways, listed at the end of \S\ref{sec:relative Massey}.

\subsubsection{Algebraic models for spaces}
%\subsubsection{Cochains determine group ring}
%To explain how letter-braiding invariants arise, let $\Gamma$ be a discrete group and let $X$ be any topological space or simplicial set with $\pi_1(X,*) \cong \Gamma$.
Mandell \cite{mandell} showed that when $X$ is a nilpotent space of finite type, its cochains $C^*(X;\Z)$ with their $E_\infty$-algebra structure determine the homotopy type of $X$, and in particular $\pi_1(X,*)$.
More constructively, building on Adams' classical theorem \cite{adams1956cobar}, Rivera--Zeinalian \cite[Corollary 9.2]{rivera-cobar} explains that the coBar construction on the singular chains $C_*^{Sing}(X)$ models the chains of the based loop space $\Omega_* X$ and thus its $\H_0$ is computes the group ring $\Z[\pi_1(X,*)]$. However, their proofs use rather abstract homotopy theory, and do not provide any means for measuring words $w\in \pi_1(X,*)$.

In contrast, our proofs are essentially elementary, with algebraic rather than homotopy theoretic flavor. Consequently, we get a simple  constructive proof of the following.
\begin{cor}[{Compare Mandell \cite{mandell}, and Dwyer \cite[Theorem 2.1]{dwyer1975exotic}} for the special case of nilpotent spaces]
Let $(X,*)$ be a pointed connected simplicial set (or topological space) with $\pi_1(X,*)$ finitely generated but perhaps infinitely related, and let $A$ be any PID. Then the filtered coalgebra of $A$-valued functionals of finite type $\varinjlim \Hom_{A}(A[\pi_1(X,*)]/I^n,A)$ is determined by the $2$-truncation of the augmented simplicial (or singular) cochain algebra $C^{*\leq 2}_\Delta(X;A)$, via letter-braiding.
\end{cor}
Compare this further to Rivera--Zeinalian \cite{Rivera-chains_and_pi1}, who prove that the singular chains of a connected space $X$ determine the entire fundamental group -- no nilpotency hypotheses necessary. However, one major obstacle in applying their theory here is that one cannot simply pass to a small chain model for calculations, as their construction is not quasi-isomorphism invariant (instead, they use a finer notion of $\Omega$-equivalence).

%*** Work of Baues, Brown, Berrick-Cohen, six-author UASS paper
%shows that UASS can be derived from lower central series
%of a group model for loop space. *** (Aug ideal versus LCS...?)

\subsection{Andreadakis--Johnson homomorphisms}
Work of Kawazumi \cite{kawazumi2005cohomological} on Johnson filtrations indicates that Magnus expansions define generalizations of the classical Johnson homomorphism of the Torelli group. Since our setup is dual to Magnus expansions, and is defined for general groups, it lets us define dual Johnson homomorphisms in general.

For any group $\Gamma$, the automorphism group $\Aut(\Gamma)$ is filtered by subgroups $J^A(k) \leq \Aut(\Gamma)$ respectively fixing $A[\Gamma]/I^{k+1}$ pointwise -- a generalized Johnson filtration. Via the letter-braiding pairing, one obtains a Johnson-like homomorphism on these groups.
\begin{theorem}\label{thm:intro-johnson}
    Denote $H = \H^1(\Gamma;A)$, and let $K_n \leq \H^{\otimes n}$ be the tensors on which the $n$-th Massey product is defined and vanishes. Then for $\varphi\in J^A(k)$, the endomorphism $(T - \varphi^*(T))$ of $\H^0_{Bar}(\Gamma;A)_{\leq k+1}$ 
    determines an $\Aut(\Gamma)$-equivariant `dual Johnson' homomorphism 
    \[J^A(k)\to \Hom_A(K_{k+1},H),\] 
    Furthermore, its kernel always contains $J^A(k+1)$, and is exactly this subgroup whenever $\Gamma$ is finitely generated and $A[\Gamma]/I^{k+2}$ is $A$-torsion free, e.g. when $A$ is a field.
\end{theorem}
With $A=\Z$ and $\Gamma = F_n$ a free group, $K_k = H^{\otimes k}$ and the integral Johnson filtration $J^{\Z}(k)\leq \Aut(F_n)$ is the classical Andreadakis filtration. In this case our dual Johnson homomorphism maps $J(k) \leq \Aut(F_n)$ to $(H^*)^{\otimes k+1}\otimes H$ (compare this with the Johnson homomorphism in \cite[Theorem 6.1]{kawazumi2005cohomological}). And when $\Gamma=\pi_1(\Sigma_g,*)$, the mapping class group $\operatorname{Mod}(\Sigma_{g,1})$ acts on $\pi_1(\Sigma_g,*)$, sending the classical Johnson filtration into the our mod $A$ Johnson filtration. We have not found tools for calculating the dual Johnson homomorphism, and leave it as an open problem to relate it with the classical ones.

Theorem \ref{thm:intro-johnson} is proved in \S\ref{sec:Johnson}. With it, we get a new proof of the following classical theorem of Burnside.
\begin{theorem}\label{thm:aut of p-groups}
    If $\Gamma$ is a finite $p$-group, then $J^{\F_p}(1)\leq \Aut(\Gamma)$ is a $p$-group. In particular, if $\varphi\in \Aut(\Gamma)$ does not have $p$-power order, then $\varphi$ acts nontrivially on the maximal elementary abelian quotient $\Gamma^{ab}\otimes \F_p$.
\end{theorem}

\subsection{Related work}
A number of constructions of a similar flavor to ours appeared in the near and far past. Most close in spirit are the Chen integrals \cite{chen1971pi1}, which attach $\R$- or $\C$-valued invariants to words $w\in \Gamma$ via integration: if $M$ is a smooth manifold with $\pi_1(M,*)\cong \Gamma$ and $\tilde{w}\maps [0,1] \to M$ is a loop representing the homotopy class of $w$, Chen took
\[
\langle (\alpha_1|\ldots|\alpha_n)\, , \, w\rangle = \int_{0\leq t_1\leq \ldots\leq t_n \leq 1} \tilde{w}^*\alpha_1(t_1)\wedge \ldots \wedge \tilde{w}^*\alpha_n(t_n)
\]
for any tensor of de-Rham $1$-forms in the bar construction $\H^0( Bar(\Omega^*_{dR}(M) )$. He further proved that this defines an isomorphism of Hopf algebras with $\varinjlim \Hom_{\R}(\R[\Gamma]/I^n , \R)$. This construction informed later integral invariants, such as the Kontsevich integral, but it has not been made so explicit as to produce combinatorial invariants of words in groups.

Our weight reduction algorithm (Proposition \ref{prop:weight-reduction}) gives a discrete analog of Chen's iterated integral, see Corollary \ref{cor:iterated sum}. With it, the resulting invariants are made combinatorial, while also being defined over general PIDs. In addition, weight reduction manifests substantial freedom to choose the order in which $1$-forms integrated, of which Chen's `time-ordered' formalism is but one choice: for most words, other orderings yield more efficient calculations.

Concurrently with Chen, Stallings \cite{Stallings_augmentation_quotients} studied quotients $I^n/I^{n+1}$ of the augmentation ideal in the group ring $\Z[\Gamma]$. He introduced a spectral sequence converging to $I^\bullet/I^{\bullet+1}$, which is conceptually dual to the spectral sequence we use in the proof of this paper's main technical theorem, Theorem \ref{thm:completeness}.

On the interaction between Massey products on $\H^1(\Gamma;R)$ and relations in $\Gamma$, work of Dwyer \cite{Dwyer} interprets `Massey systems' as homomorphisms from a group of upper triangular matrices $\Gamma\to U_n(R)$, and Fenn--Sjerve \cite{other-Massey-paper} relate the evaluation of Massey products on $\H^2(\Gamma;\Z)$ with Magnus coefficients of relations.
Dually, Garoufalidis--Levine \cite{Stavros_Levine_vanishing_cup} relate Massey products with lower central series quotients of a group. They show that a homomorphism $F_n\to \Gamma$ that induces an isomorphism on $H_1$, further induces an isomorphism modulo the $k$-th step of the lower central series and a surjection on the $k+1$-st quotient iff all Massey products of order $\leq k$ vanish. At the $k+1$-st step, the difference between the two groups is measured by a homological dual Massey product. Our setting pushes this interaction further, by considering braiding invariants past the degree at which the first Massey products appear.

Recently, Porter--Suciu \cite{porter-suciu:cup1,porter2023cup} introduced an integral version of rational homotopy theory for studying fundamental groups. Their construction makes explicit use the Steenrod's cup-1 operation, and requires a binomial structure on the coefficient ring. With this, they develop a theory of $1$-minimal models that give explicit access to fundamental groups, including their torsion, starting from a model for singular cochains over a binomial ring. Around the same time, Horel \cite{horel2022binomial} found that binomial rings are a good setup for integral homotopy theory, and exhibited a fully faithful functor between nilpotent spaces and cosimplicial binomial rings.
Our  theory is not trying to pick out $\pi_1$ from the group ring $A[\pi_1]$, since braiding invariants are already defined at the group ring level, and therefore we can avoid the cup-1 operations and binomial structures.

\subsection{Acknowledgements}
This work grew out of a closely related project, joint with Aydin Ozbek, Ben Walter and Dev Sinha \cite{joint}. I thank them for their insights, questions, feedback and support. Special thanks to Dev Sinha and Greg Friedman for their careful reading and detailed feedback on an earlier draft. Aydin Ozbek inspired and taught me much before his untimely passing -- he left a noticeable mark on this work. I also wish to thank Richard Hain for his comments, and for introducing me to lots of relevant references.

This paper was partially supported by the National Science Foundation under Grant No. DMS-1929284 while the author was in residence at the Institute for Computational and Experimental Research in Mathematics (ICERM) in Providence, RI, during the Combinatorial Algebraic Geometry Reunion program.

\section{(Co)algebraic preliminaries}

\label{sec:prelims}
Let $A$ be a PID and recall the category of graded $A$-modules with objects $B = \bigoplus_{n\in \Z} B_n$ where every $B_n$ is an $A$-submodule. The \emph{shift} of $B$ by $k$ is the graded modules $B[k]$ given by $
(B[k])_n := B_{n+k}.$
Morphisms $f\maps B\to B'$ are $A$-linear maps that preserve grading. A homomorphism $f\maps B\to B'[k]$ is said to be a morphism of degree $k$
%(when the degree is unspecified, it shall be taken to be $0$).

Also recall that the tensor product of two graded modules has graded components
\[
(B\otimes B')_n = \bigoplus_{i\in \Z} B_i\otimes_A B'_{n-i}.
\]
Throughout this paper we consider conilpotent associative coalgebras.

\begin{definition}
    A \emph{differential graded coalgebra} is a graded module $B$ equipped with a coassociative comultiplication
    \[
    \Delta\maps B\to B\otimes B, \text{ a counit } \eta\maps B \to A
    \]
     and a coderivation $d\maps B\to B$ of degree $+1$ obeying the Koszul sign convention (see \cite[\S1.5.3]{algebraic_operads}). That is,
     \[
     (\Delta\otimes \Id)\circ \Delta = (\Id\otimes \Delta)\circ \Delta, \quad (\Id\otimes \eta)\circ \Delta = (\eta\otimes \Id)\circ \Delta = \Id \quad\text{ and }
     \]
     $\quad \Delta\circ d = (d\otimes \Id + (-1)^{n} \Id \otimes d)\circ \Delta$
     when evaluated on $B_n$.
    
    A homomorphism of dg coalgebras is a morphism $h\maps B\to B'$ such that
    \[
    (h\otimes h)\circ \Delta_B = \Delta_{B'} \circ h, \quad \eta_{B'}\circ h = \eta_B,  \quad \text{and} \quad h \circ d_B = d_{B'}\circ h.
    \]
    Denote iterated coproducts by $\Delta^0 = \Id$ and
    \[
    \Delta^{n+1} = (\Delta \otimes \Id^{\otimes n})\circ \Delta^n.
    \]
\end{definition}
The ring $A$ itself has a tautological dg coalgebra structure, with $\Delta\maps A\to A\otimes_A A$ the canonical isomorphism, $\eta\maps A\to A$ the identity and $d=0$. 

\begin{definition}
    An \emph{augmentation} of a dg coalgebra $B$ is a coalgebra homomorphism $\epsilon\maps A\to B$. Denote $1\in B$ for the element $\epsilon(1)$. This splits $B$ as $\ol{B}\oplus A\cdot 1$, and we call $\ol{B} = \ker(\eta)$ the \emph{reduced} coalgebra of $B$. The \emph{reduced comultiplication} $\ol{\Delta}\maps \ol{B}\to \ol{B}$ is
    \[
    \ol{\Delta}(b) = \Delta(b) - 1\otimes b - b\otimes 1.
    \]

Define
    \[
    \ol{\Delta}^{n+1} = (\ol{\Delta}\otimes \Id^{\otimes n})\circ \ol{\Delta}^n.
    \]
    
    We say that $B$ is \emph{conilpotent} if for every $b\in B$ we have $\ol{\Delta}^N(b) = 0$ for some $N\geq 0$.
\end{definition}

Filtrations of coalgebras will play a significant role in this work. These are most naturally understood through quotients rather than subspaces.
\begin{definition}
    Suppose $(B,\Delta,\eta,d)$ is a differential coalgebra. A \emph{filtration} on $B$ is a sequence of subcomplexes
    \[
    0=F_{-1} B \leq F_0 B\leq F_1 B\leq \ldots \leq F_p B \leq \ldots \leq B
    \]
    such that every quotient $B/F_p B$ is $A$-torsion free, and for all $i\leq p$ the comultiplication factors as indicated by
    \[
    \xymatrix{
    B \ar[r]^-\Delta \ar[d]_{q_p} & B\otimes B \ar[d]^{q_i\otimes q_{p-i}} \\
    \sfrac{B}{F_p B} \ar@{-->}[r]^-{\exists \Delta^p_i} & \sfrac{B}{F_i B} \otimes \sfrac{B}{F_{p-i} B}.
    }
    \]
    Equivalently, $(q_i\otimes q_{p-i})\circ \Delta$ vanishes on $F_p B$.

    The \emph{associated graded} coalgebra is
    \[
    \gr{B} = \bigoplus \gr[p]{B} \text{ where } \gr[p]{B}:= \ker\left( \sfrac{B}{F_{p-1}}\onto \sfrac{B}{F_{p} B}\right) \cong \sfrac{F_p B}{F_{p-1} B}
    \]
    equipped with the well-defined coproduct $\Delta\maps \bigoplus_{p} \gr[p]B \to \bigoplus_{i\leq p} \gr[i]B\otimes \gr[p-i]B
    $ consisting of
    \[
    \gr[p]B \subseteq \sfrac{B}{F_{p-1}B} \xrightarrow{\oplus\Delta^{p-1}_i}\bigoplus_{i\leq p-1} \sfrac{B}{F_i B} \otimes \sfrac{B}{F_{p-1-i}B},
    \]
    whose image lands in $\bigoplus_{i\leq p} \gr[i]B \otimes \gr[p-i] B$.
\end{definition}
\begin{remark}
    Since all $A$-modules in question are torsion free, the sequence
    \[
    0\to \gr[i]B \otimes \gr[p-i] B \to \sfrac{B}{F_{i-1} B} \otimes \sfrac{B}{F_{p-i-1}B} \to  \left(\sfrac{B}{F_{i} B} \otimes \sfrac{B}{F_{p-i-1}B}\right)\oplus \left(\sfrac{B}{F_{i-1} B} \otimes \sfrac{B}{F_{p-i}B}\right)
    \]
    is left-exact. Thus, a quick diagram chase shows that $\Delta^{p-1}_i$ restricts to $\gr[p]B \to \gr[i]B \otimes \gr[p-i]B$, giving a coassociative coproduct on $\gr{B}$.
\end{remark}
\begin{lemma}
    Suppose $B=(B,\Delta,\eta,d)$ is a dg-coalgebra. Let $\H^*(B)$ denote its cohomology. Suppose that $\H^{0}(B)$ is torsion-free, then $\H^0(B)$ has a natural coalgebra structure
    \[
    \H^0(B) \xrightarrow{\H^0(\Delta)} \H^0(B\otimes B) \overset{\sim}\longto \H^0(B)\otimes \H^0(B)
    \]
    where the second arrow is inverse to the K\"unneth isomorphism.
    
    Moreover, given a filtration $F_\bullet \leq B$ such that all $\H^0(\sfrac{B}{F_p B})$ are torsion-free, the cohomology coalgebra $\H^0(B)$ is naturally filtered by
    \[
    F_p \H^0(B) := \ker\left( \H^0(B)\to \H^0(\sfrac{B}{F_p B}) \right).
    \]
\end{lemma}
\begin{proof}
    The K\"unneth formula shows that whenever $\H^0(B)$ is torsion-free, the natural homomorphism $\H^0(B)\otimes \H^0(B) \to \H^0(B\otimes B)$ is invertible. Its inverse gives a coproduct structure on $\H^0(B)$. Coassociativity follows immediately from the naturality of the K\"unneth map.

    Now, we claim that $F_\bullet \H^0(B)$ is a filtration of coalgebras. Indeed,  the quotients $\sfrac{\H^0(B)}{F_p \H^0(B)}$ are torsion-free as submodules of $\H^0(\sfrac{B}{F_p B})$, and using the K\"unneth formula again, 
    \[\xymatrixcolsep{5pc}
    \xymatrix{
    \H^0(B) \ar[r]\ar[d] & \H^0(B\otimes B) \ar[d] & 
 \H^0(B)\otimes \H^0(B) \ar[l]_{\sim} \ar[d] \\
 \H^0(\sfrac{B}{F_p B}) \ar@{..>}[r]^-{\H^0(\Delta_i^p)} & \H^0(\sfrac{B}{F_i B} \otimes \sfrac{B}{F_{p-i} B}) & \H^0(\sfrac{B}{F_{i} B}) \otimes \H^0(\sfrac{B}{F_{p-i} B}) \ar[l]_{\sim}
    }
    \]
    so after inverting the K\"unneth isomorphisms we indeed find that the comultiplication factors through $\H^0(\sfrac{B}{F_{p} B})$, as required.
\end{proof}
The filtered coalgebra structure on $\H^0(B)$ is obviously functorial.

\begin{definition}
    Any augmented coalgebra $B$ with no $A$-torsion is equipped with the  \emph{weight filtration} as follows. Let $F_p B = \ker\left( \ol{\Delta}^p\maps B\to B^{\otimes p+1} \right)$ where $\ol{\Delta}$ extends to $B$ via the projection $B\to \ol{B}$ given by $b \mapsto b-\eta(b)\cdot 1$. Coassociativity of $\ol{\Delta}$ implies that the comultiplication factors through $\sfrac{B}{F_p B}$, as required. Clearly, any homomorphism of coalgebras preserves this filtration.
\end{definition}

\subsection{Bar construction}\label{sec:bar construction}
Coalgebras play a central role in our theory through the following construction.
\begin{definition}[Bar construction]
    Let $(C^*, d,\mu)$ be an augmented associative dg algebra, and for $c\in C^i$ denote $|c|=i$. The augmentation ideal of $C$ is denoted by $\ol{C}$. The \emph{Bar construction} of $C$,
    denoted $Bar(C)$, is the totalization of the 2nd quadrant bicomplex
    \[
    Bar^{-p,q}(C) = (\ol{C}^{\otimes p})^q,
    \]
    where the vertical $q$-grading is the internal cohomological grading of the tensor product complex $\ol{C}^{\otimes p}$. Tensors $\ol{C}^{\otimes p}$ are said to have \emph{weight} $p$.
    
    \begin{remark} An equivalent way to construct the totalization $Bar(C)$ is as the tensor (co)algebra
    \[
   T(\desus \ol{C}) = \bigoplus_{p=0}^{\infty} (\desus \ol{C})^{\otimes p},
    \]
    where $\desus \ol{C} := \ol{C}[1]$ is the \emph{desuspension} of $\ol{C}$, i.e. it is the shifted complex $\ol{C}^{i+1}$. Indeed, totalization $(\tot E)^* = \bigoplus_{p=0}^\infty E^{-p,*+p}$ effectively identifies horizontal shifts with shifts in cohomological degree, equivalent to a desuspension operation.
\end{remark}

    The differential on the bar construction is the sum of two differentials. The `vertical' differential $d_1$ of bidegree $(0,1)$ is the extension of $d$ to a graded coderivation, 
    \[
    d_1(c_1\otimes \ldots \otimes c_p) = -\sum_{i=1}^p (-1)^{(i-1)+ |c_1|+\ldots+ |c_{i-1}|} c_1\otimes \ldots \otimes dc_i \otimes \ldots \otimes c_p.
    \]
    The `horizontal' differential $d_2$ of degree $(1,0)$ is defined by multiplying consecutive elements,
    \[
    d_2(c_1\otimes \ldots \otimes c_p) = -\sum_{i=1}^{p-1}(-1)^{i + |c_1|+\ldots+|c_i|} c_1\otimes \ldots\otimes \mu(c_i,c_{i+1}) \otimes \ldots \otimes c_p.
    \]
    Both differentials satisfy the Koszul sign convention, and they anti-commute. The (total) \emph{Bar differential} is $d_{Bar} = d_1+d_2$. Let $\H^*_{Bar}(C)$ denote its cohomology, which we call the \emph{Bar cohomology of $C$}.
\end{definition}

\begin{definition}[Coalgebra structure]
    The bar construction is naturally a conilpotent dg coalgebra, with respect to the deconcatenation coproduct
    \begin{equation}\label{eq:coproduct definition}
        \Delta(c_1|\ldots| c_p) = \sum_{i=0}^{p} (c_1|\ldots|c_i)\otimes (c_{i+1}|\ldots|c_p).
    \end{equation}
    where we replaced the tensor symbols by bars. The augmentation is given by the $p=0$ summand $(\desus\ol{C})^{\otimes 0} = A$ concentrated in cohomological degree $0$, so that
    \[
    \ol{Bar}(C) = \bigoplus_{p\geq 1} (\desus\ol{C})^{\otimes p}
    \]
    with reduced coproduct omitting the $i=0$ and $p$ summands from \eqref{eq:coproduct definition}. Consequently, the weight of a tensor $x$ agrees with the order of vanishing of $\ol{\Delta}^n(x)$.
\end{definition}
\begin{remark}
    When $C$ is a commutative algebra, there is a further product structure on $Bar(C)$, given by the shuffle product, making $Bar(C)$ into a Hopf algebra. This structure is significant in the context of rational homotopy theory, giving rise for example to Malcev Lie algebras of fundamental groups.
\end{remark}
Bar cohomology is a quasi-isomorphism invariant of a dg algebra $C$, in the following sense. 
\begin{lemma}[{\cite[Proposition 2.2.3]{algebraic_operads}}]\label{lem:quasi-iso invariant}
    If $C\to C'$ is a quasi-isomorphism of augmented dg algebras, and both underlying $A$-modules are torsion-free, then the induced map
    \[
    Bar(C) \to Bar(C')
    \]
    is a quasi-isomorphism.
\end{lemma}
\begin{proof}
    The argument of \cite[Proposition 2.2.3]{algebraic_operads} applies verbatim, with the observation that torsion-free $A$-modules are flat. Recounting the argument, one filters the bicomplexes $Bar(C)$ and $Bar(C')$ by weight. Their associated graded consist of tensor power complexes
    \[
    (C^{\otimes p},d_1) \to ((C')^{\otimes p},d_1).
    \]
    Since tensor products of complexes of flat modules are quasi-ismorphism invariant, the map induces an isomorphism of the $E^1$-pages of the associated spectral sequence. It therefore induces an isomorphism on bar cohomology.
\end{proof}

\begin{definition}
    An augmented dg algebra $C^*$ is said to be \emph{connected} if it is quasi-isomorphic to an augmented algebra $\tilde{C}$ with $\tilde{C}^0 = A$.
\end{definition}
\begin{lemma}
    Let $C^*$ be a dg algebra whose underlying $A$-module is torsion-free. Let $Bar_{\leq p}(C)$ denote the sub-coalgebra of tensors of weight $\leq p$. Then the quotients
    \[
    Bar_{> p}(C) := Bar(C)/Bar_{\leq p}(C)
    \]
    are torsion-free so $Bar_{\leq p}(C)$ defines a natural filtration on $Bar(C)$.

    If $C$ is connected, then $\H^0(Bar_{>p}(C))$ is torsion-free for all $p\geq 0$, so that these induce a natural filtration on $\H^0_{Bar}(C)$.
\end{lemma}
\begin{proof}
    The fact that $Bar_{>p}(C)$ is torsion-free is clear, and it gives a filtration of $Bar(C)$ as a special case of the filtration of a conilpotent coalgebra by weight.

    To see that in the connected case $\H^0(Bar_{>p}(C))$ is torsion-free, assume without loss of generality that $C^0 = A$. Then the augmentation ideal $\ol{C}$ is concentrated in positive degrees, and consequently $Bar^0(C)$ is the lowest degree in which $Bar(C)$ is nontrivial. It follows that $\H^0(Bar_{>p}(C))$ is a sub-coalgebra of $Bar_{>p}(C)$, characterized by the vanishing of $d_{Bar}$, thus torsion-free.
\end{proof}
We summarize the main results of this section as follows.
\begin{cor}
    The assignment $C^* \mapsto \H^0_{Bar}(C)$ is a functor which sends connected dg algebras, whose underlying $A$-module is torsion-free, to filtered coassociative coalgebras. Quasi-isomorphisms of algebras are sent to isomorphism of coalgebras.
\end{cor}

\subsection{Group rings} \label{sec:group ring prelims}
Let $\Gamma$ be a group. We recall basic definitions.
\begin{definition}
    The \emph{group ring} $A[\Gamma]$ is the module of $A$-linear combinations of elements $w\in \Gamma$. It forms an augmented (Hopf) algebra, with product
    \[
    m\maps A[\Gamma]\otimes A[\Gamma] \to A[\Gamma]
    \]
    extending $m(w_1\otimes w_2) = w_1w_2$ to an $A$-linear map, with coproduct $\Delta(w) = w\otimes w$, and with antipode $s(w) = w^{-1}$ for every $w\in \Gamma$.

    The augmentation $\epsilon\maps A[\Gamma]\to A$ sends every $w\in \Gamma$ to $1$, so that 
    \[
    \epsilon(\sum \lambda_i \cdot w_i) = \sum \lambda_i.
    \]
    Let $I\lhd A[\Gamma]$ denote the augmentation ideal, which is the kernel of $\epsilon$, generated by $(w-1)$ for $w\in \Gamma$. Then $A[\Gamma]$ is filtered by powers $\ldots \leq I^2\leq I \leq A[\Gamma]$.
\end{definition}
The association $\Gamma \mapsto A[\Gamma]$ is a functor, sending groups to augmented (Hopf) algebras.

Recall the dimension series of a group.
\begin{definition}[\cite{losey1960dimension}]
    Let $\Gamma$ be a group and $A$ any ring. The $A$-\emph{dimension series} of $\Gamma$ modulo $A$ is the descending central series of subgroups
    \[
    \ldots \leq \gamma_k (\Gamma;A) \leq \ldots \leq \gamma_2(\Gamma;A) \leq \gamma_1(\Gamma;A) = \Gamma
    \]
    defined by
    \[
    w\in \gamma_k(\Gamma;A) \iff (w-1)\in I^k \subseteq A[\Gamma].
    \]
\end{definition}
Recall that the \emph{lower central series} of $\Gamma$, defined recursively as $\gamma_1\Gamma = \Gamma$ and $\gamma_{k+1}\Gamma = [\gamma_k\Gamma, \Gamma]$, is the minimal central series. It follows that $\gamma_k\Gamma \subseteq \gamma_k (\Gamma;A)$ for every ring $A$.
\begin{example}
    For $A=\Z$, Passi \cite{passi} described classes of groups for which $\gamma_k\Gamma = \gamma_k(\Gamma;\Z)$. However, starting with $k=4$ there exist groups for which $\gamma_k\Gamma \subsetneq \gamma_k (\Gamma;A)$, as found by Rips \cite{rips}.
    For $A = \Q$, the dimension series is constructed recursively by the rule 
    \[\gamma_{k+1}(\Gamma;\Q) = \{ w\in \Gamma \mid \exists n\geq 1 \text{ s.t. } w^n\in [\gamma_k(\Gamma;\Q),\Gamma] \},\]
    and for $A = \F_p$, the dimension series is the \emph{lower exponent $p$ central series} constructed recursively by
    \[\gamma_{k+1}(\Gamma;\F_p) = [\gamma_k(\Gamma;\F_p),\Gamma]\cdot\gamma_k(\Gamma;\F_p)^p.\]
\end{example}

\begin{definition}
    Given an $A$-module $M$, let $M^*$ denote $ \Hom_A (M,A)$, the module of $A$-linear functionals on $M$.
    For every $p\geq 1$ define an evaluation pairing
    \[
    (M^*)^{\otimes p} \otimes M^{\otimes p} \to A \quad \text{via}\quad \langle \varphi_1\otimes \ldots\otimes\varphi_p \;,\; m_1\otimes \ldots\otimes m_p\rangle = \prod_{i=1}^p \varphi_i(m_i).
    \]
\end{definition}
    Given a function $f\maps \Gamma\to A$, its \emph{linear extension} to $A[\Gamma]$ is the $A$-linear map
    \[
    f\maps A[\Gamma]\to A \quad \text{defined by}\quad f\left( \sum \lambda_i \cdot w_i \right) = \sum \lambda_i \cdot f(w_i),
    \]
    defining an isomorphism $A[\Gamma]^* \cong A^\Gamma$.

\subsection{Tensors and infinite cartesian powers}\label{sec:infinite powers}
    The following technical fact will be crucial for constructing elements in group rings that evaluate nontrivially on  invariants of words.
    \begin{lemma}[{\cite[Theorem 1(a)]{goodearl}}] \label{lem:tensor of infinite powers}
    Let $A$ be a commutative Noetherian ring. Then for any two sets $S$ and $T$, the natural map evaluation $A^S\otimes A^T \to A^{S\times T}$ is injective.
\end{lemma}
\begin{proof}
    \cite[Theorem 1]{goodearl} shows that for any $A$-module $M$ the natural map
    \[
    M\otimes A^T \to M^T
    \]
    is injective. Taking $M = A^S$, the claim follows from the adjunction $(A^S)^T \cong A^{S\times T}$.
\end{proof}

\section{Free groups}

\label{sec:free groups}
We start by defining  letter-braiding polynomials of words in free groups. In \S\ref{sec:introducing relations} we will find which of these remain well-defined once relations have been introduced.  Our first approach in \S\ref{sec:topological construction} is  topological, with few calculations, which are then the focus of \S\ref{sec:weight reduction}.

\subsection{Topological construction}\label{sec:topological construction}
The invariants we construct are grounded in the bar construction for simplicial cochains on  simplicial sets. 

\begin{definition}
    
Given a pointed simplicial set $X_\bullet$, let $C^*_{\cup}(X;A)$ denote its simplicial cochain complex with coefficients in $A$.  This is the complex of functions on nondegenerate simplices of $X$ (or dual to the normalized chains). The cup product makes this into an associative algebra, and the basepoint provides an augmentation. Denote
\[
Bar(X;A) := Bar(C^*_\cup(X;A)),
\]
and $\H^*_{Bar}(X;A)$ for its cohomology.
\end{definition}

The following is a standard application of the Quillen adjunction between  the realization and singular functors and the quasi-isomorphism invariance of the bar construction.

\begin{lemma}\label{lem:induced maps on bar cohomology}
    The assignment $X \mapsto Bar(X;A)$ is a contravariant functor from  simplicial sets to filtered dg coalgebras.

    A pointed continuous map on geometric realizations $g\maps |X|\to |Y|$ induces a homomorphism of filtered coalgebras $g^*\maps \H^*_{Bar}(Y;A) \to \H^*_{Bar}(X;A)$ that furthermore satisfies
    \[
    (g_2\circ g_1)^* = g_1^* \circ g_2^* \quad\text{ and }\quad f\sim g \text{ homotopic }\implies f^* = g^*.
    \]
    If $g$ is the geometric realization of a simplicial map $|f|$, then $f^* = g^*$ on bar cohomology.
\end{lemma}

The definition of letter braiding requires calculations on simplicial models for circles and wedges of circles. 

\begin{definition}\label{notation:Xs}

Let $\Delta[1]$ be the standard $1$-simplex, with endpoints $\partial_0\Delta[1]$ and $\partial_1\Delta[1]$.
   
    The \emph{standard circle} is the pointed simplcial set
    \[
    S^1_\bullet := \Delta[1]/(\partial_0 \Delta[1] \sim \partial_1\Delta[1]).
    \]

Let $F_S$ be a free group on basis $S$.
    The \emph{presentation complex of $F_S$} is the pointed simplicial set
    \[
    (X_S)_\bullet = \bigvee_{s\in S} S^1_\bullet,
    \]
    having one standard loop for every generator $s\in S$ of $F_S$.
\end{definition}

The letter-braiding number is defined given the following choice of isomorphism.
\begin{lemma}
    There is a canonical isomorphism of filtered augmented coalgebras
    \[
    \H^*_{Bar}(S^1_\bullet; A) \cong A[t],
    \]
    with $|t|=0$, and $A[t]$ considered as the cofree coalgebra generated by $t$.
\end{lemma}
\begin{proof}
    Let $\tau\in C^1(S^1_\bullet;A)$ be the cochain that assigns the value $1$ to the unique nondegenerate $1$-simplex of $S^1$. The complex $C^*_\cup(S^1_\bullet;A)$ is the augmented algebra $A\cdot 1 \oplus A\tau$ with $d(\tau)=0 = \tau^2$. Therefore, its bar construction is the cofree coalgebra
    \[
    Bar(A[\tau]/(\tau^2)) = T(A\cdot \desus\tau) = \bigoplus_{p\geq 0} A(\desus\tau)^{\otimes p}
    \]
    with vanishing differentials. Denoting $t=\desus\tau$, this coalgebra is polynomial in $t$.
\end{proof}

    More generally, for the wedge of circles $X_S$, the cochain algebra naturally consists of a unit along with
    \[
    C^1_\cup(X_S;A) \cong A^S,
    \]
    the $A$-valued functions on the set $S$ with trivial cup product. This is in turn naturally interpreted as $\H^1(X_S;A) \cong \Hom(F_S,A)$.
    As in the last observation, the bar construction is concentrated in degree $0$, with vanishing differentials. The bar cohomology is thus naturally identified with the cofree coalgebra
    \[
    \H^*_{Bar}(X_S;A) \cong T(A^S) = \bigoplus_{p\geq 0} \left(A^S \right)^{\otimes p}.
    \]

\begin{definition}
    Recall that elements $w\in F_S$ are in bijection with based homotopy classes of pointed maps $|S^1|\to |X_S|$. Abusing notation, let $w$ denote a representative of the respective homotopy class.

    The \emph{letter-braiding polynomials} of the element $w\in F_S$ are the composition
        \[
        L_{\bullet}(w) \maps T( A^S ) \cong H^0_{Bar}(X_S;A) \overset{w^*}\longto \H^0_{Bar}(S^1;A) \cong A[t].
        \]
    The \emph{letter-braiding numbers} of $w$ are the linear coefficient of the respective polynomials
        \[
        \ell_\bullet(w) := \frac{d}{dt}L_\bullet(w)|_{t=0}\maps T(A^S) \longto A.
        \]
    \end{definition}
    Since induced maps on bar cohomology depend only on homotopy class, the polynomials $L_\bullet(w)$ are independent of the choice of representative for $w$.

    \begin{example}\label{ex:braiding of a generator}
        The easiest calculation of braiding polynoimals is that of a generator $s\in S$. In this case, given functionals $\alpha_1,\ldots,\alpha_n\in A^S$, we have
        \[
        L_{(\alpha_1|\ldots|\alpha_n)}(s) = \alpha_1(s)\alpha_2(s)\ldots\alpha_n(s) \cdot t^n \quad\text{ and }\quad L_1(s) = 1.
        \]
        Indeed, let $i_s\maps S^1_\bullet \to (X_S)_\bullet$ be the inclusion of the edge labeled by $s$. On $1$-cochains it defines the pullback $i_s^*\maps C^1(X_S;A)\to C^1(S^1;A)$ given by restriction $A^S \to A^{\{s\}}$. That is, $i_s^*(\alpha) = \alpha(s)\cdot \tau$, where $\tau\in A^{\{s\}}$ is the functional $s\mapsto 1$. Recalling that the desuspension of $\tau$ in the bar complex of $S^1_\bullet$ is $t$, the general calculation follows.
        
        The braiding number is therefore
        \[
        \ell_{(\alpha_1|\ldots|\alpha_n)}(s) = \begin{cases}
            \alpha_1(s) & n=1\\
            0 & \other.
        \end{cases}
        \]
    \end{example}

    Recall that $A[t]$ is cofreely cogenerated by $t$, which satisfies $\Delta(t^n) = \sum_{i=0}^{n}t^i\otimes t^{n-i}$. With this, the coproduct of braiding polynomials is readily computable.
    \begin{theorem}
        The letter-braiding polynomial is a homomorphism of augmented filtered coalgebras, meaning that $L_{(T+a\cdot T')}(s) = L_T(w)+a\cdot L_{T'}(w)$ and if $\Delta(T) = \sum_i T_1^i\otimes T_2^i$ then
        \[
        \Delta L_T(w) = \sum_i L_{T_1^i}(w) \otimes L_{T_2^i}(w) \quad \text{and}\quad L_1(w) = 1.
        \]

        Furthermore, the braiding polynomial is determined by the braiding number, in  that
        \[
        L_T(w) = \eta(T)\cdot 1 + \sum_{k=0}^{\infty} (\ell(w)^{\otimes k+1})_{\ol\Delta^k(T)} \cdot t^{k+1}
        \]
        where $\ol\Delta^k$ is the iterated reduced coproduct $(\ol\Delta\otimes \operatorname{Id}^{\otimes k-1})\circ\ldots\circ(\ol\Delta\otimes \operatorname{Id})\circ \ol\Delta$, and $(\ell(w)\otimes \ell(w))_{T_1\otimes T_2} = \ell_{T_1}(w)\cdot \ell_{T_2}(w)$.
    \end{theorem}
    \begin{proof}
        Functoriality of the bar construction implies that the letter braiding polynomial $L(w)$ of any $w\in F_S$ is a coalgebra homomorphism. Thus, for every $T$ with $\Delta(T) = \sum_i T_1^i\otimes T_2^i$,
        \[
        (\Delta \circ L(w))_T = (L(w)\otimes L(w))\circ \Delta(T) = \sum_i (L(w)\otimes L(w))_{T_1^i\otimes T_2^i} = \sum_i L_{T_1^i}(w) \otimes L_{T_2^i}(w)
        \]
        as claimed.
        
        Now, suppose $L_T(w) = \sum_{k=0}^N a_k t^k$. Applying the reduced coproduct $k$ times to $t^{k+1}$ gives the pure tensor $\underbrace{(t\otimes\cdots\otimes t)}_{k+1 \text{ terms}}$. On the other hand,
        \[
        a_k \cdot (t\otimes \cdots \otimes t) + \text{ higher order terms } = \ol\Delta^k L_T(w) = L(w)^{\otimes k+1}({\ol\Delta^k(T)}),
        \]
        so the coefficient $a_k$ is the product of the linear terms of the respective braiding polynomials of the factors of $\ol{\Delta}^k(T)$.
    \end{proof}

    We next unpack the sense in which braiding numbers are functorial with respect to group homomorphisms of free groups.  
    On bar cohomology
    \[
    \H^0_{Bar}(X_S;A) \cong T(A^S)
    \]
    is the cofree coalgebra generated by functions $S\to A$. By the universal property of cofree coalgebras, coalgebra homomorphisms $f\maps (B,\Delta) \to \H^0_{Bar}(X_S;A)$ are in bijection with $A$-linear maps $\ell f\maps B\to A^S$. Thus to give such a homomorphism one must only provide an element $\ell f_b \in A^S$, that is a function $\ell f_b \maps s\mapsto \ell f_b(s)\in A$ for every $b\in B$, depending linearly on $b$.
    \begin{theorem}\label{thm:naturality}
        Given a homomorphism of free groups $h\maps F_{S_1}\to F_{S_2}$, there is a functorial pullback $\H^*_{Bar}(X_{S_2};A) \to \H^*_{Bar}(X_{S_1};A)$ coinciding with the coalgebra homomorphism
        \[
        h^*\maps T(A^{S_2}) \to T(A^{S_1})
        \]
        determined by the functions $\ell h^*_T\in A^{S_1}$ for every $T\in T(A^{S_2})$ given by
        \[
        \ell h^*_T\maps s \mapsto \ell_{T}(h(s)) \quad \text{for all} \quad s\in S_1.
        \]
        These homomorphisms satisfy $(h\circ h')^* = (h')^*\circ h^*$.
    \end{theorem}
    \begin{example}\label{ex:homomorphism from cyclic group}
        Let $S_1 = \{s\}$ generate a cyclic group $F_{\{s\}}$. Then $T(A^{\{s\}}) \cong A[t]$ where $t$ corresponds to the functional $s\mapsto 1$. An element $w\in F_{S_2}$ corresponds to $h\maps F_{\{s\}} \to F_{S_2}$, and the induced map $h^*\maps T(A^{S_2}) \to A[t]$ is exactly the letter-braiding polynomial $L_\bullet(w)$.
    \end{example}
    \begin{example}\label{ex:pullback of product}
        Let $S_{1} = \{e_1,e_2\}$ generated a free group of rank two, so that $A^{S_1} = A\langle E_1,E_2\rangle$ generated by the `dual' basis $E_i(e_j) = \delta_{ij}$.
        A homomorphism $h\maps F_{S_1}\to F_{S_2}$ corresponds to two elements $(w_1,w_2) = (h(e_1),h(e_2))$, and its pullback is
        \begin{align} \label{eq:pullback of two words}\begin{split}
        h^*(T) =  & \; \eta(T)\cdot 1 + \ell_{T}(w_1)\cdot E_1 + \ell_{T}(w_2)\cdot E_2 \\
           + \sum & \; \ell_{T_1^i}(w_1)\ell_{T_2^i}(w_2)\cdot E_1|E_2 + \ell_{T_1^i}(w_2)\ell_{T_2^i}(w_1)\cdot E_2|E_1 \\ + & \; \ell_{T_1^i}(w_1)\ell_{T_2^i}(w_1)\cdot E_1|E_1 + \ell_{T_1^i}(w_2)\ell_{T_2^i}(w_2)\cdot E_2|E_2 \\
         + & \text{ weight } \geq 3
        \end{split}
        \end{align}
        where $\ol\Delta(T) = \sum_i T_1^i\otimes T_2^i$.

        Composing $h$ with the homomorphism $F_{\{s\}}\to F_{S_1}$ sending $s\mapsto e_1 e_2$, functoriality implies \[
        L_{T}(w_1\cdot w_2) = \eta(T) \cdot 1 + \ell_{T}(w_1)\cdot  L_{E_1}(e_1 e_2) + \ell_{T}(w_2) \cdot L_{E_2}(e_1 e_2) + \Sigma\, \ell_{T_1^i}(w_1)\ell_{T_2^i}(w_2) \cdot L_{E_1|E_2}(e_1 e_2) + \ldots 
        \]
        thus computing the letter-braiding polynomial of a product $w_1\cdot w_2$ from polynomials of $F_2$.
    \end{example}
    \begin{proof}[Proof of Theorem \ref{thm:naturality}]
        A group homomorphism $h\maps F_{S_1}\to F_{S_2}$ corresponds to a unique based homotopy class of pointed maps $h\maps |X_{S_1}|\to |X_{S_2}|$. We show that the induced map $h^*$ on bar cohomology has weight $1$ component as stated. 
        
        For every $T\in \ol{T}(A^{S_2})$ consider its pullback $h^*_T\in \ol{T}(A^{S_1})$ and let $\ell h^*_T\in A^{S_1}$ be its image under the projection onto the weight $1$ component $\ell\maps \ol{T}(A^{S_1})\onto A^{S_1}$.
        Fixing $s\in S_1$, we calculate $\ell h^*_T(s)$ by considering the inclusion of the circle labeled by $s$
        \[
        i_s\maps S^1 \into \bigvee_{s'\in S_1}S^1 = X_{S_1}.
        \]
        
        On geometric realizations,
        \[
        |S^1| \overset{|i_s|}\longrightarrow |X_{S_1}| \overset{h}\longrightarrow |X_{S_2}|
        \]
        has composition representing the homotopy class of $h(s)\in F_{S_2}$. By definition of the letter-braiding polynomial, this gives $(i_s^*\circ h^*)(T) = L_{T}(h(s))$, having linear coefficient $\ell_{T}(h(s))$. 
        On the other hand, as explained in Example \ref{ex:braiding of a generator}, $i_s^*(\alpha_1|\ldots|\alpha_n) = \alpha_1(s)\ldots\alpha_n(s) \cdot t^n$ and has vanishing linear coefficient unless $n=1$.
        It follows that
        \[
        \ell_{T}(h(s)) = \ell( i_s^* \circ h^*(T) ) = \ell h_T^*(s),
        \]
        which proves the claimed description of $h^*$. Functoriality is automatic by Lemma \ref{lem:induced maps on bar cohomology}.
\end{proof}

\subsection{Weight reduction algorithm}
\label{sec:weight reduction}
We give a combinatorial construction to calculate the letter-braiding invariants of an element $w\in F_S$,  Let $[n]:= \{1,\ldots,n\}$ and $[n]_*:= \{0,\ldots, n\}$, considered modulo $n+1$.
\begin{definition}\label{def:strings}
    A \emph{word} in the generators of $F_S$ is a map $\tilde{w}\maps [n] \to S\cup S^{-1}$.
    The \emph{loop} associated to $\tilde{w}$ is the simplicial circle 
    \[
    S^1(\tilde{w})_\bullet = \left( [n]_* \times \Delta[1] \right) / \sim
    \]
    where $\sim$ glues the endpoints of the intervals in order to get a circle, but with orientation determined by the sign of the respective letter, as indicated by the following examples:
    \[
    \xymatrix{
    1\ar[r]^{s_1} & 2 \ar[r]^{s_2} & 3 & 4 \ar[l]_{s_3^{-1}} & 5 \ar[l]_{s_4^{-1}} \ar[r]^{s_5} & 0 \ar@/^2pc/[lllll] \\
    &&\phantom{4}&&&&
    }
    \xymatrix{
    1 & 2 \ar[l]_{s_1^{-1}}\ar[r]^{s_2} & 3 \ar[r]^{s_3} & 4 \ar[r]^{s_4} & 5 & 0 \ar@/^2pc/[lllll] \ar[l]_{s_5^{-1}} \\
    &&&\phantom{4}&&&
    }
    \]
The first segment $0\to 1$ does not correspond to a letter and is always oriented in increasing order.  We refer to it as the \emph{standard segment}.

Let $|s| = |s^{-1}| = s$ for every $s\in S$. There is a canonical simplicial map $\tilde{w}$ from $ S^1(\tilde{w})_\bullet$
to the wedge of circles $ (X_S)_\bullet$, sending the $i$-th segment to the circle indexed by $|\tilde{w}(i)|\in S$.
The orientation of the segment matches that of the target circle. The (unlabeled) standard segment is collapsed onto the basepoint.  We abuse notation\footnote{Maps from simplicial circles to $X_S$ are in obvious bijection with words in $S\cup S^{-1}$. Here we insist that all such maps start with a collapsed standard segment for reasons that will become apparent next.} and refer to this simplicial map by $\tilde{w}$.

For brevity, let $\H^*_{Bar}(\tilde{w};A)$
denote $ \H^*_{Bar}(S^1(\tilde{w});A)$.
\end{definition}

\begin{definition}
    Let the \emph{standard collapse} $\std \maps S^1(\tilde{w})\to S^1$ be the pointed simplicial map defined by collapsing all edges except for the standard segment $0\to 1$, which is sent to the nondegenerate segment of $S^1$.
    It induces an isomorphism of filtered coalgebras
    \[
    L^{-1}\maps {A[t]} = \H^0_{Bar}(S^1;A) \to H^0_{Bar}(\tilde{w};A)
    \]
    whose inverse we denote by $L\maps H^0_{Bar}(\tilde{w};A)\to {A[t]}$.
\end{definition}

These maps exist and are well-defined because the collapse map $S^1(\tilde{w})\to S^1$ is a homotopy equivalence on geometric realizations, and so induces an isomorphism on bar cohomology.

\begin{cor}\label{cor:combinatorial model}
    Let $\tilde{w}\maps [n]\to S\cup S^{-1}$ be a word in the generators of $F_S$, and let $w\in F_S$ denote the element it represents. The letter-braiding polynomials $L_\bullet(w)$ of $w$ are given by the composition
        \[
        {T}( A^S ) \overset{\tilde{w}^*}\longto \H^0_{Bar}(\tilde{w};A) \overset{L}{\longto} {A[t]}.
        \]
    \end{cor}
    \begin{proof}
        We have a commutative diagram of pointed simplicial sets
        \[
        \xymatrix{
        & S^1(\tilde{w}) \ar[ld]_{\std}^{\sim} \ar[rd]^{\tilde{w}} &\\
        S^1 \ar[d]^{\sim} & & X_S \ar[d]_{\sim} \\
        \operatorname{Sing}(|S^1|) \ar[rr]^{\operatorname{Sing}(w)} & & \operatorname{Sing}(|X_S|)
        }
        \]
        where the arrows decorated with $\sim$ induce quasi-isomorphisms on simplicial cochain algebras. It follows that the two paths $S^1(\tilde{w}) \to \operatorname{Sing}(|X_S|)$ induce the same map on bar cohomology. By the construction of $w^*\maps {T}(A^S) \to A[t]$ given in Lemma \ref{lem:induced maps on bar cohomology}, it agrees with $L_\bullet(w)$.
    \end{proof}

We unpack the last corollary and give an explicit description of $L_\bullet(w)$. We use the \emph{sign function} $(-1)^\bullet\maps S\cup S^{-1} \to \{ \pm 1 \}$, defined by $(-1)^s=1$ and $(-1)^{s^{-1}}=-1$ for every $s\in S$.
\begin{definition}\label{def-lem:chains on a word}
    Given a word $\tilde{w}\maps [n] \to S\cup S^{-1}$, let $C^*(\tilde{w};A)$ be the following dg algebra.
    \[
    C^0(\tilde{w};A) = C^1(\tilde{w};A) := \{ a\maps [n]_* \to A \} \text{ and otherwise } C^i = 0,
    \]
    with $d:C^0\to C^1$ given by $(df)_i :=(-1)^{\tilde{w}(i)}(f_{i+1} - f_{i})$ and $i+1$ is computed modulo $n$. Elements of $C^0(\tilde{w};A)$ are called \emph{functions} and denoted by Latin letters, such $f$ and $g$.
    Elements of $C^1(\tilde{w};A)$ are called \emph{forms}, and denoted by Greek letters.
    
    The product is defined on functions by $(f \smile g)_i = f_i\cdot g_i$, and between functions and forms as follows
    \begin{align}
       \begin{aligned}
       \text{ if } \tilde{w}(i)\in S: \\
            (f \smile \alpha)_i &= f_i\cdot \alpha_i \\
            (\alpha \smile f)_i &=  f_{i+1}\cdot \alpha_i
       \end{aligned} 
       \text{ and, }
        \begin{aligned}
       \text{ if } \tilde{w}(i)\in S^{-1}: \\
            (f \smile \alpha)_i &= f_{i+1}\cdot \alpha_i \\
            (\alpha \smile f)_i &=  f_{i}\cdot \alpha_i
       \end{aligned} 
    \end{align}
    with $i+1$ computed modulo $n$, and the case $i=0$ satisfying the same formulas as $\tilde{w}(i)\in S$.    
    The augmentation is given by the evaluation $f \longmapsto f_0$.
\end{definition}

The following is immediate from definitions.

\begin{lemma}
$C^*(\tilde{w};A)$
     is isomorphic to $C^*_{\cup}(S^1(\tilde{w});A)$ with its standard cup product.
\end{lemma}

\begin{comment}
    
\begin{proof}
    The nondegenerate simplices of $S^1(\tilde{w})$ are precisely $n$ vertices and $n$ edges, arranged cyclically. Label the vertices in cyclic order by $[n]_*$, starting with the basepoint corresponding to $0\in [n]_*$. Then thinking of the edge between $i$ and $i+1$ as labeled by $i\in [n]_*$, we get a bijection between edges and $[n]_*$, where the standard segment corresponds to $0\in [n]_*$. Therefore the cochain groups $C^*_\cup(S^1(\tilde{w});A)$ are the functions on $[n]_*$ in degrees $0$ and $1$.

    Quick examination of the examples in Definition \ref{def:strings} shows that the simiplicial differential agrees with the one in $C^*(\tilde{w};A)$. For the cup products, recall that for any $1$-simplex $e_i$,
    \[
    f\smile \alpha( e_i ) = f(\partial_0(e_i)) \cdot \alpha(e_i) \quad\text{ and }\quad \alpha\smile f(e_i) = \alpha(e_i)\cdot f(\partial_1(e_i)).
    \]
    These formulas specialize to the cup products given in $C^*(\tilde{w};A)$.
\end{proof}

\end{comment}

%We next explain how this cochain model facilitates calculation.

%\subsection{Calculation via weight reduction}\label{sec:weight reduction}

The combinatorial model of the letter-braiding polynomial using $C^*(\tilde{w};A)$ supports an explicit algorithm for calculation, distinct from Fox calculus, using weight reduction as introduced by Monroe--Sinha \cite{monroe-sinha}. 
%Let us present this algorithm.

%Let $\tilde{w}\maps [n]\to S\cup S^{-1}$ be a word representing some element in $F_S$, and recall that $C^0(\tilde{w};A)$ and $C^1(\tilde{w};A)$ are both the set of functions $[n]_*=\{0,\ldots,n\}\to A$.

We build forms from the following special forms in $C^1(\tilde{w};A)$.

\begin{definition}
 For $0\leq i \leq n$, the \emph{indicator form} $\delta_i\maps [n]_*\to A$ is
        \[
        \delta_{ij} = \begin{cases} 1 & i=j \\
        0 &\other,
        \end{cases}
        \]
        where we also denote $\delta_0$ by $t$.
 The \emph{orientation form} $dx\in C^1(w;A)$ is
    \[
    (dx)_i := \begin{cases}
        +1 & w(i)\in S \\
        -1 & w(i) \in S^{-1}\\
        0 & i=0.
    \end{cases}
    \]
    
    Given a function $f\in [n]\to A$, let $fdx$ denote the form
    \[
    (fdx)_i := f_i \cdot (dx)_i \quad \text{with} \quad (fdx)_0=0.
    \]
    \end{definition}

Our definition of $f dx$ is \emph{not} an instance of the cup product on $C^*(\tilde{w};A)$.

\begin{definition}
    
    For every $1\leq i\leq j\leq n+1$, the \emph{integral} of $fdx\in C^1(\tilde{w};A)$ is
    \[
    \int_i^j fdx := f_i+ f_{i+1} + \ldots + f_{j-1},
    \]
    and $\int_w fdx = f_1+\ldots+f_n$.

    The \emph{cobounding} of a form $fdx\in C^1(\tilde{w};A)$ is the function $d^{-1}(fdx)\in C^0(\tilde{w};A)$ defined by
    \[
    d^{-1}(fdx)_j := \int_1^j fdx - \int_w fdx  = -\left( f_{j} +f_{j+1} + \ldots + f_{n} \right)
    \]
    and $d^{-1}(fdx)_0 = 0$.
\end{definition}

The following facts are easy to check.
\begin{lemma}
    Any form $\alpha\in C^1(\tilde{w};A)$ can be written uniquely as
\[
\alpha = fdx - a_0 \cdot \delta_0.
\]
Given a function $F\in C^0(\tilde{w};A)$, its differential is
    \[
    dF = F'dx - F_1\cdot \delta_0,
    \]
    where $(F')_i = F_{i+1}-F_i$. In particular, the cobounding of $fdx$ belongs to the augmentation ideal of $C^*(\tilde{w};A)$ and satisfies
    \[
    d( d^{-1}(fdx) ) = fdx - \left(\int_w fdx \right) \cdot \delta_0.
    \]
    Thus $fdx$ is exact if and only if $\int_w fdx = 0$.
\end{lemma}

\begin{prop}[0-linking and weight reduction]\label{prop:weight-reduction}
    Every element $T\in Bar^0(\tilde{w};A) = T(C^1(\tilde{w};A))$ is cohomologous to a polynomial in $\delta_0$.  The isomorphism
    \[
    L\maps \H^0_{Bar}(\tilde{w};A) \to {A[t]}
    \]
    defining the letter-braiding polynomial sends a polynomial in $\delta_0$ to one in $t$ by substitution $\delta_0 \mapsto t$. In the notation below, we identify $\delta_0$ with $t$.
    
    The polynomial can be computed via the following weight-reduction procedure. Given a tensor $(\alpha_1|\ldots|\alpha_n)$ such that every $\alpha_i$ is either $t$ or $fdx$, let its \emph{modified weight} be $n-\#(t \text{ entries})$.
    \begin{enumerate}
        \item  If the modified weight is $0$, then the tensor is already a polynomial in $t$, and we halt. Otherwise, pick one $\alpha_k = fdx$ (preferably one with $\int_w fdx = 0)$ and proceed.
        \item The tensor $(\alpha_1|\ldots|\alpha_{k-1}|fdx|\alpha_{k+1}|\ldots|\alpha_n)$ is cohomologous to 
        \begin{align*} %\label{eq:weight-reduced expression}
        \begin{split}
        \left( \int_w fdx \right)   \cdot &(\alpha_1|\ldots|\alpha_{k-1}|t|\alpha_{k+1}|\ldots|\alpha_n )
        \\
         - &(\alpha_1|\ldots|\alpha_{k-2}|  (\alpha_{k-1}\smile d^{-1}(fdx))|\alpha_{k+1}|\alpha_{k+2}|\ldots|\alpha_n)\\
        + &(\alpha_1|\ldots|\alpha_{k-2}|  \alpha_{k-1}|(d^{-1}(fdx) \smile \alpha_{k+1})|\alpha_{k+2}|\ldots|\alpha_n ). 
        \end{split}
        \end{align*}
        By convention, if $(fdx)$ occurs in one of the tensor's extremal terms, omit the undefined cup product.
    These cup products are given by the following rules:
    \begin{align*}
        t \smile d^{-1}(fdx) = -\left(\int_w fdx\right) \cdot t \quad &\text{ and }\quad  d^{-1}(fdx)\smile t = 0 \\
        gdx \smile d^{-1}(fdx) = -\left(\int_{\partial_1(\bullet)}^{n+1}fdx\right)g\cdot dx  \quad &\text{ and }\quad
            d^{-1}(fdx)\smile gdx  = -\left(\int_{\partial_0(\bullet)}^{n+1}fdx\right)g\cdot dx,
    \end{align*}
    where $\partial_0(\bullet)$ (or $\partial_1(\bullet)$) is the function on segments that returns the initial vertex (or terminal vertex), and we define $\int_0^{n+1}fdx = 0$.
    \item The summands obtained in Step (2) have strictly smaller modified weight than the original. Collect like terms and return to Step (1) for each resulting summand.
    \end{enumerate}
\end{prop}

%\begin{remark}\label{rmk:explain why letter-braiding}
    We refer to this weight-reduction procedure as \emph{$0$-linking}, since general linking of submanifolds is obtained by cobounding and intersecting
    \[
    \operatorname{linking}(M_1,M_2) = \partial^{-1}(M_1)\cap M_2,
    \]
    which is Poincar\'{e} dual to computing $d^{-1}(-)$ followed by the cup product of cochains. This manifests above with linking of $0$-chains (dual to $1$-forms) in a $1$-manifold (the circle $S^1(w)$).

\begin{example}\label{ex:body-diagram weight reduction}
    Example \ref{ex:intro-free group} depicts the higher $0$-linking of letters in $[xy,x^{-2}]\in F_{\{x,y\}}$. The picture describes the weight-reduced form $X\smile d^{-1}(Y\smile d^{-1}(X\smile d^{-1}(X)))$: the drawn intervals at every level depict the contributions of letters to the respective cobounding function, and the dots record cup products of forms and cobounding functions. Explicitly, an occurrence of the letter $x$ contributes $+1$ to the cobounding function $d^{-1}X$ in all subsequent positions, so a solid line is drawn from every occurrence of $x$ and extending to its right, representing addition of $+1$ to the cobounding function. Similarly, occurrences of $x^{-1}$ contribute $-1$ to the cobounding function, as represented by a dotted line. For cup products, a solid dot represents a form with mass of $+1$ on an interval, while an empty dot represents mass $-1$.
    
    For example, at the first step the cobounding function $d^{-1}X$ takes value $+1$ over the 2nd segment and $-1$ on the last three. Its cup product with $X$ assigns mass $+1$ to the 3rd segment and $-1$ to the last.
\end{example}
    
    A special case of these higher linking invariants of letters appeared recently in Monroe--Sinha \cite{monroe-sinha}: their forms $\Lambda_a, \Lambda_b,\ldots $ are the indicator forms we would call $A, B ,\ldots$ above (compare the visual depiction of letter linking in \cite[Example 2.11]{monroe-sinha} with our depiction of letter braiding in Example \ref{ex:intro-free group}).
    The weight-reduction algorithm given above is based on their work and is also the starting point of a sibling paper to this one based on the Harrison complex \cite{joint}.  This algorithm gives an alternative to Fox calculus, that is both faster and simpler to carry out by hand.

\begin{proof}
    The isomorphism $L^{-1}\maps {A[t]}\to \H^0_{Bar}(\tilde{w};A)$ is induced by the simplicial map $S^1(\tilde{w})\to S^1$ collapsing all but the standard $0$-th segment. On cochains, this induces the map $t \mapsto \delta_0$. It follows that any cocycle in $Bar^0(\tilde{w};A) = {T}(C^1(\tilde{w};A))$ is cohomologous to a polynomial in $\delta_0$, obtained by the substitution $t\leftrightarrow \delta_0$.

    To compute this polynomial, one can follow the steps described in the proposition. The coboundary
    \begin{align*}
    d_{Bar}(\alpha_1|\ldots|\alpha_{k-1}|d^{-1}(fdx)|\alpha_{k+1}|\ldots|\alpha_n) = - (\alpha_1|\ldots| \bigg[fdx - (\int_w fdx)\cdot \delta_0 \bigg]|\ldots|\alpha_n) \\
    - (\alpha_1|\ldots|\alpha_{k-1}\smile d^{-1}(fdx)|\alpha_{k+1}|\ldots|\alpha_n) + (\alpha_1|\ldots|\alpha_{k-1}|d^{-1}(fdx)\smile \alpha_{k+1}|\ldots|\alpha_n)
    \end{align*}
    shows that $(\alpha_1|\ldots|fdx|\ldots|\alpha_n)$ and the expression in Step (2) are cohomologous. The cup products described in Step (3) follow directly from Definition \ref{def-lem:chains on a word}.
\end{proof}
With this, one gets an explicit procedure for computing letter-braiding polynomials.
\begin{cor}
    Given $(\alpha_1|\ldots|\alpha_n)\in (A^S)^{\otimes n}$, the letter-braiding polynomial of $w\in F_S$ is obtained by applying the weight-reduction procedure of Proposition \ref{prop:weight-reduction} to the tensor
    \[
    (\tilde{w}^*(\alpha_1)|\ldots|\tilde{w}^*(\alpha_n)) \in C^1(\tilde{w};A)^{\otimes n}
    \]
    where $\tilde{w}^*(\alpha)$ is the composition $[n] \overset{\tilde{w}}\to S\cup S^{-1} \overset{\alpha}\to A$ extended by zero to $0\in [n]_*$.
\end{cor}
\begin{proof}
    In Corollary \ref{cor:combinatorial model} we saw that $L_\bullet(w)$ is obtained as the composition $L\circ w^*$. The weight-reduction procedure gives a formula for $L$, while $w^*$ is induced by $S^1(w)\to X_S$ which pulls back 1-cochains $\alpha \in A^S$ as claimed.
\end{proof}

%\begin{remark}\label{rmk:Monroe-Sinha}
    Monroe--Sinha \cite{monroe-sinha} consider the above weight-reduction procedure in the special case of tensors $\tilde{w}^*(\alpha_1|\ldots|\alpha_n)$ for a word $w\in F_S$ in the $n$-th step of the lower central series. In that case, they shows that all integrals $\int_w fdx$ vanish in every iteration of the weight-reduction process. As a result, their braiding polynomials never involve $\delta_0$ until the very last step, resulting in a linear polynomial, as echoed in Corollary \ref{cor:degree bound} below.
%\end{remark}

\begin{example}\label{ex:linear invariants}
    For tensors of weight $1$, letter-braiding factors through the abelianization in the following sense: suppose $\alpha\in A^S  \cong \H^1(F_S;A) \cong \Hom_{grp}((F_S)^{ab},A)$, then $\forall w\in F_W$
    \[
    L_\alpha(w) = \alpha(w) \cdot t.
    \]
    Indeed, writing the pullback as $\tilde{w}^*(\alpha) = fdx$, Proposition \ref{prop:weight-reduction} shows that
    \[
    \tilde{w}^*(\alpha) \sim \left(\int_w fdx\right) \cdot t
    \]
    and one can check that the integral agrees with $\alpha(w)$. Thus, we need only to count occurrences of letters in $w$ ignoring their ordering, hence factoring through the abelianization.
\end{example}

    Next consider the universal commutator $[x,y] = xyx^{-1}y^{-1}$ in $F_2$. Using the weight reduction procedure, its letter-braiding number on $(\alpha_1|\alpha_2)\in (A^{\{x,y\}})^{\otimes 2}$ is
    \[
    \ell_{(\alpha_1|\alpha_2)}([x,y]) = \alpha_1(x)\alpha_2(y) - \alpha_1(y)\alpha_2(x)
    \]
    as indicated by the following weight-reduction diagram.
    \[
    \begin{aligned}
    \begin{tikzpicture}
\draw[->] (0.05,0) -- (1.45,0);
\draw[->] (1.55,0) -- (2.95,0);
\draw[<-] (3.05,0) -- (4.45,0);
\draw[<-] (4.55,0) -- (5.95,0);
\node at (0.75, -0.25) {$^{\phantom{-1}}x^{\phantom{-1}}$};
\node at (2.25, -0.25) {$^{\phantom{-1}}y^{\phantom{-1}}$};
\node at (3.75, -0.25) {$x^{{-1}}$};
\node at (5.25, -0.25) {$y^{{-1}}$};

\node at (0.75, 0.25) {$+\alpha_1(x)$};
\node at (2.25, 0.25) {$+\alpha_1(y)$};
\node at (3.75, 0.25) {$-\alpha_1(x)$};
\node at (5.25, 0.25) {$-\alpha_1(y)$};

\draw[->] (1.45,0.1) -- (1.45,.6) -- (6,.6);
\draw[->, densely dotted] (3.05,0.1) -- (3.05,0.5) -- (6,.5);
\draw[->] (2.95,0.1) -- (2.95,0.85) -- (6,.85);
\draw[->, densely dotted] (4.55,0.1) -- (4.55,0.75) -- (6,.75);

\node at (2.25,0.6){$\bullet$};
\node at (2.25,1) {$+\alpha_2(y)$};
\node at (3.75,0.85){$\bullet$};
\node at (3.75,0.6){$\circ$};
\node at (3.75,0.5){$\bullet$};
\node at (3.75,0.85){$\bullet$};
\node at (3.75,1.25) {$-\alpha_2(x)$};
\node at (5.25,0.6){$\circ$};
\node at (5.25,0.5){$\bullet$};
\node at (5.25,0.85){$\circ$};
\node at (5.25,0.75){$\bullet$};
\end{tikzpicture}
    \end{aligned}
    \]
    More generally, let $\gamma_n F_S$ be the $n$-th step of the lower central series of $F_S$, recalled in \S\ref{sec:group ring prelims}.
    
    \begin{prop}[Invariants of commutators]\label{prop:invariants of commutators}
    For a nested commutator $w\in \gamma_n F_S$ and any generator $s\in S$, the commutator $[w,s]$ has invariants
    \begin{equation}\label{eq:invariants of commutators}
        \ell_{(\alpha_0|\alpha_1|\ldots|\alpha_n)}([w,s]) = \ell_{(\alpha_0|\ldots|\alpha_{n-1})}(w)\cdot\alpha_n(s) - \alpha_0(s) \cdot \ell_{(\alpha_1|\ldots|\alpha_n)}(w).
    \end{equation}
\end{prop}
The proof is deferred to Remark \ref{rmk:invariant of commutator} below.

Weight reduction can be treated as a recursive procedure. The $k$-th step of the weight-reduction process shows that any tensor $(\alpha_1|\ldots|\alpha_n)\in Bar^0(\tilde{w};A)$ is cohomolgous to a linear combination of tensors of the form $(\alpha_1|\ldots|\alpha_{n-k}|t^d)$ for $d\leq k$. Suppose we already completed a weight-reduction calculation for $(\alpha_1|\ldots|\alpha_{n-k})$, then we can use the following.
\begin{lemma}\label{lem:recursive}
    If $T\in Bar^0(\tilde{w};A)$ has weight reduction $P(t)\in A[t]$, then $T|t^d$ has weight reduction $P(t)\cdot t^d$.
\end{lemma}
\begin{proof}
 If $T = P(t) + d_{Bar}(\eta)$ for some polynomial $P(t)\in {A[t]}$ and for $\eta\in Bar^{-1}(w;A)$, then we claim that
    \[  T|\underbrace{t|\ldots|t}_{d \text{ times}} = P(t)\cdot t^d + d_{Bar}(\eta|\underbrace{t|\ldots|t}_{d \text{ times}}).
    \]
    Indeed, since every $f\in \ol{C}^0(\tilde{w};A)$ satisfies $f(0)=0$, it follows that $f\smile t = 0$. Thus,
    \[
    d_{Bar}(\eta|\underbrace{t|\ldots|t}_{d \text{ times}}) = d_{Bar}(\eta)|\underbrace{t|\ldots|t}_{d \text{ times}}
    \]
    demonstrating the claimed equality.
\end{proof}
One particular choice of weight reduction gives rise to the following discrete analog of Chen's iterated integrals in \cite{chen1971pi1}.

\begin{cor}[Letter braiding as iterated sum]\label{cor:iterated sum}
    For functions $\alpha_1,\ldots,\alpha_r \in A^S$, the letter-braiding number of a word $\tilde{w}\maps [n] \to S\cup S^{-1}$ is
    \begin{equation}
        \ell_{(\alpha_1|\ldots|\alpha_r)}(w) = \sum_{i_1<_w i_2<_w\ldots<_w i_r} \alpha_1(\tilde{w}(i_1))\ldots\alpha_r(\tilde{w}(i_r))
    \end{equation}
    where $i<_w j$ is the relations $i<j$ when $\tilde{w}(i)\in S$, and $i\leq j$ when $\tilde{w}(i)\in S^{-1}$.
\end{cor}
In particular, when the $\alpha_j$'s are indicator forms, this sum counts the (signed) number of occurrences of a sequence of generators and their inverses, in order, within a word. Both the signs and the ``edge cases" specified by $<_w$ are critical for invariance. For example, in $F_{\{x,y,z\}}$ the invariant $\ell_{X|Y|Z}(w)$ counts the number of times the generator $x^{\pm}$ appears before $y^{\pm}$ which in turn appears before $z^{\pm}$ in the word $w$ (see Example \ref{ex:intro-free group}).
\begin{proof}
    First, observe that for any $\alpha\in A^S$, writing  $\tilde{w}^*\alpha = fdx$ we have
    \[
    \tilde{w}^*(\alpha)_i = \alpha(\tilde{w}(i)) = f(i).
    \]
    The result thus follows from the next characterization of weight-reduced polynomials: let $\sim$ denote the relation of differing by coboundary, then for all functions $f_1,\ldots,f_r$ we have $(f_1 dx|\ldots|f_r dx ) \sim \ell\cdot t + t^2 P(t)$ for some polynomial $P(t)$ and
    \begin{equation}\label{eq:iterated sum}
        \ell = \sum_{\substack{1\leq i_1\leq  i_2\leq \ldots\leq i_r \leq n\\ \text{s.t. }\forall j:\,
        \partial_1(i_j) \leq i_{j+1}}} f_1(i_1)f_2(i_2)\ldots f_r(i_r),
    \end{equation}
    where $\partial_1(i)$ is the terminal vertex of the $i$-th edge in $S^1(\tilde{w})$. Indeed, the condition $\partial_1(i)\leq j$ is equivalent to $i<_w j$.

    The equality is proved by induction on $r$, with $r=1$ being the equality in Example \ref{ex:linear invariants}. Now, denote the iterated sum in \eqref{eq:iterated sum} by $\int_{w} (f_1dx|\ldots|f_r dx)$ and assume by induction that for all $(f_1,\ldots,f_r)$
    \[
    (f_1 dx|\ldots|f_r dx) \sim t^2\cdot P(t) + t\cdot \int_{w} (f_1dx|\ldots|f_r dx).
    \]
    Then for $r+1$ we have
    %by \eqref{eq:weight-reduced expression}
    \[
    f_1 dx|\ldots|f_r dx | f_{r+1}dx \sim \left(\int_{w}f_{r+1}dx \right)\cdot (f_1 dx|\ldots|f_r dx |t) - (f_1 dx|\ldots|f_r dx\smile d^{-1}( f_{r+1}dx ) ).
    \]
    Lemma \ref{lem:recursive} shows that the first term has polynomial invariant divisible by $t^2$ and thus does not contribute to the linear coefficient. And by the cup product formulas in Proposition \ref{prop:weight-reduction},
    \[
    f_r dx \smile d^{-1}(f_{r+1}dx) = -\bigg( \sum _{i_{r+1} \geq \partial_1(\bullet)} f_{r+1}(i_{r+1})\bigg)f_r(\bullet) dx.
    \]
    The inductive hypothesis now implies that the linear coefficient is
    \[
    \sum_{\substack{1\leq i_1\leq \ldots\leq i_r \leq n \\
    \forall j:\, \partial_1(i_j) \leq i_{j+1}}} f_1(i_1)\ldots f_{r-1}(i_{r-1}) \cdot \bigg(f_r(i_r)\cdot \sum_{\substack{i_{r}\leq i_{r+1} \\ \partial_1(i_{r}) \leq i_{r+1}}} f_{r+1}(i_{r+1})\bigg). \qedhere
    \]
\end{proof}

\subsection{Properties of letter-braiding}
The following example specializes the iterated sum formula to the case of the universal 2-letter word.
\begin{example}[The universal product]\label{ex:universal product}
    Let $\{e_1,e_2\}$ generate a free group of rank $2$. We compute the letter-braiding invariants of the word $\tilde{w} = e_1e_2$. Denote the dual basis to $e_1$ and $e_2$ in $A^{\{e_1,e_2\}}$ by $\langle E_1, E_2 \rangle$, that is, $E_1(e_1) = E_2(e_2) = 1$ and $E_1(e_2)=E_2(e_1) = 0$. Then ${T}(A^{\{ e_1,e_2\}}) \cong A\langle E_1,E_2\rangle$ is the ring of noncommutative polynomials in $E_1$ and $E_2$.

    The simplicial map corresponding to the word $e_1e_2$ is
    \[
    \begin{split}
    \xymatrix{
    1 \ar[r]^{e_1}_1 & 2 \ar[r]^{e_2}_2 & 0 \ar@/^2pc/[ll]_0
    }    
    \end{split} \quad \xrightarrow{e_1e_2}\quad \underbrace{S^1}_{e_1} \vee \underbrace{S^1}_{e_2} \;=\; X_{\{e_1,e_2 \}}
    \]
    collapsing the $0$-th edge to the basepoint. Its pullback $(e_1e_2)^*\maps C^1(X_{\{e_1,e_2\}};A) \to C^1(e_1e_2;A)$ sends $E_1$ to the indicator form $\delta_1$ and $E_2$ to the indicator $\delta_2$.   
    \begin{prop}
        Let $T\in A\langle E_1,E_2\rangle$ be any noncommutative monomial in $E_1$ and $E_2$. The letter-braiding number of the word $E_1E_2\in F_{\{e_1,e_2\}}$ on $T$ is
        \[
        \ell_{T}(e_1e_2) = \begin{cases}
            1 & \text{ if } T = E_1, E_2 \text{ or }E_1|E_2 \\
            0 & \other.
        \end{cases}
        \]
    \end{prop}

\begin{proof}
    By the iterated sum formula in Corollary $\ref{cor:iterated sum}$,
    \[
    \ell_{(E_{j_1}|\ldots|E_{j_r})}(e_1e_2) = \sum_{1\leq i_1<\ldots< i_r \leq 2} E_{j_1}(e_{i_1})\ldots E_{j_r}(e_{i_n}).
    \]
    If $n>2$ then this sum is empty, and the invariant vanishes. When $n=2$ the invariant is $E_{j_1}(e_1)E_{j_2}(e_2)$, which vanishes unless $(E_{j_1}|E_{j_2}) = (E_1|E_2)$ whose invariant is $1$. The remaining case of $n=1$ is clear.
\end{proof}
\end{example}

We can now calculate the letter-braiding number of a product.
\begin{cor}\label{cor:invariant of products}
    If $w_1,w_2\in F_S$ are two elements, and $T\in \ol{T}(A^S)$ in the augmentation ideal satisfies $\ol\Delta(T) = \sum_i T_1^i\otimes T_2^i$, then
    \begin{equation}\label{eq:invariant of product}
        \ell_{T}(w_1\cdot w_2) = \ell_{T}(w_1) + \ell_{T}(w_2) + \sum_i \ell_{T_1^i}(w_1)\cdot \ell_{T_2^i}(w_2).
    \end{equation}
\end{cor}
\begin{proof}
    Define a group homomorphism $h\maps F_{\{e_1,e_2\}} \to F_S$ by $e_i\mapsto w_i$. Now, the composition
    \[
    \xymatrix{
    F_{\{s\}} \ar[r]^-{e_1\cdot e_2} & F_{\{e_1,e_2\}} \ar[r]^-{h} & F_S
    }
    \]
    sends $s$ to $w_1\cdot w_2$. Consider the functorial pullback maps defined in Theorem \ref{thm:naturality}.
    
    On the one hand, Example \ref{ex:homomorphism from cyclic group} shows that the pullback along the composition is the letter-braiding polynomial $L_{T}(w_1\cdot w_2)$.
    On the other hand, Example \ref{ex:pullback of product} gives
    \[
    h^*(T) = \ell_{T}(w_1)\cdot E_1 + \ell_{T}(w_2)\cdot E_2 + \sum_{i}\ell_{T_1^i}(w_1) \ell_{T_2^i}(w_2) \cdot E_1|E_2  + \text{ all other monomials}.
    \]
    Applying Example \ref{ex:homomorphism from cyclic group} again, now to $(e_1 e_2)^*$, gives the letter-braiding polyonmial of the word $e_1 e_2$ on any noncommutative monomial in $E_1$ and $E_2$. But these were computed in Example \ref{ex:universal product}, and the braiding numbers $\ell_{\bullet}(e_1 e_2)$ vanish on all monomials other than the ones listed in the last equation. Substituting the nonzero braiding numbers proves \eqref{eq:invariant of product}.
\end{proof}

\begin{example}
    The product formula \eqref{eq:invariant of product} gives a way to compute $\ell_{T}(s^{-1})$, for inverses of generators, distinct from Formula \eqref{eq:iterated sum}. Indeed, if $T = (\alpha_1|\ldots|\alpha_r)$ then
    \[
    0 = \ell_{T}(ss^{-1}) = \ell_{T}(s) + \ell_{T}(s^{-1}) + \sum_{i=1}^{r-1}\ell_{(\alpha_1|\ldots|\alpha_i)}(s)\ell_{(\alpha_{i+1}|\ldots|\alpha_r)}(s^{-1})
    \]
    and since Example \ref{ex:braiding of a generator} shows that $\ell_{(\alpha_1|\ldots|\alpha_i)}(s)$ vanishes unless $i=1$, we get
    \[
    \ell_{(\alpha_1|\ldots|\alpha_r)}(s^{-1}) = -\alpha_1(s)\ell_{(\alpha_2|\ldots|\alpha_r)}(s^{-1}).
    \]
    It thus follows by induction that
    \begin{equation}
        \ell_{(\alpha_1|\ldots|\alpha_r)}(s^{-1}) = (-1)^n \alpha_1(s)\ldots\alpha_r(s).
    \end{equation}

    From this, one can check more generally that for any $w\in F_S$, its inverse has
    \begin{equation}
        \ell_{T}(w^{-1}) = L_T(w)|_{t=-1},
    \end{equation}
    the evaluation of the braiding polynomial for $w$ at $(-1)$.
\end{example}
    We will see in \S\ref{sec:group rings} that the product formula \eqref{eq:invariant of product} is equivalent to the statement that the (unreduced) coproduct $\Delta$ on $\H^0_{Bar}(F_S;A)$ is dual to the convolution product of the group ring $A[F_S]$.

\begin{cor}\label{cor:degree bound}
    Let $\gamma_k F_S$ be the $k$-th term in the lower central series. If $w\in \gamma_k F_S$, and $T\in \ol{T}(A^S)$ is a reduced tensor of weight $r$, then the letter-braiding polynomial satisfies
    \[
    \deg L_{T}(w) \leq \frac{r}{k}.
    \]
    In particular,
    \[
    L_{T}(w) = \begin{cases}
        \ell_T(w)\cdot t &\text{ if } r<2k \\
        0 & \text{ if } r<k.
    \end{cases}
    \]
\end{cor}
\begin{proof}
%\cite[Theorem 5.1]{FoxI} proves that the recursion in \ref{cor:invariant of products} implies that $\ell_w(T)=0$ whenever $w\in \gamma_{k}F_S$ and $T$ has weight $n<k$ (Fox's proof only considered functionals defined over $\Z$, but the same argument applies to general commutative rings).
For $T$ of weight $1$ the claim is obvious, since all $\alpha\in \H^1(F_S;A)$ vanish on $[F_S,F_S]$.

Now, assume by induction on $n$ that any tensor $T'$ of weight $r'< r$ has $\deg L_{T'}(w) \leq r'/k$, so in particular if $L_{T'}(w)\neq 0$ then $r'\geq k$. Now, if $T$ has weight $r$ and $\deg L_T(w) = d$, then $\ol\Delta^{d}L_T(w) \neq 0$. But since $L_\bullet(w)$ is a coalgebra homomorphism, this means that 
\[
\ol\Delta^d(T) = \sum_i T_1^i\otimes \ldots \otimes T_d^i \implies 0\neq \ol\Delta^d(L_T(w)) = \sum_i L_{T_1^i}(w)\otimes \ldots \otimes L_{T_d^i}(w).
\]
Fix a summand $T_1^i\otimes \ldots \otimes T_d^i$ on which $L_\bullet(w)$ does not vanish. Since all $T_j^i$ are tensors of weight $< r$, the induction hypothesis implies that they must each have weight $\geq k$. It follows that $r\geq d\cdot k$.
\end{proof}
\begin{cor}
    If $T\in \ol{T}(A^S)$ is a reduced tensor of weight $< 2k$, then the function $w\mapsto \ell_T(w)$ is a group homomorphism $\gamma_k F_S \to A$. %Moreover, when $F_S$ is finitely generated and $A$ is a $\Q$-algebra, every homomorphism $\sfrac{\gamma_k F_S}{\gamma_{2k}F_S} \to A$ is obtained in this way\footnote{Over $\Z$, tensors $T$ only realize a full-rank sublattice inside $\Hom(\sfrac{\gamma_k F}{\gamma_{2k}F}, \Z)$.}.
\end{cor}
\begin{proof}
    Using the product formula $\ell_{T}(w_1\cdot w_2) = \ell_{T}(w_1) + \ell_{T}(w_2) + \sum_i\ell_{T_1^i}(w_1)\cdot \ell_{T_2^i}(w_2)$, it suffices to show that $\ell_{T_1^i}(w_1)\ell_{T_2^i}(w_2)=0$ for every $i$. This follows since their weights satisfy $\wt(T_1)+  \wt(T_2) \leq \wt(T) < 2k$, so one of the weights must be smaller than $k$ and thus vanish on $\gamma_k F_S$.
    %To see that any homomorphism is realized by some $T$, filter by abelian groups
    %\[
    %0 \leq \sfrac{\gamma_{2k-1}F}{\gamma_{2k}F} \leq \ldots \leq \sfrac{\gamma_{k+i}F}{\gamma_{2k}F} \leq \ldots \leq \sfrac{\gamma_{k}F}{\gamma_{2k}F}
    %\]
    %whose $i$-th successive quotient is $\sfrac{\gamma_{k+i}F}{\gamma_{k+i+1}F}$. \cite[Theorem 5.1]{monroe-sinha} shows that every group homomorphism $\sfrac{\gamma_{k+i}F}{\gamma_{k+i+1}F}\to \Q$ is realized by $\ell_\bullet(T_{k+i})$ for some $T_{k+i}\in \ol{T}^{k+i}(\Q^S)$, which automatically vanishes on $\gamma_{k+i+1}F$. This fact generalizes to any $\Q$-algebra $A$, since the domain is finitely generated. Therefore, start by realizing $h$ on $\gamma_{2k-1}F$ by some $T_{2k-1}$. Working recursively, having picked $T_{k+i+1}+\ldots+T_{2k-1}$, it defines a homomorphism $h'\maps \gamma_{k+i}F\to A$, and we pick $T_{k+i}$ so as to realize $h-h'$. Since $T_{k+i}$ vanishes on deeper commutators, adding it to the existing sum corrects the homomorphism on $\gamma_{k+i} F$ without altering its values on deeper commutators. Continuing this process, $T_k+\ldots+ T_{2k-1}$ defines $h$ on all of $\gamma_k F$.
\end{proof}
In fact, over $\Q$, the results in \cite{monroe-sinha} imply that every homomorphism $\gamma_k F_S/\gamma_{2k}F_S \to \Q$ is obtained as the letter-braiding invariant for a unique $T\in \ol{T}(\Q^S)$ of weight $< 2k$.

\section{General groups}

\label{sec:general groups}
We extend  letter-braiding invariants to general groups, first recalling definitions. 
Define the fundamental group of a pointed simplicial set $(Y_\bullet,*)$ as that of its geometric realization. Let $\pi_1(Y,*)=\Gamma$. Letter-braiding invariants of $\Gamma$ are assigned to the coalgebra $\H^0_{Bar}(Y;A)$. Every $w \in \Gamma$ determines a continuous map $w\maps |S^1|\to |Y|$, well-defined up to homotopy, and thus induces a pullback
\[
w^*\maps \H^0_{Bar}(Y;A) \to \H^0_{Bar}(S^1;A).
\]
Using the canonical isomorphism $L\maps \H^0_{Bar}(S^1;A) \to A[t]$, we define the following.
\begin{definition}\label{def:braiding polynomial for general group}
    Fix $T \in \H^0_{Bar}(Y;A)$. For every $w\in \Gamma$  the \emph{letter-braiding polynomial} of $w$ associated to $T$ is
    \[
    L_T(w) = L(w^*(T)) \in  A[t].
    \]
The \emph{letter-braiding number} $\ell_T(w)\in A$ is the linear coefficient of $L_T(w)$.
\end{definition}
By relating these letter-braiding polynomials to the ones in \S\ref{sec:free groups}, we will get formulas for computing them and obstructions for their existence.

The following claim shows that the particular space that realizes a group $\Gamma$ as its fundamental group is immaterial. Consequently, any such space $(Y,*)$ gives a distinct viewpoint on the braiding invariants of words in $\Gamma$, in some cases by a geometric description of the ring $\H^*(Y;A)$.
\begin{lemma}\label{lem:independence of model}
    If $(Y,*)$ and $(Y',*')$ are connected pointed simplicial sets and $f\maps (|Y|,*)\to (|Y'|,*')$ is a pointed map of their realizations that induces an isomorphism on $\pi_1$, then the pullback
    \[
    f^* \maps \H^0_{Bar}(Y';A)\to \H^0_{Bar}(Y;A)
    \]
     given by Lemma \ref{lem:induced maps on bar cohomology} is an isomorphism.
\end{lemma}

\begin{proof}
    Denote $\Gamma = \pi_1(Y',*')$ and consider the map $i:Y'\to K(\Gamma,1)$ into an aspherical simplicial set inducing an isomorphism on fundamental groups. The claim will follow if we prove that the composition $(i\circ f):Y\to K(\Gamma,1)$ induces an isomorphism on $H^0_{Bar}(-;A)$ since then both maps
    \[
    H^0_{Bar}(Y';A) \xleftarrow{(i\circ \Id)^*} H^0_{Bar}(K(\Gamma,1);A) \xrightarrow{(i\circ f)^*} H^0_{Bar}(Y;A)
    \]
    will have been shown to be isomorphisms and the right-pointing map factors as $f^*\circ i^*$.

    Consider the homotopy fiber $F\to Y\to K(\Gamma,1)$. By the long exact sequence in homotopy,
    \[
    \cancel{\pi_2(K(\Gamma,1))}\to \pi_1(F) \to \pi_1(Y)\xrightarrow{\sim} \pi_1(K(\Gamma,1)) \to \pi_0(F)\to \pi_0(Y)\xrightarrow{\sim} \pi_0(K(\Gamma,1))
    \]
    it follows that $F$ is simply-connected and, in particular $H^1(F;A)=0$. The Serre spectral sequence of the fibration defining $F$ implies the five-term exact sequence in cohomology with coefficients in $A$
    \[
    0\to H^1(K(\Gamma,1)) \to H^1(Y)\to\ \cancel{H^0(K(\Gamma,1);H^1(F))} \to H^2(K(\Gamma,1)) \to H^2(Y),
    \]
    using $H^0(F;A)=A$. This sequence shows that the induced restriction $(i\circ f)^*:H^*(K(\Gamma;1))\to H^*(Y)$ is an isomorphism for $*=1$ and an injection for $*=2$.
    At the level of cochain complexes this means that the pullback $(i\circ f)^*\maps C^*(K(\Gamma,1);A)\to C^*(Y;A)$ is $1$-connected, so that after desuspension it becomes $0$-connected on $\desus{\ol{C}(-;A)}$.

    Since the complexes here consist of torsion-free and thus flat $A$-modules, it follows that on every tensor power 
    \[
    ((i\circ f)^*)^{\otimes p}\maps \desus{\ol{C}(K(\Gamma,1);A)}^{\otimes p} \to \desus{\ol{C}(Y;A)}^{\otimes p}
    \]
    is also $0$-connected.
    Passing to bar constructions, we imitate the proof that $Bar(-)$ is a quasi-isomorphism invariant (Lemma \ref{lem:quasi-iso invariant}) and consider the spectral sequence of the filtration by weight. Its $E_0$ page is simply the direct sum of complexes $\desus{\ol{C}(-;A)}^{\otimes p}$, and thus $i\circ f$ induces a $0$-connected map on the $E_0$-page. Since this connectivity property passes through the pages of the spectral sequence, it follows that the map 
    \[
    (i\circ f)^*\maps Bar(K(\Gamma,1);A) \to Bar(Y;A)
    \]
    is already $0$-connected. It therefore induces an isomorphism on $\H^0_{Bar}(-;A)$ and an injection on $\H^1_{Bar}(-;A)$.
\end{proof}
In the present work, our preferred choice of simplicial set with fundamental group $\Gamma$ is the standard simplicial model for the classifying space.
\begin{definition}
    Recall, the simplicial bar construction of a group $\Gamma$ is the simpilicial set 
    \[
    B\Gamma_\bullet\maps [p] \mapsto \Gamma^p
    \]
    with face maps $\partial_0$ and $\partial_p$ forgetting the respective extremal coordinate, internal face maps given by multiplication in $\Gamma$ of adjacent terms, and degeneracies inserting copies of the unit. Its fundamental group is canonically isomorphic to $\Gamma$.    
    Denote
    \[
    Bar(\Gamma;A) := Bar(B\Gamma;A), \text{ and } \H^*_{Bar}(\Gamma;A) = \H^*(Bar(B\Gamma;A)).
    \]
\end{definition}
The functoriality of the construction $\Gamma \mapsto B\Gamma$ implies the following.
\begin{lemma}\label{lem:Bar(G) is functorial}
    The assignment $\Gamma \mapsto Bar(\Gamma;A)$ is a contravariant functor from the category of groups to the category of bigraded dg-coalgebras. Thus $ \Gamma \mapsto \H^0_{Bar}(\Gamma;A)
    $
    is a contravariant functor taking values in the category of filtered coalgebras.

Moreover, any connected $(Y,*)$ admits a natural isomorphism
\[\H^0_{Bar}(\pi_1(Y,*);A) \overset{\sim}\to  \H^0_{Bar}(Y;A).
\]
\end{lemma}
\begin{proof}
    Recalling that $|B\Gamma_\bullet|$ is a $K(\Gamma,1)$ space, homotopy classes of pointed maps $(X,*)\to |B\Gamma_\bullet|$ are in bijection with group homomorphisms $\pi_1(X,*)\to \Gamma$. In particular, for any topological space $(X,*)$ there exists a pointed map $i\maps X \to |B(\pi_1(X,*))_\bullet|$, unique up to homotopy, inducing the identity on $\pi_1$. Applied to $X = |Y|$, Lemma \ref{lem:induced maps on bar cohomology} gives a canonical pullback
 \[
 i^*\maps \H^0_{Bar}(\pi_1(Y,*);A) \to  \H^0_{Bar}(Y;A),
 \]
 which by Lemma \ref{lem:independence of model} is an isomorphism. Functoriality and naturality are clear.
\end{proof}
\begin{comment}
\begin{example}
    Given a set $S$, the simplicial set $X_S$ from \ref{notation:Xs} has free fundamental group $F_S$, so that the Lemma \ref{lem:Bar(G) is functorial} gives a canonical isomorphism
    \[
    \H^0_{Bar}(F_S;A) \cong \H^0_{Bar}(X_S;A) = {T}(A^S).
    \]
\end{example}
\end{comment}

\subsection{Letter braiding and vanishing (Massey) products}\label{sec:invariants from Massey}

The pages of the spectral sequence of the bicomplex $Bar(\Gamma;A)$ give approximations to $H^0_{Bar}(\Gamma; A)$. The first page the tensor (co)algebra on $H^1(\Gamma; A)=\Hom_{grp}(\Gamma,A)$, which is already the full answer when $\Gamma$ is free as we saw in \S\ref{sec:free groups}. The second page of the sequence in degree $0$ is the $0$-th bar cohomology of the graded ring $H^*(B\Gamma; A)$.  We see immediately
that, in weight two, $H^0_{Bar}(\Gamma; A)$ is
computed by the kernel of the cup product 
$H^1(B\Gamma; A) \otimes H^1(B\Gamma; A) 
\to H^2(B\Gamma; A)$

In general constructing letter-braiding invariants for $\Gamma$ relies on understanding the cohomology ring, including its Massey products, of $B\Gamma$ or any other $X$ with $\pi_1(X,*)\cong \Gamma$ in 
degrees one and two. We will observe that introducing relations to a group simultaneously gives obstructions for braiding invariants to be well-defined as well as creates Massey products in this cohomology ring, both encoded by $d_{Bar}\maps Bar^0\to Bar^1$.

\begin{example} \label{ex:relation1}
    Let $F_{\{e,f\}}$ be the free group on two generators, with letter-braiding invariants indexed by noncommutative polynomials $\Z\langle E,F\rangle$. For example, $\ell_E(-)$ counts the signed number of $e$'s in words.
    
    Now, introduce the relation $e=f$, so a letter $e$ may be replaced by an $f$. With this, the functional $\ell_E(-)$ is no longer well-defined, e.g. $eefe = fffe$ under this relation, but
    \[
    \ell_{E}(eefe) = 3 \neq 1 = \ell_{E}(fffe).
    \]
    So $E$ and $F$ are lost as invariants of the group. However, $E+F$ still defines an invariant.

    At the level of the bar complex, the effect of introducing the relation is to add a $2$-cell $D^2_{e^{-1}f}$, so that $dE = -dF = [D^2_{e^{-1}f}]$ are no longer closed. But their sum $E+F$ is closed, and represents a braiding invariant.
\end{example}
\begin{example}\label{ex:relation2}\label{ex:abelian group}
    A more interesting example is obtained by introducing the relation $ef=fe$ to $F_{\{e,f\}}$, modeled by a $2$-torus. Now, $E$ and $F$ still define invariants of words, but $(E|F)$ does not: it counts occurrences of $e^{\pm}$ \emph{before} $f^{\pm}$, which changes upon application of the relation.
    
    Accordingly, $(E|F)$ is not a bar cocycle, because $E\smile F \neq 0$ in cohomology. Anti-symmetry of the cup product means that $(E|F) + (F|E)$ is an invariant.
\end{example}
\begin{remark}
    The loss of the invariant $(E|F)$ in the last example is reflected by the fact that in the abelian group $F_{\{e,f\}}/\langle [e,f] \rangle$, any two distinct elements are already distinguished by their evaluation on $E$ and $F$, so no additional invariants are needed. Still, higher weight invariants exist and serve to distinguish linear combinations of elements in the group ring, as discussed in Example \ref{ex:abelian group ring}.
\end{remark}

On the cohomological side, the bar cohomology features cochains with vanishing products, as explained next.
\begin{example}\label{ex:relation3}
    Let $Y$ be any pointed simplicial set. Given $1$-cochains $\alpha_1,\ldots,\alpha_n \in C^1(Y;A)$, let's examine when the tensor $T=(\alpha_1|\ldots|\alpha_n)\in Bar^0(Y;A)$ is a cocycle. The differential in the bar construction is $d_{Bar} = d_1+d_2$, with $d_1$ internal to $C^*(Y;A)$ and $d_2$ given by cup products of cochains. So $d_{Bar}T=0$ iff 
\begin{itemize}
    \item every $\alpha_i$ is closed in $C^1(Y;A)$, and
    \item every cup product $\alpha_i\smile \alpha_{i+1}$ vanishes in $C^2(Y;A)$.
\end{itemize}
When these conditions are satisfied, $T$ represents a class in $\H^0_{Bar}(Y;A)$.

 Now, let $w\maps S^1(w)\to Y$ be a simplicial map representing some $w\in \pi_1(Y,*)$. Then 
 \[
 w^*(T) = (w^*(\alpha_1)|\ldots|w^*(\alpha_n)) \in Bar^0(w;A)
 \]
 can be weight-reduced to a polynomial $L_T(w)\in {A[t]}$, thus giving the letter-braiding polynomial of $w$.
 Most importantly, due to relations in $\pi_1(Y,*)$ there might be different words representing the same group element. But, Lemma \ref{lem:induced maps on bar cohomology} guarantees that the letter-braiding polynomial only depends on the pointed homotopy class $S^1\to Y$, and is thus an invariant of $w$.
\end{example}

To better understand the condition $d_{Bar}(T) = 0$, consider the following two examples.
\begin{example}\label{ex:relation4}
    Suppose $\alpha,\beta\in C^1(Y;A)$ are $1$-cocycles whose cup product vanishes in cohomology. Then $d_{Bar}(\alpha|\beta) = -\alpha\smile\beta = d\eta$ for some $\eta\in C^1(Y;A)$ but might not be $0$. However, the in-homogeneous linear combination $(\alpha|\beta) + \eta$ satisfies,
    \[
    d_{Bar}\maps (\alpha|\beta) + \eta \longmapsto -\alpha\smile \beta - d\eta = 0
    \]
    giving a Bar $0$-cocycle of weight $2$. We thus get a well-defined letter-braiding invariant, by weight reducing $(w^*(\alpha)|w^*(\beta)) + w^*(\eta) \in Bar^0(w;A)$, while neither of the summands is by itself a homotopy invariant.
    
    Summarising, cup products need only vanish at the level of cohomology to give invariants, but one must then include a correction term of lower weight `witnessing' this.
\end{example}
\begin{example}\label{ex:relation5}
    To construct a weight $3$ bar cocycle, consider closed $1$-cochains $\alpha,\beta,\gamma\in C^1(Y;A)$ such that both $\alpha\smile \beta = d\eta$ and $\beta\smile \gamma = d\zeta$ vanish in cohomology. Then
    $$d_{Bar}((\alpha|\beta|\gamma) - (\eta|\gamma) - (\alpha|\zeta)) = ((\eta\smile\gamma) + (\alpha\smile \zeta))$$ might still be nonzero. This last expression represents the \emph{triple Massey product} $\langle \alpha,\beta,\gamma\rangle$, and it is a $2$-cocycle in $C^2(Y;A)$. If it vanishes in cohomology, then a secondary correction term $\tau$ exists, making the inhomogeneous sum $((\alpha|\beta|\gamma) - (\eta|\gamma) - (\alpha|\zeta) + \tau)$ into a bar cocycle and defining a letter-braiding invariant.
\end{example}
More generally, tensors $(\alpha_1|\ldots|\alpha_n)$ extend to bar cocycles iff all their consecutive cup and higher Massey products vanish in cohomology.
For general linear combinations of tensors, this is all best housed in the spectral sequence we discussed above, given formally as follows.

%Making this precise, recall that $Bar^*(Y;A)$ is bigraded and consider its filtration by columns. Then its associated graded coalgebra $\operatorname{gr} Bar^*(Y;A)$ has the simpler differential $d_1$ -- the tensor product differential of $\bigoplus_p C^*(Y;A)^{\otimes p}$ up to sign.
\begin{lemma}[Weight spectral sequence]\label{lem:Bar spectral sequence}
    Suppose $(Y,*)$ is connected. There is a natural second quadrant spectral sequence converging to $\H^*_{Bar}(Y;A)$, which in total degree $0$ admits a natural inclusion of graded coalgebras 
    \[
    E_r^{-\bullet,\bullet}\into  T(\H^1(Y;A)) \quad\text{where}\quad E_r^{-p,p} \into \H^1(Y;A)^{\otimes p}
    \]
    on every page. On the $E_1$-page this inclusion is an equality, and the differential given by the cup product in cohomology
    \[
    [\alpha_1]\otimes \ldots\otimes [\alpha_p] \mapsto -\sum_{i=1}^{p-1} [\alpha_1]\otimes\ldots\otimes [\alpha_i\smile \alpha_{i+1}] \otimes\ldots\otimes[\alpha_p].
    \]
    Later differentials are given by higher Massey products, so that cycles in $E_r^{-p,p}$ are tensors of weight $p$ in $\H^1(Y;A)$ with vanishing Massey products of size $\leq r$.
    For example\footnote{In practice, it is not likely that pure tensors have well-defined Massey products, and only linear combinations of them survive to the $E_r$-page of the spectral sequence. However, see Proposition \ref{prop:global Massey}.},
    \begin{equation}\label{eq:higher differentials}
    d_r\maps [\alpha_1]\otimes \ldots\otimes [\alpha_p] \mapsto - \sum_{i=1}^{p-r} [\alpha_1]\otimes\ldots\otimes \langle \alpha_i,\ldots, \alpha_{i+r}\rangle \otimes\ldots\otimes[\alpha_p],
    \end{equation}
    
    This spectral sequence is natural in pointed simplicial maps $(Y,*)\to (Y',*')$.
\end{lemma}
\begin{proof}
    The spectral sequence comes from the filtration of $Bar^*(Y;A)$ by weight, which is clearly natural in maps of pointed simplicial sets. The $E_0$-page of this spectral sequence consists of columns $\ol{C}^*(Y;A)^{\otimes p}$ shifted down $p$ degrees, with the usual tensor product differential. Since $\tilde\H^0(Y;A)=0$, the K\"unneth spectral sequence shows that in degree $<p$ this tensor product complex has vanishing cohomology, and in degree $p$ its cohomology is $\cong \H^1(Y;A)^{\otimes p}$.
    
    It follows that $E_1^{-p,q}$ vanishes in negative total degree, and $E_1^{-p,p} \cong \H^1(Y;A)^{\otimes p}$, as claimed. In particular, no differentials land in total degree $0$, and thus all pages of the spectral sequence have
    \[
    E_r^{-p,p}\into E_1^{-p,p}
    \]
    characterized by the vanishing of the first $r$ differentials. The cup product on cohomology and higher Massey products coincide with the differentials of this spectral sequence essentially by their definition, see \cite{may-matric-massey}.
\end{proof}
\begin{cor} \label{cor:leading terms}
    For every weight $n\geq 1$ there is a natural `leading term' map
    \[
    \lead \maps \H^0_{Bar}(\Gamma;A)_{\leq n} \to \H^1(\Gamma;A)^{\otimes n}
    \]
    with image all tensors with vanishing cup and Massey products, and kernel $\H^0_{Bar}(\Gamma;A)_{\leq n-1}$.
\end{cor}
%\begin{remark}[\textbf{Globally defined Massey product}] \label{rmk:global massey}

    In many contexts, the celebrated Homotopy Transfer Theorem \cite{kadeishvili} gives globally defined Massey products
    \begin{equation}\label{eq:global massey}
    \mu_r \maps \H^1(Y;A)^{\otimes r} \to \H^2(Y;A) \quad \forall r\geq 2
    \end{equation}
    such that the $r$-th differential of the above weight spectral sequence agrees with $\mu_r$ as written in \eqref{eq:higher differentials}. For example, this is always possible when $A$ is a field. More generally we have the following.
\begin{prop}\label{prop:global Massey}
    Let $(Y,*)$ be a finite $2$-complex with a single $0$-cell $*$, so that $\H^1(Y;A) \into C^1(Y;A)$ is the submodule of closed cochains. If this inclusion is split, then global Massey products in \eqref{eq:global massey} are defined and compute the differentials of the weight spectral sequence via \eqref{eq:higher differentials}.

    When $A=\Z$, the same holds even when $(Y,*)$ is countable and $\H^1(Y;A)$ and $\H^2(Y;A)$ are countably generated. When $A$ is a field, global Massey products exist with no hypotheses.
\end{prop}
\begin{proof}
    Under the hypotheses, the complex $C^*(Y;A)$ is a complex of free $A$-modules $C^0 \xrightarrow{0} C^1 \xrightarrow{d} C^2$ with $C^0 = A\cdot 1$. There is already a projection $q\maps C^2 \onto \coker(d) \cong \H^2$, so that a splitting $p\maps C^1 \to \H^1$ defines a quasi-isomorphism
    \[
    \xymatrix{
    \big(A\cdot 1 \ar[r]^0 \ar@{=}[d] & \phantom{\big(}C^1 \phantom{\big)} \ar[r]^d \ar[d]^p & C^2\big) \ar[d]^q \\
    \big(A\cdot 1 \ar[r]^0 & \phantom{\big(}\H^1 \phantom{\big)} \ar[r]^0 & \H^2\big).
    }
    \]
    Since $C^*$ is a bounded-below complex of projectives, Theorem 2 of  \cite{petersen-homotopy-transfer} and its following paragraph show that the cup product structure on $C^*$ can be transferred to a square-zero coderivation on $Bar(\H^*(Y;A))$ making the induced map
    \[
    Bar(C^*(Y;A)) \to Bar(\H^*(Y;A))
    \]
    into a quasi-isomorphism. Such a coderivation in degree $0$ has components $\mu_r \maps \H^1(Y;A)^{\otimes r} \to \H^2(Y;A)$, extended to $Bar^0(\H^*(Y;A))$ via the Leibnitz rule.
    Now, the two weight spectral sequences are equal from their $E_1$-pages and on, but with differentials of $Bar(\H^*(Y;A))$ restricted from \eqref{eq:higher differentials}, as claimed.
    
    When $A=\Z$, countable generation of $\H^1$ and $\H^2$ implies that there exists a quasi-isomorphic dg subalgebra $\tilde{C}^* \leq C^*$ that is countably generated. Since every $C^i \cong \Z^{\N}$ is a Baer--Specker group, Specker's theorem \cite{specker-countable-is-free} shows that the countably generated subgroup $\tilde{C}^i$ is free. From this point \cite[Theorem 2]{petersen-homotopy-transfer} applies, and the same argument proceeds as above.

    When $A$ is a field, every module is free and every injection is split, so the above argument applies without further hypotheses.
\end{proof}
%\end{remark}
%\begin{remark}
The weight spectral sequence gives an inductive process for calculating invariants of high weight: if one has calculated all weight $< n$ braiding invariants, then to get all invariants of weight $n$ it takes only to find those tensors in $\H^1(\Gamma;A)^{\otimes n}$ with vanishing cup and Massey products. Every such tensor extends to a weight $n$ braiding invariant, and while the extension is not unique, any two extensions differ by a lower weight invariant (which has already been found).    
%\end{remark}
\begin{comment}
\begin{example}
    Take $\Gamma = F_S$ a free group and consider the simplicial set $X_S$ presenting it. In $C^*(X_S;A)$, all cup products already vanish at the chain level, and therefore the spectral sequence mentioned above collapses: $\H^1(X_S;A)\cong A^S$ and
    \[
    E_1^{-p,p} = E_{\infty}^{-p,p} = \ol{T}^p(A^S).
    \]
    The filtration induced on $\H^0_{Bar}(F_S;A)$ is the usual weight filtration on the tensor coalgebra. Later on, Corollary \ref{cor:evaluation on high commutators} will show that this filtration is compatible with the lower central series filtration of $F_S$ in the sense that
    \[
    L_w\maps \H^0_{Bar}(F_S;A)_{\leq p-1} \to 0 
    \]
    whenever $w\in \gamma_p F_S$. Equivalently, $L_w$ is well-define on the filtered quotients $\H^0_{Bar}(F_S;A)_{> p-1}$.
\end{example}
\end{comment}
\begin{example}[\textbf{Surface groups}] \label{ex:surface grps}
    Let $\Sigma_g$ be a closed orientable surface of genus $g$, with fundamental group $\Gamma\cong \langle a_1,\ldots,a_g,b_1,\ldots b_g \mid [a_1,b_1]\cdots[a_g,b_g] \rangle$. Since $H^1(\Sigma_g;\Z)=\Z^{2g}$ with the symplectic pairing $\smile\maps \H^1\otimes \H^1\to \H^2\cong \Z$, letter-braiding invariants of words in the surface group, in low weights, are linear combinations of:
    \begin{itemize}
        \item $A_1,\ldots,A_g,B_1,\ldots,B_g \in H^1(\Sigma_g;\Z)$ counting signed occurrences of $a_i$ and $b_i$. These detect the image of a word in the abelianization.
        \item $A_i|A_j, B_i|B_j$ for all $(i,j)$ and $A_i|B_j, B_i|A_j$ for all $i\neq j$. These count pairs of generators in a particular order, e.g. $A_i|A_j$ for $i\neq j$ is the signed count of times $a_i$ appears \emph{before} $a_j$. They detect and distinguish between commutators $[a_i,a_j],[b_i,b_j]$ and $[a_i,b_j]$. The invariant $A_i|A_i$ coincides with the binomial $\binom{A_i}{2}$, thus it distinguishes between $2a_i$ and $a_i^2$ in the group ring, and similarly for $B_i|B_i$.
        \item $(A_i|B_i - A_1|B_1)$ and $(B_i|A_i + A_1|B_1)$ for all $i=1,\ldots,g$. These linear combinations of counts of ordered pairs of generators are unchanged by application of the surface relation, while the individual count of pairs $A_i|B_i$ are not invariant.
        \item $(A_i|B_i|A_i) - (A_1|B_1|A_i) - (A_i|B_1|A_1)$ for all $i = 1,\ldots,g$ are invariants on double-commutators, detecting e.g. that $[[a_i,b_i],a_i]\neq 1$. There are many other less interesting invariants in weight $3$, such as $(A_i|B_j|A_k)$ for $i\neq j \neq k$.
    \end{itemize}
    We give a full treatment modulo shuffles across all weights in \cite{joint}.
    \end{example}

\subsection{Properties of letter-braiding}\label{sec:properties general groups}
Formal properties of letter-braiding polynomials proved in \S\ref{sec:free groups} continue to hold for general groups. Most useful  is the following.
\begin{lemma}[\textbf{Push-pull formula}]
    Let $h\maps \Gamma_1\to \Gamma_2$ be a group homomorphism, then for every $w\in \Gamma_1$ and $T\in \H^0_{Bar}(\Gamma_2;A)$,
    \begin{equation}\label{eq:push-pull}
        L_{h^*(T)}(w) = L_{T}(h(w))
        \quad \text{ and } \quad \ell_{h^*(T)}(w) = \ell_{T}(h(w)).
    \end{equation}
\end{lemma}
\begin{proof}
    This is a direct application of Lemma \ref{lem:induced maps on bar cohomology} to the composition $S^1\xrightarrow{w} B\Gamma_1\xrightarrow{h} B\Gamma_2$.
\end{proof}

In particular, the push-pull formula immediately reveals behavior with respect to products.

\begin{cor}[\textbf{Product formula}]
    Suppose $w_1,w_2\in \Gamma$ are two elements, then
    \begin{equation*}%\label{eq:product formula general groups}
        \ell_\bullet({w_1\cdot w_2}) = \ell_\bullet({w_1}) + \ell_\bullet({w_2}) + (\ell_\bullet({w_1})\otimes \ell_\bullet({w_2}))\circ \ol{\Delta}.
    \end{equation*}
\end{cor}
\begin{proof}
    Let $h\maps F_{\{e,f\}}\to \Gamma$ send $h(e) = w_1$ and $h(f) = w_2$. By push-pull and Corollary \ref{cor:invariant of products},
    \[
    \ell_{T}(w_1w_2) = \ell_{h^*(T)}(ef) = \ell_{h^*(T)}(e) + \ell_{h^*(T)}(f) + (\ell_\bullet(e)\otimes \ell_\bullet(f))_{ \ol\Delta (h^*(T))}
    \]
    and since $h^*$ is a coalgebra homomorphism, $\ol\Delta(h^*(T)) = (h^*\otimes h^*)\circ \ol\Delta(T)$. Thus,
    \[
    \ell_{T}(w_1w_2) = \ell_{h^*(T)}(e) + \ell_{h^*(T)}(f) + (\ell_{h^*(\bullet)}(e) \otimes \ell_{h^*(\bullet)}(f))_{\ol\Delta(T)},
    \]
    which proves the claim by another application of the push-pull formula.
\end{proof}
\begin{cor}\label{cor:bounded degree general groups}
    If $w\in \Gamma$ belongs to the $k$-th term of the lower central series, $\gamma_k\Gamma$, and $T\in H^0(\Gamma;A)$ has weight $\leq r$, then 
    \[
    \deg L_T(w) \leq \frac{r}{k}.
    \]
    In particular, $\ell_T(w) = 0$ whenever $ r < k$.
\end{cor}
In the proof below, we make use of the following standard notation.
\begin{definition}
    A set of elements $s_1,\ldots,s_m\in \Gamma$ determines a group homomorphism $F_m\to \Gamma$ by $e_i \mapsto s_i$. Given a word $w\in F_m$, we denote its image in $\Gamma$ by $w(s_1,\ldots,s_m)$. An element of this form will be called \emph{a word} in $s_1,\ldots,s_m$.
\end{definition}

\begin{proof}[Proof of Corollary \ref{cor:bounded degree general groups}]
    Suppose $w$ can be written as a $k$-fold commutator in the elements $s_1,\ldots,s_m\in \Gamma$. That is, there exists a word $\tilde{w}\in \gamma_k F_m$ such that $\tilde{w}(s_1,\ldots,s_m) = w$. Denote $h\maps F_m\to \Gamma$ the homomorphism with $e_i\mapsto s_i$, so by push-pull
    \[
    L_{T}(w) = L_{h^*(T)}(\tilde{w}).
    \]
    But since $h^*$ preserves the weight filtration, the weight of $h^*(T)$ is also $\leq r$, and by \ref{cor:evaluation on high commutators} it follows that $\deg L_{h^*(T)}(\tilde{w})\leq r/k$.
\end{proof}

Since we already know how to calculate braiding numbers on free groups, it remains to understand pullbacks $h^*\maps \H^0_{Bar}(\Gamma;A)\to \H^0_{Bar}(F_S;A)$ for homomorphisms $h\maps F_S \to \Gamma$.

\begin{theorem}[\textbf{Composition formula}]\label{thm:pullback to free group}
    Suppose $h\maps F_S\to \Gamma$ is a group homomorphism from a free group. Then the induced coalgebra homomorphism
    \[
    h^* \maps\H^0_{Bar}(\Gamma;A) \to \H^0_{Bar}(F_S;A) \cong T(A^S)
    \]
    is determined by its weight $1$ component $\ell h^*\maps \H^0_{Bar}(\Gamma;A) \to A^S$, which is
    \begin{equation*}
        \ell h^*(T)\maps (s \mapsto \ell_{T}(h(s))).
    \end{equation*}
    It follows that the letter-braiding numbers of a word $w(s_1,\ldots,s_r)\in \Gamma$ are
    \begin{equation} \label{eq:composition formula}
    \ell_{T}(w(s_1,\ldots,s_r)) = \sum_{k\geq 1} \sum_{i_1,\ldots,i_k} \ell_{(S_{i_1}|\ldots|S_{i_k})}(w) \cdot (\ell_\bullet({h(s_{i_1})})\otimes \cdots \otimes \ell_\bullet({h(s_{i_k})}))_{\ol{\Delta}^{k-1} (T)},
    \end{equation}
    where for all $1\leq i \leq n$ the function $S_i\in A^{S}$ is the indicator function $S_i\maps s_j \mapsto \delta_{ij}$ for $j\in S$.
\end{theorem}
\begin{proof}
Theorem \ref{thm:naturality} above is the special case when $\Gamma$ is a free group, and the same proof holds in the present generality after replacing $X_{S_2}$ by $B\Gamma$.

Finally, the evaluation of $\ell_T$ on $w(s_1,\ldots,s_r)$ is the same as the evaluation of $h^*(T)$ on $w\in F_S$. Since $h^*$ is a coalgebra homomorphism, the coefficient of $(S_{i_1}|\ldots|S_{i_k})$ in $h^*(T)$ is the sum of products of the coefficients of the individual $S_{i_j}$'s appearing in the factors of $\overline{\Delta}^{k-1}(h^*(T))$. These are computed by the linear parts $\ell h^*(T)$ that we already understood in the previous paragraph and in Theorem \ref{thm:naturality}. Therefore, the coefficient of $(S_{i_1}|\ldots|S_{i_k})$ is given by
$\lambda_{(i_1,\ldots,i_k)} := (\ell_\bullet({h(s_{i_1})})\otimes \cdots \otimes \ell_\bullet({h(s_{i_k})}))_{\ol{\Delta}^{k-1} (T)}$. Thus, evaluating the linear combination $h^*(T) = \sum_{i_1,\ldots,i_k} \lambda_{(i_1,\ldots,i_k)} (S_{i_1}|\ldots|S_{i_k})$ on $w$ is just as claimed.
\end{proof}

Application of the composition formula \eqref{eq:composition formula} only requires understanding $\ell_\bullet({h(s)})$ for a particular finite set of generators. The following common setup is simpler still, with evaluation given by simple weight reduction of only leading terms.

\begin{cor}\label{cor:evaluation on high commutators}
    Consider a group $\Gamma$ and suppose $w(s_1,\ldots,s_d)\in \Gamma$ belongs to the $n$-th step of the lower central series, presented by some word $w\in \gamma_n F_S$. Then on weight $\leq n$ invariants, the letter-braiding number $\ell_\bullet({w(s_1,\ldots,s_d)})\maps \H^0_{Bar}(\Gamma;A) \to A$ is computed via the weight reduction algorithm of Proposition \ref{prop:weight-reduction} applied to the leading terms
    \[
    \H^0_{Bar}(\Gamma;A)_{\leq n} \xrightarrow{\lead} \H^1(\Gamma;A)^{\otimes n} \xrightarrow{(\operatorname{res}|_S){}^{\otimes n}} (A^S)^{\otimes n} \xrightarrow{\ell_\bullet(w)} A,
    \]
    where $\operatorname{res}|_S\maps \H^1(\Gamma;A)\to A^S$ is the restriction of homomorphisms $\H^1(\Gamma;A) \cong \Hom_{grp}(\Gamma,A)$ to $S\subseteq \Gamma$.
\end{cor}
\begin{proof}
    By naturality of the weight spectral sequence in Lemma \ref{lem:Bar spectral sequence} and the `leading term' map of Corollary $\ref{cor:leading terms}$, there is a commutative diagram
    \[
    \xymatrix{
    \H^0_{Bar}(\Gamma;A)_{\leq n} \ar[r]^-{\operatorname{lead}} \ar[d]_{h^*} &  \H^1(\Gamma;A)^{\otimes n} \ar[d]_{\operatorname{res|_S{}^{\otimes n}}}  & \\
    \H^0_{Bar}(F_S;A)_{\leq n} \ar[r]^-{\operatorname{lead}} & (A^S)^{\otimes n} \ar[r]^-{\ell(w)} & A.
    }
    \]
    Since $w$ is an $n$-fold commutator, its letter-braiding numbers vanish on weights $\leq n-1$ by Corollary \ref{cor:bounded degree general groups}, so factor through the leading term map. In other words, the composition of the bottom two arrows coincides with $\ell_\bullet(w)\maps \H^0_{Bar}(F_S;A)\to A$ on weights $\leq n$. The claim now follows form the commutativity of the diagram, along with 
    \[
    \ell_{T}(w(s_1,\ldots,s_d)) = \ell_{h^*(T)}(w).\qedhere
    \]
\end{proof}

\begin{example}[Pure braid group]\label{ex:braids}
    The cohomology ring of Artin's pure braid group $PB_n$ was calculated by Arnol'd \cite{Arnold} to be
    \[
    \H^*(PB_n;A) \cong \frac{\Lambda^*_A(\omega_{ij} \mid 1\leq i<j\leq n)}{( \omega_{ij}\wedge \omega_{jk}+ \omega_{jk}\wedge \omega_{ik} + \omega_{ik}\wedge \omega_{ij}  \mid 1\leq i<j<k \leq n)},
    \]
    where $\omega_{ij}$ measures the winding of the $i$-th strand around the $j$-th.
    \begin{cor}
        The tensors $T_{ijk}=(\omega_{ij}|\omega_{jk})+(\omega_{jk}|\omega_{ik})+(\omega_{ik}|\omega_{ij}) $  for $1\leq i<j<k\leq n$ define additive invariants on commutator braids $\ell_\bullet(T_{ijk})\maps [PB_n,PB_n]\to A$. They are given on basic commutator braids $[\sigma_1,\sigma_2]$ by
        \[
        \ell_{(\omega_{i_1j_1}|\omega_{i_2j_2})}([\sigma_1,\sigma_2]) = \omega_{i_1j_1}(\sigma_1)\omega_{i_2j_2}(\sigma_2) - \omega_{i_1j_1}(\sigma_2)\omega_{i_2j_2}(\sigma_1).
        \]
    \end{cor}
    \begin{proof}
        This follows immediately from the last corollary and Proposition \ref{prop:invariants of commutators}.
    \end{proof}
    Arnol'd shows, moreover, that when $A=\Q$ all Massey products vanish. It follows that cup products are the only obstruction preventing tensors $\H^1(PB_n;\Q)^{\otimes n}$ from defining letter-braiding invariants.
    \begin{cor}\label{cor:weight3 braids}
        Denoting $\omega_{ij}$ by $ij$, then for $i<j<k<l$ the sum of tensors in each of the following pairs of parentheses
        \[
    \begin{pmatrix}
    ij|ik|ij & ij|ij|jk & ij|jk|ik\\
    jk|ij|ij &  & jk|ik|ik \\
    ik|jk|ij & ik|ik|jk & ik|ij|ik 
    \end{pmatrix}
    \text{and }
    \begin{pmatrix}
 & ij|jl|jk &  &  & ij|kl|jl & ij|jk|kl \\
 &  & jk|il|ik & jk|kl|il &  & jk|ik|kl \\
ik|kl|ij &  &  &  &  & ik|ij|kl \\
il|ik|ij & il|ij|jk & il|jk|ik &  &  &  \\
 & jl|il|jk &  & jl|jk|il &  &  \\
kl|il|ij &  &  & kl|jl|il & kl|ij|jl & 
 \end{pmatrix}
    \]
    give well-defined additive invariants on double-commutators\footnote{These invariants are technically defined over $\Q$, but they obviously take integer values. Had we considered their extension to the entire braid group, we would possibly encounter fractional invariants due to lower-weight correction terms with rational coefficients, but these corrections vanish on double commutators.} $[[PB_n,PB_n],PB_n] \to \Z$, where
    \[
    \begin{aligned}
        \ell_{(\omega_1|\omega_2|\omega_3)}([[\sigma_1,\sigma_2],\sigma_3]) 
 = \; &\left[ \omega_1(\sigma_1)\omega_2(\sigma_2) - \omega_1(\sigma_2)\omega_2(\sigma_1)\right] \omega_3(\sigma_3) \\  - &\omega_1(\sigma_3)\left[ \omega_2(\sigma_1)\omega_3(\sigma_2) - \omega_2(\sigma_2)\omega_3(\sigma_1) \right].
    \end{aligned}
    \]
    \end{cor}
    \begin{proof}
        The two grids of tensors have the property that within a column all tensors have the same right-most factor, while the other two columns exhibit either the Arnol'd relation or anti-symmetry of the cup product. It follows that the contraction of the first two columns using the cup product gives $0$.
        
        At the same time, within every row the left-most factor is the same while the other two factors exhibit Arnol'd or anti-symmetry. Thus the contraction of the last two columns also gives $0$. Together each of the two sums defines a bar cocycle.

        Calculation of the braiding number of a double-commutator follows from the same formula for the universal word $[[x,y],z]\in F_3$ and \eqref{eq:invariants of commutators}.
    \end{proof}
    As an application, let $\sigma_{ij}$ be the Artin generators of $PB_n$, which satisfy $\omega_{ij}(\sigma_{kl}) = \delta_{ik}\delta_{jl}$.
    \begin{cor}
        The double commutator braids $[[\sigma_{12},\sigma_{23}],\sigma_{12}]$ and $[[\sigma_{12},\sigma_{23}],\sigma_{34}]$ are not triple-commutators in any pure braid group $PB_n$ in which they are defined.
    \end{cor}
    \begin{proof}
        Indeed, the former braid pairs nontrivially with the $3\times 3$ invariant of Corollary \ref{cor:weight3 braids} with $(i,j,k)=(1,2,3)$, while the latter braid pairs nontrivially with the $6\times 6$ invariant with $(i,j,k,l)=(1,2,3,4)$. Any triple commutator would give invariant $0$.
    \end{proof}
    We have computed a weight $4$ invariant that contains $1\cdot (12|23|34|23)$, along with $80$ other tensor summands, which detects that $[[\sigma_{12},\sigma_{23}],[\sigma_{23},\sigma_{34}]]\in PB_n$ is not a $4$-fold commutator.
    
    \begin{problem}
        Find a basis of leading terms of elements in $\H^0_{Bar}(PB_n;\Q)$, perhaps as cord diagrams, and determine its coproduct expansion.
    \end{problem}
    
    Note that over $\Z$ the pure braid groups are not formal (see \cite{salvatore-braid}), and thus a priori $\H^0_{Bar}(PB_n;\Q)$ appears essentially larger than $\H^0_{Bar}(PB_n;\Z)$. However, the above example suggests that the two coalgebras of invariants capture the same information. Indeed,
\end{example}
\begin{prop}
        Let $K$ be the fraction field of $A$. When $\Gamma$ is a finitely generated group, the inclusion $\H_{Bar}^0(\Gamma;A)\otimes K \into \H^0_{Bar}(\Gamma;K)$ is an isomorphism. In particular, Massey products on $\H^1(\Gamma;A)$ that are $A$-torsion do not essentially restrict the coalgebra of braiding invariants: every invariant defined over $K$ is realized by one over $A$ up to some scalar multiple.
\end{prop}
\begin{proof}
    Since $\Gamma$ is finitely generated, there exists a simplicial set $Y$ with a finite $1$-skeleton and $\pi_1(Y,*)\cong \Gamma$. It follows that
    \[
    C^*(Y;A)\otimes K \to C^*(Y;K)
    \]
    is an isomorphism in degrees $\leq 1$, and it is always an injection in degree $2$, i.e. this map is $1$-connected. The same argument as in the proof of Lemma \ref{lem:independence of model} now shows that $Bar^*(Y;A)\otimes K \to Bar^*(Y;K)$ is $0$-connected.

    Since $K$ is flat over $A$, there is an isomorphism $\H^0_{Bar}(Y;A)\otimes K \cong \H^0(Bar(Y;A)\otimes K)$. By Lemma \ref{lem:independence of model}, the former is isomorphic to $\H^0_{Bar}(\Gamma;A)\otimes K$ , while the latter is $\H^0_{Bar}(\Gamma;K)$ using the previous paragraph.
\end{proof}
This comparison isomorphism no longer holds for infinitely generated groups.
\begin{example}
    Consider the group $F_\infty$ modulo relations $x_1 = (x_2)^2 = (x_3)^3 = \ldots$. Then $\H^1(\Gamma;\Z) = 0$ but $\H^1(\Gamma;\Q) = \Q$. Therefore, letter-braiding invariants over $\Z$ detect nothing, while over $\Q$ they distinguish $x_1$ from $1$.
\end{example}

\section{Group rings}
\label{sec:group rings}
As discussed above, letter-braiding numbers define functionals $\ell_T(\bullet)\maps \Gamma \to A$, for every $T\in \H^0_{Bar}(\Gamma;A)$. These extend $A$-linearly to the group ring $A[\Gamma]$.
\begin{definition}
    The letter-braiding pairing
    \[
    \langle -, \bullet \rangle \maps \H^0_{Bar}(\Gamma;A)\otimes A[\Gamma] \to A
    \]
    is the $A$-linear extension of $\langle T,w \rangle = \ell_T(w)+\eta(T)$ from $\Gamma$ to $A[\Gamma]$, where $\eta \maps\H^0_{Bar}(\Gamma;A)\to A$ is the counit. It further extends to a pairing between their tensor powers as in \S\ref{sec:group ring prelims}.
\end{definition}
The various properties proved in \S\ref{sec:general groups} are most naturally stated in terms of this pairing, matching the functorial filtered coalgebra structure of $\H^0_{Bar}(\Gamma;A)$ with the filtered algebra structure of $A[\Gamma]$, as stated in Theorem \ref{thm:intro-pairing}.

\begin{remark}[\textbf{Diagonal}]
    Recall that $A[\Gamma]$ is furthermore equipped with a diagonal coproduct $\Delta(w) = w\otimes w$ for $w\in \Gamma$, making it into a Hopf algebra. In contrast, the bar construction of a commutative dg algebra has a shuffle product, making it into a Hopf algebra as well. Chen \cite{chen1971pi1} matches these structures on $\R[\Gamma]$ and $\H^0_{Bar}(\Omega_{dR}(M))$ under his pairing by iterated integrals. However, since our cochain algebras are rarely commutative, the shuffle product does not commute with $d_{Bar}$ and does not induce a Hopf algebra structure on $\H^0_{Bar}(\Gamma;A)$.
    
    Instead, there is an $E_{\infty}$-algebra structure on cochains, as instrumental in Mandell's theory \cite{mandell}. For our low dimensional applications, only the first layer of this homotopy commutativity structure is needed: the cup-1 operation, as central to Porter--Suciu work on 1-minimal models \cite{porter-suciu:cup1}.
    \begin{problem}
        Describe the product structure on $\H^0_{Bar}(\Gamma;A)$ dual to the coproduct on $A[\Gamma]$.
    \end{problem}
\end{remark}

For the proof of Theorem \ref{thm:intro-pairing}, we first extend the product formula 
to group rings.
\begin{lemma}\label{lem:product formula group ring}
    Let $\epsilon\maps A[\Gamma]\to A$ denote the canonical augmentation $\epsilon(\sum \lambda_i w_i) = \sum \lambda_i$. Then for any $x,y\in A[\Gamma]$ and $T\in \ol{\H^0_{Bar}(\Gamma;A)}$ in the augmentation co-ideal $\eta(T)=0$,
    \begin{equation*} %\label{eq:product formula group rings}
        \langle T,xy\rangle = \epsilon(y)\cdot \langle T,x\rangle + \epsilon(x)\cdot \langle T,y\rangle + \langle \ol{\Delta}(T),x\otimes y\rangle.
    \end{equation*}
\end{lemma}
\begin{proof}
    Let $x = \sum \lambda_i w_i$ and $y = \sum \mu_j w_j$. Then
    \[
    \langle T, xy\rangle  = \sum \lambda_i \mu_j \cdot \ell_{T}(w_iw_j) = \sum \lambda_i \mu_j \cdot \left( \ell_{T}(w_i) + \ell_{T}(w_j) + (\ell_{\bullet}(w_i)\otimes \ell_\bullet({w_j}))_{\ol\Delta(T)} \right).
    \]
    This is exactly $(\sum \mu_i) \cdot \langle T, x\rangle + (\sum \lambda_i)\cdot \langle T, y\rangle + \langle \ol{\Delta}T, x\otimes y\rangle$, as claimed.
\end{proof}
%Next, introduce some notation.
\begin{definition}
    For words $w_0,\ldots,w_n\in \Gamma$, let the \emph{multi-evaluation} 
    \begin{equation*}
        \ell_{\bullet}(w_0|\ldots|w_n)\maps \H^0_{Bar}(\Gamma;A) \to A
    \end{equation*}
    be the composite $
    (\ell_{\bullet}(w_0)\otimes \ldots \otimes \ell_\bullet({w_n})) \circ \ol\Delta^{n}.$
\end{definition}
\begin{prop} \label{prop:properties of multi-eval}
    Let $w_0,\ldots,w_n\in \Gamma$ be any sequence of elements. The multi-evaluation $\ell_\bullet({w_0|\ldots|w_n})$ satisfies the following.
    \begin{enumerate}

        \item \label{case:multi in terms of single} It is a linear combination of letter-braiding invariants of products,
        \begin{equation*}%\label{eq:multi-eval as grp elements}
        \ell_\bullet({w_0|\ldots|w_n}) = \sum_{k=1}^n  (-1)^{n-k} \sum_{0\leq i_1< \ldots < i_k \leq n}  \ell_\bullet({w_{i_1}\cdots w_{i_k}}).
        \end{equation*}
        \item In terms of the letter-braiding pairing,
        \begin{equation*}%\label{eq:multi-eval as group ring}
            \ell_T({w_0|\ldots|w_n}) = \langle T \;,\; (w_0-1)\cdots (w_n-1)\rangle.
        \end{equation*}
        \item For every $0\leq i<n$,
        \begin{equation*}%\label{eq:multi-eval of coproduct}
        (\ell_\bullet({w_0|\ldots| w_i}) \otimes \ell_\bullet({w_{i+1}|\ldots| w_n})) \circ \ol\Delta = \ell_\bullet({w_0|\ldots|w_n}).
        \end{equation*}
    \end{enumerate}
\end{prop}
\begin{proof}
    The third statement follows immediately from coassociativity of $\Delta$. Next, the first and the second equations are evidently equivalent, so let us prove the latter by induction on $n$. Consider $x= (w_0-1)\cdots (w_n-1)\in A[\Gamma]$ and assume that
    $\ell_\bullet({w_0|\ldots|w_n}) = \langle \bullet, x\rangle$. Using the product formula for group rings,
    \[
    \langle \bullet, {x(w_{n+1}-1)}\rangle = \cancel{\epsilon(w_{n+1}-1)} \cdot \langle \bullet, x\rangle + \cancel{\epsilon(x)} \cdot \langle \bullet, w_{n+1}-1\rangle + \langle \ol\Delta(\bullet), x\otimes (w_{n+1}-1)\rangle,
    \]
    \[
    = \langle \ol\Delta(\bullet), x\otimes w_{n+1}\rangle.
    \]
    The last expression is $(\ell_\bullet({w_0|\ldots|w_n})\otimes \ell_\bullet{(w_{n+1})})\circ \ol\Delta$ and thus agrees with $\ell_\bullet({w_0|\ldots|w_{n+1}})$.
\end{proof}

\begin{remark}\label{rmk:invariant of commutator}
    With Proposition \ref{prop:properties of multi-eval} we can prove the recursive formula for invariants of commutators.
    \begin{proof}[Proof of Proposition \ref{prop:invariants of commutators}]
        In the group ring $A[F_S]$ there is an identity
        \begin{align*}
        [w,s]-1 &= wsw^{-1}s^{-1}-1 = wsw^{-1}s^{-1}(1-s) + wsw^{-1}(1-w) + w(s-1) +(w-1) \\ &=    
        (w-1)(s-1) - ([w,s]-1)(s-1) -(wsw^{-1}-1)(w-1).
        \end{align*}
        Therefore, by Proposition \ref{prop:properties of multi-eval},
        \[
        \langle T, [w,s]\rangle = \ell_{T}(w|s) - \ell_{T}([w,s]|s) -\ell_{T}(wsw^{-1}|w).
        \]
        We note the following: $\ell_{T'}(s)=0$ on any homogeneous tensor of weight $\geq 2$, and $\ell_{T'}([w,s])=0$ on any tensor of weight $\leq n$. Similarly $\ell_{T'}(w)=0$ on tensors of weight $\leq n-1$. When $T = \alpha_0|\ldots|\alpha_n$, these observations leave only the terms
        \[
        \langle (\alpha_0|\ldots|\alpha_n)\,,\, [w,s]\rangle = \ell_{(\alpha_0|\ldots|\alpha_{n-1})}(w)\ell_{\alpha_n}(s) -\ell_{\alpha_0}(wsw^{-1})\ell_{(\alpha_1|\ldots|\alpha_n)}(w).
        \]
        But since $\ell_{\alpha_0}(\bullet)$ is an additive homomorphism, we have $\ell_{\alpha_0}(wsw^{-1}) = \ell_{\alpha_0}(s)$.
    \end{proof}
\end{remark}

\begin{proof}[Proof of Theorem \ref{thm:intro-pairing}]
    The push-pull formula \eqref{eq:push-pull} shows that the pullback $h^*$ along a group homomorphism is dual to the pushforward $h_*$.
    Now, observing the the braiding pairing of the unit $1\in \H^0_{Bar}(\Gamma;A)$ is $\langle 1, x\rangle = \epsilon(x)$ for every $x\in A[\Gamma]$, the product formula for group rings (Lemma \ref{lem:product formula group ring}) becomes \begin{align*}
        \langle \Delta(T), x\otimes y\rangle &= \langle (T\otimes 1) + (1\otimes T) + \ol{\Delta}(T) \,,\, x\otimes y\rangle\\  &= \langle T,x\rangle \cdot \epsilon(y) + \epsilon(x) \cdot \langle T,y\rangle + \langle \ol{\Delta}(T), x\otimes y \rangle = \langle T ,xy\rangle.
    \end{align*}
Let $I\lhd A[\Gamma]$ denote the augmentation ideal. We prove that weight $\leq n-1$ invariants vanish on $I^{n}$ by induction on $n$. The fact that the only weight $0$ invariant up to scaling is $\langle 1,x\rangle = \epsilon(x)$ and $I = \ker(\epsilon)$ gives the base case. For the step, suppose $x\in I^k$ and $y\in I^{n-k}$ for $0<k<n$. Then
\[
\langle T, xy\rangle = \langle \Delta(T), x\otimes y\rangle = \sum_i \langle T_1^i,x\rangle\cdot \langle T_2^i,y\rangle,
\]
where $\Delta(T) = \sum_i T_1^i\otimes T_2^i$. But for every $i$ the weights of $T_1^i$ and $T_2^i$ must sum to at most $n-1$, so at either the first has weight $< k$ or the second has weight $<n-k$. With the induction hypothesis, it follows that every summand in the above sum vanishes.

Now, let $I\lhd A[\Gamma]$ be the augmentation ideal and let $T\in \H^0_{Bar}(\Gamma;A)$ have weight $< n$. we claim that $\ell_T(\bullet)$ vanishes identically on $I^n$. By induction on $n$, if $T$ has weight $0$ then it vanishes on $I$ by definition, since $\langle 1, x\rangle = \epsilon(x)$. Now let $T$ have weight $\leq n-1$ and suppose by induction that $(I^{n-1})^{\perp} = \H^0_{Bar}(\Gamma;A)_{<n-1}$, then if $x\in I^{n-1}$ and $y\in I$, we have
\[
\langle T, xy\rangle = \sum_i\langle T_1^i, x\rangle \langle T_2^i, y\rangle
\]
where $\Delta(T) = \sum_i T_1^i\otimes T_2^i$. But since for every $i$, the weights of $T_1^i$ and $T_2^i$ sum to $\leq n-1$, one of the respective evaluations on $x$ and $y$ must vanish. Since $I^n$ is generated by elements of the form $xy$, it follows that $\ell_T(\bullet)\equiv 0$ on $I^n$.

Conversely, if $0\neq T\in \H^0_{Bar}(\Gamma;A)$ has weight $n$, then there exists some $x\in I^n$ such that $\langle T,x\rangle \neq 0$. Indeed, let $S\subseteq \Gamma$ be a generating set, so that the restriction $\H^1(\Gamma;A)\to A^S$ is injective, and consider the leading term
\[
\lead(T) \in \H^1(\Gamma;A)^{\otimes n}.
\]
On tensor powers, restriction gives an injection $\H^1(\Gamma;A)^{\otimes n} \into (A^S)^{\otimes n}$, and by Lemma \ref{lem:tensor of infinite powers} there is a further inclusion $(A^S)^{\otimes n} \into A^{S\times \cdots \times S}$. Therefore, if $T\neq 0$ then there exist generators $(s_1,\ldots,s_n)$ such that
\[
\langle T \,,\, (s_1-1)\cdots(s_n-1)\rangle = \langle \ol\Delta^n(T) \,,\, (s_1\otimes\ldots\otimes s_n)\rangle = \langle \lead(T) \,,\, (s_1\otimes\ldots\otimes s_n)\rangle \neq 0.
\]
This fact already shows that $\H^0_{Bar}(\Gamma;A)\to \Hom_A(A[\Gamma],A)$ is injective. Completeness of letter-braiding is more difficult and will be the subject of \S\ref{sec:completeness}.
\end{proof}
\begin{example}\label{ex:abelian group ring}
    We saw in Example \ref{ex:abelian group} that in $F_{\{e,f\}}$ subject to the relation $ef = fe$, all words are already distinguished by the invariants $E$ and $F$. Yet $\H^0_{Bar}(\Z\langle e,f\rangle;\Z)$ contains further invariants, such as $(E|E)$ and $(E|F) + (F|E)$. Their role is now made clear: the two elements $ef$ and $e+f$ in $\Z[\Gamma]$ can not be separated using neither $E$ nor $F$, while
    \[
    \langle (E|F)+(F|E) \,,\, ef\rangle = 1 \quad \langle (E|F)+(F|E)\,,\, e+f\rangle = 0.
    \]
    To wit, higher weight invariants in $\H^0_{Bar}(\Gamma;A)$ distinguish between group ring elements.
\end{example}

\subsection{Fox calculus and Magnus expansion}\label{sec:magnus}
For coefficients in $A = \Z$, the letter-braiding number is equivalent to Fox calculus \cite{Fox1} and dual to the Magnus expansion \cite{Magnus_1935}, in the following sense.
\begin{theorem}\label{thm:magnus}
    Let $F_n$ be the free group generated by $\{x_1,\ldots,x_n\}$, and let $X_1,\ldots X_n\in \Z^{n}$ be the dual basis. Fix $w\in F_n$. Then,
    \begin{itemize}
        \item for every sequence of indices $(i_1,\ldots i_k) \in \{1,\ldots,n\}^k$ there is an equality
        \[
        \ell_{(X_{i_1}|\ldots|X_{i_k})}(w) = \epsilon\left( \frac{\partial}{\partial x_{i_1}}\ldots\frac{\partial}{ \partial x_{i_d}}w \right)
        \]
        where $\frac{\partial}{\partial x_i}$ is the Fox derivative on $\Z[F_n]$ associated to $x_i$, and $\epsilon\maps \Z[F_n]\to \Z$ is the standard augmentation.
        \item There is an equality of noncommutative power series
        \[
        M(w) = 1+ \sum_{1 \leq i_1,\ldots,i_k \leq n} \ell_{(X_{i_1}|\ldots|X_{i_k})}(w) t_{i_1}\cdots t_{i_k}
        \]
        where $M(-)$ is the Magnus expansion of \cite{Magnus_1935}.
    \end{itemize}
\end{theorem}
\begin{proof}
    On generators $x_j$, we computed in Example \ref{ex:braiding of a generator} that 
    \[
    \ell_{(X_{i_1}|\ldots|X_{i_k})}(x_j) = \begin{cases}
        1 & \text{ if } (X_{i_1}|\ldots|X_{i_k}) = X_j \\
        0 & \other.
    \end{cases}
    \]
    It follows that $\ell_\bullet(w)$ agrees with the augmentation of Fox derivatives on generators of $F_n$.

    Now, for products $w = u\cdot v$ \cite[Equation (3.2)]{Fox1} shows that Fox derivatives satisfy the recurrence relation as given in Corollary \ref{cor:invariant of products} for $\ell_{\bullet}(uv)$. Since the two agree on generators, they must be equal.

    Characterization of the coefficients of Magnus' expansion in terms of the Fox derivatives is well-known, see e.g. \cite[Discussion following (3.4)]{Fox1}.
\end{proof}
\begin{remark}\label{rmk:magnus}
In light of the last theorem, our setup can be understood as an extension of the Magnus expansion on free groups in the following ways.
\begin{itemize}
    \item The group $\Gamma$ needs not be free, where now the only well-defined Fox derivatives are inhomogeneous linear combinations that exhibit vanishing cup and Massey products.
    \item The group $\Gamma$ need not be finitely generated.
    \item The coefficient ring $A$ need not be $\Z$. E.g. it may be a field of characteristic $p>0$.
    \item The calculus is independent of a presentation for $\Gamma$. Instead, the only necessary input is a presentation of the cohomology $\H^1(\Gamma;A)\cong \Hom(\Gamma,A)$ with its Massey products
    \[\mu_n\maps \H^1(\Gamma;A)^{\otimes n} \to \H^2(\Gamma;A).\]
\end{itemize}
\end{remark}

Let us illustrate a major improvement afforded by allowing coefficients in general PIDs. While rational homotopy theory gives insight into the ranks of the lower central series quotients of a group, it can never detect words in finite groups. More generally, any invariant theory with values in a torsion-free ring $A$ will necessarily have $\H^1(\Gamma;A) = 0$ for any finite group $\Gamma$.
\begin{example}\label{ex:heisenberg}
Take $A = \F_2$ and consider the Heisenberg group $\Gamma = \mathscr{H}_3(\F_2)$ of matrices
\[
\left\{ \begin{pmatrix}
    1 & X & Z \\
    & 1 & Y \\
    & & 1
\end{pmatrix} \text{ s.t. } X,Y,Z\in \F_2 \right\}  \subseteq \operatorname{Gl}_3(\F_2).
\]
It has $8$ elements, a presentation $\langle x,y,z \mid x^2,y^2,z^2,[x,y]z^{-1} \rangle$, and $\H^1(\Gamma;\F_2) = \langle X, Y\rangle$ given by the sub-diagonal entries.

On the group cochain level, the functions $X,Y$ and $Z$ are defined in $C^1(B\Gamma;\F_2)$, but $Z$ is not closed. Instead, $dZ = X\smile Y$, so that the cup product vanishes on cohomology:
\[
dZ(w,w') = Z(ww')-Z(w)-Z(w') = X(w)Y(w') = (X\smile Y)(w,w').
\]
A well-defined invariant is therefore $(X|Y)+Z$, which is indeed unchanged by application of any of the relations.

Now, $X$ and $Y$ detect the image of any group element in the abelianization. But, more interestingly, the invariant $X|Y+Z$ detects the commutator $[x,y]=z$,% = \begin{psmallmatrix}    1&0&1\\0&1&0\\0&0&1\end{psmallmatrix}$.
\end{example}
More generally, finite $p$-groups are well-tailored for $\F_p$-valued letter-braiding invariants.

\section{Completeness of letter-braiding}
\label{sec:completeness}
We next show that letter-braiding invariants furnish a \emph{complete} coalgebra of functionals that are continuous with respect to the augmentation ideal topology on the group ring. Recall that the group ring $A[\Gamma]$ is filtered by powers of the augmentation ideal $I\lhd A[\Gamma]$. The $I$-adic completion is inverse limit
    \[
    A[\Gamma]^{\wedge}_I = \varprojlim A[\Gamma]/I^n
    \]
    topologized so that all quotients $A[\Gamma]^{\wedge}_I\onto A[\Gamma]/I^n$ are continuous, when the target has the discrete topology. In particular, a functional $\varphi\maps A[\Gamma]^{\wedge}_I\to A$ is \emph{continuous} if and only if it vanishes on some power $I^N$.

\begin{theorem}[\textbf{Completeness of letter-braiding invariants}] \label{thm:completeness}
\begin{comment}
    Let $\Gamma$ be any finitely generated group. For every $A$-linear functional $\varphi\maps A[\Gamma]\to A$ that vanishes on some power of the augmentation ideas $I^n$, there exists a unique $T\in \H^0_{Bar}(\Gamma;A)$ of weight $\leq n$ such that $\langle T, \bullet\rangle = \varphi(\bullet)$.
\end{comment}
    Let $\Gamma$ be a finitely generated group. The letter-braiding pairing defines a natural isomorphism of filtered coalgebras
    \[
    \H^0_{Bar}(\Gamma;A) \xrightarrow{\sim} \Hom_{cts}(A[\Gamma]^{\wedge}_I , A) = \varinjlim_n \Hom_A( A[\Gamma]/I^n , A)
    \]
    where the $A[\Gamma]^{\wedge}_I$ is the $I$-adic completion of the group ring equipped with the $I$-adic topology, and $\Hom_{cts}(-,-)$ denotes continuous linear functionals.
\end{theorem}
In other words, for every functional $\varphi\maps A[\Gamma]\to A$ that vanishes on some power $I^n$, there exists a unique $T\in \H^0_{Bar}(\Gamma;A)$ of weight $\leq n$ such that
\[
\langle T,\bullet\rangle = \varphi(\bullet).
\]

\begin{remark}\label{rmk:infinite generation}
    When the group $\Gamma$ is infinitely generated there is no hope that letter-braiding would be an isomorphism in general. For example, $\H^0_{Bar}(\Gamma;A)$ is always a coalgebra, while duals of infinite dimensional algebras are not -- this is due to the failing of $(V\otimes V)^* \cong V^*\otimes V^*$. Instead, in this case the coalgebra $\H^0_{Bar}(\Gamma;A)$ is universal conilpotent dual of $\Gamma$, as stated in Theorem \ref{intro-thm:infinite type universal pairing} and proved at the end of this section.
\end{remark}

\begin{comment}
    
In terms of functionals on the group itself, the letter-braiding invariants are not complete unless division is allowed. In a forthcoming sequel, joint with A. Ozbek, D. Sinha, B. Walter, we discuss the relations between the Malcev Lie algebra of a group and the Harrison complex of its commutative cochain model -- the commutative version of the associative bar complexes considered here.

Over the remainer of this section we develop the tools needed for the proof of completeness. After completing the proof, we discuss natural restatements in a more categorical language.
\end{comment}

\subsection{Introducing relations}\label{sec:introducing relations}

Building towards a proof of Theorem \ref{thm:completeness}, we quantify the effect of introducing relations to a group. The main result of this subsection is that the letter-braiding invariants of a quotient group $\Gamma_0/\langle R\rangle$ are exactly those invariants on $\Gamma_0$ that vanish in a strong sense on all $r\in R$.

Let $\Gamma_0$ be a group and $R\subseteq \Gamma_0$ a subset that will serve as relations for $
\Gamma := \Gamma_0/\langle R\rangle$. 
The quotient map $p\maps \Gamma_0\to \Gamma$ induces a pullback of filtered coalgebras
\begin{equation*}
    p^*\maps \H^0_{Bar}(\Gamma;A) \longto \H^0_{Bar}(\Gamma_0;A).
\end{equation*}
\begin{theorem}[\textbf{Lifting criterion}] \label{thm:lifting criterion}$p^*$ is an injection, identifying $\H^0_{Bar}(\Gamma;A)$ as a sub\-coalgebra of $\H^0_{Bar}(\Gamma_0;A)$. Moreover, for $T\in \ol{\H^0_{Bar}(\Gamma_0;A)}$, the following are equivalent:
    \begin{enumerate}
        \item \label{case:lifting} $T$ is in the sub-coalgebra $p^*\H^0_{Bar}(\Gamma;A)$.
        \item \label{case:splitting lifts}One can write $\ol\Delta (T) = \sum_{i=1}^m T_1^i\otimes T_2^i$ with all $T_1^i$ and $T_2^i$ in the sub-coalgelbra $p^*\H^0_{Bar}(\Gamma;A)$, and moreover $\ell_T(r) = 0$ for all $r\in R$.
        \item \label{case:group ring factors} The group ring functional $\langle T,\bullet\rangle\maps A[\Gamma_0]\to A$ factors through the projection 
        \[p_*\maps A[\Gamma_0]\onto A[\Gamma].\]
        \item \label{case:vanishing multi-evaluations} For every sequence of group elements $w_0,\ldots,w_k \in \Gamma_0$ such that at least one of them belongs to $R$, the multi-evaluation
        \begin{equation*}
            \ell_{T}(w_0|\ldots|w_n) = 0.
        \end{equation*}
    \end{enumerate}
\end{theorem}
Injectivity of $p^*$ follows from the following more general fact.
\begin{lemma}\label{lem:injective on bar}
    Whenever $p\maps \Gamma_1\to \Gamma_2$ induces an injection $p^*\maps \H^1(\Gamma_2;A)\into \H^1(\Gamma_1;A)$, in particular when $p$ induces a surjection on abelianization $\Gamma_1^{ab}\to \Gamma_2^{ab}$, the pullback
    \[
    p^*\maps \H^0_{Bar}(\Gamma_2;A) \to \H^0_{Bar}(\Gamma_1;A)
    \]
    is injective.
\end{lemma}
\begin{proof}
    Consider the induced map of the weight spectral sequences, as described in Lemma \ref{lem:Bar spectral sequence}. On their $E_1$-pages in total degree $0$ it consists of injections
    \[
    p^*\maps \H^1(\Gamma_2;A)^{\otimes n} \into \H^1(\Gamma_1;A)^{\otimes n}.
    \]
    But since the $E_\infty$-pages include inside the respective $E_1$-pages, $p^*$ induces an injection on $E_{\infty}$. This already implies injectivity on $\H^0_{Bar}(-;A)$.
\end{proof}

\begin{proof}[Proof of \eqref{case:lifting} $\implies$ \eqref{case:splitting lifts}]
 If $T = p^*\tilde{T}$ for some $\tilde{T}\in \H^0_{Bar}(\Gamma;A)$ and $\ol\Delta(\tilde{T}) = \sum \tilde{T}_1^i \otimes \tilde{T}_2^i$, then
    \[
    \ol\Delta(T) = \ol\Delta(p^*\tilde{T}) = \sum p^* \tilde{T}_1^i \otimes p^* \tilde{T}_2^i,
    \]
    so the components of the comultiplication indeed lift. By push-pull, $\ell_{p^*\tilde{T}}(r) = \ell_{\tilde{T}}(p(r)) = 0$ since $r$ is in the kernel of $p$.
\end{proof}

\begin{proof}[Proof of \eqref{case:splitting lifts} $\implies$ \eqref{case:group ring factors}:]
    We must show that $\ell_\bullet(T)$ vanishes on the ideal 
    \[I_R := \ker\left( A[\Gamma_0]\onto A[\Gamma] \right).\]
    Fox observed in \cite[\S1]{Fox1} that the ideal $I_R$ is generated by the subset $\{ r-1 \}_{r\in R}$.
%    ideals $I'\lhd A[\Gamma_0]$ contained in the augmentation ideal are in natural bijection with normal subgroups $K\lhd \Gamma$ via the rule
%    \[
%    k\in K \iff (k-1)\in I'.
%    \]
%    Under this bijection, the group $\langle R\rangle$ normally generated by $r\in R$ corresponds to the ideal $I_R$, which is therefore generated by the subset $\{ r-1 \}_{r\in R}$.
    On every such generator, Condition \eqref{case:splitting lifts} gives
    \[
    \langle T, r-1\rangle = \ell_T(r) - \ell_T(1) = 0 - 0 = 0.
    \]
    To see that vanishing extends to any product with generators, we employ the product formula for $A[\Gamma_0]$. Suppose that $x\in I_R$ was already proven to be in the kernel of $\langle T,\bullet\rangle$ and $y\in A[\Gamma_0]$ is any element. Then by the product formula
    \[
    \langle T, xy \rangle = \epsilon(y)\cdot \cancel{\langle T, x \rangle} + \cancel{\epsilon(x)}\cdot \langle T, y \rangle + \langle \ol{\Delta}T, x\otimes y \rangle,
    \]
    where $\langle T, x \rangle=0$ by assumption and $\epsilon(x) = 0$ since $x$ belongs to the augmentation idea. But Condition \ref{case:splitting lifts} implies that $\ol\Delta(T) = \sum T_1^i\otimes T_2^i$ where all terms $T_1^i$ and $T_2^i$ are pulled back from $\Gamma = \Gamma_0/\langle R\rangle$, so $\langle T_1^i, x\rangle =0$ for all $i$. Together, this shows that $\langle T, xy \rangle = 0$.

    The same proof shows that $yx$ also pairs trivially with $T$, where now we use vanishing of $\langle T_2^i,x\rangle$. With these, $\langle T,\bullet\rangle$ vanishes on the two-sided ideal generated by $\{r-1\}_{r\in R}$, thus factors through $A[\Gamma]$.
\end{proof}

\begin{proof}[Proof of \eqref{case:group ring factors} $\implies$ \eqref{case:vanishing multi-evaluations}]
    Assume that $\langle T, \bullet \rangle$ factors through $p_*\maps A[\Gamma_0]\onto A[\Gamma]$. We must show that $\ell_{T}(w_0|\ldots|w_n) = 0$ whenever one or more of the $w_i$'s belongs to $R$. Let $w_0,\ldots,w_n$ be a sequence of group elements with $w_i = r \in R$ for some $0\leq i \leq n$. Then by Proposition \ref{prop:properties of multi-eval},
    \[
    \ell_{T}(w_0|\ldots|w_n) = \langle T \,,\, (w_0-1)\cdots (w_n-1)\rangle.
    \]
    But the product $(w_0-1)\cdots (r-1)\cdots (w_n-1)$ is sent to $0 \in A[\Gamma]$ and so this multi-evaluation vanishes.
\end{proof}
The last step of the proof is the most elaborate, first requiring substantial  terminology. 
%This is the content of the next subsection.

\subsection{Defect and obstruction to lifting}\label{sec:defect}
To complete the proof of the lifting criterion of Theorem \ref{thm:lifting criterion} we  show that if $T$ does not lift to $\H^0_{Bar}(\Gamma;A)$, then there exists a multi-evaluation containing $r\in R$ that does not vanish on $T$.
To achieve this we introduce a secondary ``defect" filtration on $\H^0_{Bar}(\Gamma_0;A)$.

In the construction below, a `space' $Y$ can be taken to mean either a topological space or a simplicial set. When working in the topological category, all maps are continuous and $S^1(r)$ denotes a topological circle. Otherwise, maps are simplicial and $S^1(r)$ is a subdivided simplicial circle as defined in \S\ref{sec:weight reduction}. Let $D^2(r)$ denote the reduced cone on $S^1(r)$ -- this is a disc with boundary $S^1(r)$.
\begin{definition}[\textbf{Presentation complex}]
    Let $(Y,*)$ be a pointed space and $R\subseteq \pi_1(Y,*)$ a subset of elements, realized by pointed maps $\{ S^1(r)\to Y \}_{r\in R}$. Then $Y_{|R}$ is the reduced mapping cone
    \[
    Y_{|R} = \operatorname{Cone}\left( \vee r\maps \bigvee_{r\in R}S^1(r) \to Y \right) = \left(\bigvee_{r\in R} D^2(r)\right) \coprod_{\vee_{r\in R}S^1(r)} Y.
    \]
Denote $\Gamma_0 := \pi_1(Y,*)$ and $\Gamma := \pi_1(Y_{|R},*)$.
\end{definition}
By Van-Kampen's theorem, the natural inclusion $p\maps Y\into Y_{|R}$ induces an isomorphism
\[
p_*\maps \pi_1(Y_{|R},*) \xrightarrow{\sim} \pi_1(Y,*)/\langle R\rangle.
\]
\begin{comment}
\begin{example}
Let $X_S$ be a presentation complex for the free group $F_S$, and fix a collection of simplicial maps $\{S^1(\tilde{r})\to X_S\}_{r\in R}$, as introduced in \S\ref{sec:combinatorial}. Then $X_{S|R}$ is the presentation complex for the group $\Gamma = \langle S|R\rangle$.
\end{example}
\end{comment}

The long exact sequence of the pair $(Y_{|R},Y)$ gives an exact sequence
\begin{equation*}
    0 \to \H^1(\Gamma;A) \xrightarrow{p^*} \H^1(\Gamma_0;A) \xrightarrow{d} A^R
\end{equation*}
where $A^R$ is the relative cohomology  $\H^2(Y_{|R},Y;A) \cong \tilde\H^2(\bigvee_{r\in R} S^2;A)$.
\begin{lemma}
    Under the identification $\H^1(\Gamma_0;A)\cong \Hom_{grp}(\Gamma_0,A)$, the connecting homomorphism $d$ is the restriction of functions $\Gamma_0\to A$ to $R\subseteq \Gamma_0$.
\end{lemma}
This exact sequence is the first instance of the lifting criterion: a functional $\alpha\maps \Gamma_0\to A$ comes from $H^1(\Gamma;A)$ if and only if it evaluated trivially on the relations.
We will establish a vast generalization of this for bar cohomology.

Observe that at the cochain level we get a surjection of augmented dg algebras
    \[
    C^*(Y_{|R};A) \xrightarrow{p^*} C^*(Y;A) \to 0.
    \]
    Its kernel is $K = C^*(Y_{|R},Y; A)$, having $\H^2(K) \cong A^R$ detecting obstructions to lifting, with vanishing cohomology in lower degrees.

For simplicity of notation, let us consider this setup in the abstract.
\begin{hypothesis}\label{hyp:defect setup}
Let $p\maps C \onto C_0$ be a surjection of connected dg algebras whose underlying $A$-modules are $A$-torsion free. Denote $K = \ker(p)$ and assume further that $\H^1(K) = 0$.
\end{hypothesis}

\begin{definition}
We next consider with complexes of bigraded modules, where the vertical grading is the usual cohomological grading, and the new horizontal grading is called \emph{defect degree}. A plain (singly-graded) complex will be considered as concentrated in defect degree $0$.

Let $\delta A$ denote the free module of bidegree $(-1,0)$. If $C^*$ is a cochain complex concentrated in defect degree $0$, the tensor $(\delta A)^{\otimes p} \otimes C^*$ is denoted $\delta^p C^*$ and is concentrated in defect degree $-p$.
\end{definition}

\begin{theorem}\label{thm:defect filtration}
    Let $p\maps C\onto C_0$ be a surjection of dg algebras satisfying Hypothesis \ref{hyp:defect setup}. Then there exists a functorial `defect degree' filtration $D_\bullet \H^0_{Bar}(C_0)$ with the following properties:
    \begin{itemize}
        \item The induced map $p\maps \H^0_{Bar}(C)\to \H^0_{Bar}(C_0)$ is an isomorphism onto defect degree zero
    \[
    H^0_{Bar}(C) \cong D_0 H^0_{Bar}(C_0).
    \]
    \item The associated graded coalgebra $\operatorname{gr}^D_\bullet \H^0_{Bar}(C_0)$ is filtered by a functorial secondary filtration $W_\bullet$, called \emph{defect weight}, such that there is a natural injective homomorphism of bigraded coalgebras -- the `defect map' --
\begin{equation*}%\label{eq:defect on coalgebras}
    \defect\maps \operatorname{gr}^{W} \operatorname{gr}^D \H^0_{Bar}(C_0) \into T( \delta\H^2(K) \oplus \H^1(C) )
\end{equation*}
    which sends defect weight $n$ and degree $p$ to tensors of weight $n$ containing exactly $p$ factors in $\delta \H^2(K)$.
        
    \item On defect degree $0$, the defect weight filtration coincides with the usual weight filtration on $\H^0_{Bar}(C)$, and the defect map on weight $n$ is $\ol{\Delta}^{n-1}\maps \H^0_{Bar}(C)_{\leq n} \to \H^1(\Gamma;A)^{\otimes n}$.
    \item On weight $1$ classes $\alpha \in \H^1(C_0) \leq \H^0_{Bar}(C_0)$ of defect degree $1$, the defect map is the coboundary of the long exact sequence
    \[
    \H^1(C) \to \H^1(C_0) \xrightarrow{d} \H^2(K).
    \]
        %\item On defect degree and weight $1$, the defect map is obtained by
        %\[
        %\defect([T]) = [ d_{Bar} \tilde{T} ] \in \H^2(K),
        %\]
        %where $\tilde{T}\in Bar^0(C)$ is any bar cochain whose image in $Bar^0(C_0)$ is cohomologous to $T$ and $d_{Bar}(\tilde{T})$ has weight $1$.
    \end{itemize}
\end{theorem}
%\begin{remark}
    Lifting the defect weight filtration to $\H^0_{Bar}(C_0)$, these are refinements of $D_\bullet$ by
    \[
    D_{p-1} = W_{p-1,p} \subseteq W_{p,p} \subseteq W_{p+1,p}\subseteq \ldots \subseteq W_{n,p} \subseteq \ldots  \subseteq D_p
    \]
    so that the defect map in the last theorem is a collection of linear maps
    \[
    \defect_{n,p}\maps W_{n,p} \longto T\left( \delta\H^2(K) \oplus \H^1(C) \right)
    \]
    for every $n\geq p$, whose image consists of tensors of weight $n$ and exactly $p$ many factors drawn from $\delta\H^2(K)$, and having
    \[
    \ker(\defect_{n,p}) = W_{n-1,p}.
    \]
    Comultiplication interacts with these subspaces by
    \begin{equation}\label{eq:coproduct on defect weight}
    \ol\Delta(W_{n,p}) \subseteq \sum_{\substack{m,q\\ 0 \leq m-q \leq n-p}} W_{m,q}\otimes W_{n-m,p-q} + \sum_{0\leq i \leq p-1} D_i\otimes D_{p-1-i}.
    \end{equation}
%\end{remark}
The next construction is key to defining the defect filtrations.
\begin{definition}
    Let $p\maps C \onto C_0$ be a surjection of connected dg algebras, with kernel $K$. Treat the inclusion $i\maps K \into C$ as a bicomplex 
    \[
     \xymatrix{
    \vdots & \vdots & \\
    \delta K^1 \ar[r]^{i_1}\ar[u]^{-d} & C^1 \ar[u]^d & =:(C,K)^{*,*}\\
    \delta K^0 \ar[r]^{i_0}\ar[u]^{-d} & C^0 \ar[u]^d &
    }
    \]
    where the horizontal grading is the \emph{defect degree}. 
    
    Extend $p$ to a bigraded projection $p\maps (C,K) \longto C_0$ by having it vanish of $\delta K$.
\end{definition}
\begin{lemma}
    There is a natural bigraded differential algebra structure on $(C,K)$, making the maps $C\into (C,K)\onto C_0$ into bigraded differential algebra homomorphisms, and the induced map $\tot (C,K)\to C_0$ a quasi-isomorphism.
\end{lemma}
\begin{proof}
    Since $K\lhd C$ is an ideal, the product on $C$ induces maps
    \[
    \delta K \otimes C \to \delta K \quad\text{and}\quad C \otimes \delta K \cong \delta C\otimes K \to \delta K,
    \]
    where we use the Koszul sign convention when moving $\delta$ past elements of $C$:
    \[
    c\cdot (\delta k) = (-1)^{|c|} \delta (c\cdot k)
    \]
    for homogeneous $c\in C$. Along with the existing product $C\otimes C\to C$ and the zero product $\delta K\otimes \delta K \to 0$, these give a bigraded algebra structure on $(C,K)$. Quick verification using the Koszul sign rule shows that both differentials are indeed derivations.

    With this product, $C\into (C,K) \onto C_0$ are homomorphisms of bigraded differential algebras. The short exact sequence of complexes $0\to K\to C\to C_0\to 0$ shows that $\tot(C,K)\to C_0$ is a quasi-isomorphism.
\end{proof}

Recalling that by Lemma \ref{lem:quasi-iso invariant} the bar construction is a quasi-isomorphism invariant of dg coalgebras, one can compute $\H^0_{Bar}(C_0)$ with the bar construction of $(C,K)$. But since $(C,K)$ is bigraded, its bar construction is now triply graded.
\begin{cor}
    Given a surjection of augmented dg algebras $C\onto C_0$ whose underlying $A$-modules are torsion free, $\H^*_{Bar}(C_0)$ is computed by the totalization of the tri-complex
    \[
    Bar(C,K) = T\left( \desus \delta K \oplus \desus \ol{C} \right)
    \]
    with $d_1$ and $d_2$ induced by the internal differential $C$ and the product structure, and an additional differential $d_0$ induced by the inclusion $i\maps K \to C$. We call the three gradings: cohomological, weight and defect degree. This tri-complex is functorial in commutative squares
    \[
    \xymatrix{
    C\ar@{->>}[r]^p \ar[d] & C_0 \ar[d] \\
    C'\ar@{->>}[r]^{p'} & C_0'.
    }
    \]
\end{cor}
\begin{proof}
    Since $\tot(C,K)\to C_0$ is a quasi-isomrophism of dg algebras whose terms are $A$-torsion free, Lemma \ref{lem:quasi-iso invariant} shows that the map $\tot Bar(\tot(C,K))\to \tot Bar(C_0)$ is a quasi-isomorphism. But the totalization of the bicomplex $Bar(\tot(C,K))$ is natutrally isomorphic to totalization of the tri-complex $Bar(C,K)$. Functoriality is clear by construction.
\end{proof}

\begin{definition}\label{def:defect}
    Let $C\onto C_0$ now be a surjection of \emph{connected} dg algebras whose underlying $A$-modules are torsion free. The \emph{defect filtration} on $\H^0_{Bar}(C_0)$ is the filtration induced by
    \[
    D_{p} Bar(C,K) = \text{subcomplex of terms with defect degree $\geq -p$}.
    \]
    This is the subcomplex in which at most $p$ tensor factors are taken from $\delta K$. Denote the induced filtration on cohomology by $D_{p}\H^0_{Bar}(C_0)$.
\end{definition}
\begin{proof}[Proof of Theorem \ref{thm:defect filtration}]
Let us consider the spectral sequence associated with the defect filtration on $\tot Bar(C,K) = T(\delta K \oplus \ol{C})$, suppressing the weight desuspension from the notation. Its $E_0$-page consists of terms
\[
E_0^{-p,\bullet,\bullet} = D_p/D_{p-1} = \bigoplus_{i_0,\ldots,i_p\geq 0} \ol{C}^{\otimes i_0} \otimes \delta K \otimes \ol{C}^{\otimes i_1} \otimes \delta K \otimes \ldots \otimes \ol{C}^{\otimes i_{p-1}} \otimes \delta K \otimes \ol{C}^{\otimes i_p}
\]
with differential given as the usual Bar differential $d_1 + d_2$. In particular, when $p=0$ one gets the usual total complex $(\tot Bar(C), d_{Bar})$, so that the $p=0$ part of the $E_1$-page is naturally isomorphic to $H^*_{Bar}(C)$.

Define the defect weight filtration on $E_0^{\bullet,\bullet,\bullet} = (T(\delta K \oplus \ol{C}), d_1+d_2)$ to be the filtration by tensor weight in $T(\delta K \oplus \ol{C})$ -- crucially, this filtration is distinct from the one induced from usual weight filtration on $Bar(C,K) = (T(\delta K \oplus \ol{C}),d_0+d_1+d_2)$, see Remark \ref{rmk:defect weight filtration}. Let us examine the spectral sequence associated with defect weight. To this end, fix $p\geq 0$ and let us denote
$B_p^{\bullet,\bullet} = (E_0^{-p,\bullet,\bullet},d_1+d_2)$.

The associated graded for the defect weight filtration on $B_p^{\bullet,\bullet}$ -- an auxiliary $E'_0$-page -- consists of the same terms but with only the internal differential $E_0' = ( E_0^{-p,\bullet,\bullet}, d_1 )$, so it is a direct sum of tensor products of the complexes $(\ol{C},d)$ and $(\delta K,-d)$. Let us apply the K\"unneth formula to the natural map
\[
\bigoplus_{i_0,\ldots,i_p\geq 0} \H(\ol{C})^{\otimes i_0} \otimes \delta\H(K) \otimes \H(\ol{C})^{\otimes i_1} \otimes \delta\H(K) \otimes \ldots \otimes \H(\ol{C})^{\otimes i_{p-1}} \otimes \delta\H(K) \otimes \H(\ol{C})^{\otimes i_p} \to E_1'.
\]
Recall our connectivity assumption, that $\ol{C}$ has cohomology only in degrees $\geq 1$ and $\H(K)$ is concentrated in degrees $\geq 2$. These imply that tensors of weight $n$ with $p$ many factors in $K$ only have cohomology in cohomological degrees $\geq n+p$. Therefore, on tensors of weight $n$ the K\"unneth map is an isomorphism in cohomological degree $\leq n+p$, where in degrees $<n+p$ the $E_1'$-page vanishes, and in degree $n+p$
\[
\bigoplus_{i_0+\ldots+i_p = n-p} \H^1(\ol{C})^{\otimes i_0} \otimes \delta\H^2(K) \otimes \H^1(\ol{C})^{\otimes i_1} \otimes \ldots \otimes \delta\H^2(K) \otimes \H^1(\ol{C})^{\otimes i_p} \cong (E_1')^{-p,-n,n+p}.
\]
In terms of total degree in $E_1'$, we conclude that $E_1'$ vanishes in total degree $< p$, and in total degree $p$ it is the part of $T(\delta \H^2(K) \oplus \H^1(C) )$ containing  exactly $p$ many factors in $\delta\H^2(K)$.

Diagrammatically, abbreviating $\H = \H^1(\ol{C})$ and $V = \delta\H^2(K)$, the $E_1'$-page takes the form
    \[\xymatrixcolsep{2pc}\xymatrixrowsep{1.5pc}
    \xymatrix{
      0 & \bigoplus V^{\otimes i}|\H|V^{\otimes j}|\H|V^{\otimes p-i-j} \ar[r]^-{d_1}\ar@{-->}[rrd]^-{d_2} &  \# & \# \\
      0 & 0 &  \bigoplus V^{\otimes i}|\H|V^{\otimes p-i} \ar[r] & \bigoplus V^{\otimes i}|\delta\H^3(K)|V^{\otimes p-i-1} \\
      0 & 0 & 0 & V^{\otimes p}
    }
    \]
from which we see that the entire auxiliary spectral sequence $(E'_r,d'_r)$ vanishes total degree $< p$, and no differentials can hit the terms in degree $p$. Thus there is an inclusion $(E'_\infty)^{-q,p+q} \into (E'_1)^{-q,p+q}$, with image characterized by the vanishing of all differentials $d'_r$.

Recalling that the auxiliary spectral sequence computes the cohomology of $\tot B_p = \tot E_0^{-p,\bullet,\bullet}$, which is the $E_1$ page of the original spectral sequence, we conclude that $\tot E_1^{-p,\bullet,\bullet}$ vanishes in total degree $< p$, and in degree $p$ it has associated graded
\[
\bigoplus_{n=0}^{\infty} \operatorname{gr}^F E_1^{-p,-n,p+n} \into \bigoplus_{i_0,\ldots,i_p\geq 0} \H^1(\ol{C})^{\otimes i_0} \otimes \delta \H^2(K) \otimes \ldots \otimes H^1(\ol{C})^{\otimes i_{p-1}} \otimes \delta \H^2(K) \otimes H^1(\ol{C})^{\otimes i_p}.
\]
Summing over all $p\geq 0$ and observing that the defect weight filtration is a filtration of coalgebras, the same argument performed in the category of coalgebras shows that the inclusion
\[
\tot(E_1^{\bullet,\bullet,\bullet})^0 = \bigoplus_{p\geq 0} \operatorname{gr}^W \tot (E_1^{-p,\bullet,\bullet})^{p} \into T(\delta \H^2(K) \oplus \H^1(\ol{C}))
\]
is in fact a homomorphism, respecting defect degree and weight.

Return now to the original spectral sequence. As just proved, its $E_1^{-p,\bullet} = \tot (E_1^{-p,\bullet,\bullet})$ term begins as follows.
    \[\xymatrixcolsep{2pc}\xymatrixrowsep{1.5pc}
    \xymatrix{
     0 & E_1^{-2,2} \ar[r]^-{d_1}\ar@{-->}[rrd]^-{d_2} &  E_1^{-1,2} & \H^2_{Bar}(C) \\
     0 & 0 &  E_1^{-1,1} \ar[r] & \H^1_{Bar}(C) \\
     0 & 0 & 0 & \H^0_{Bar}(C)
    }
    \]
    Since no differentials can hit the diagonal in total degree $0$, one obtains a natural inclusion $E_{\infty}^{-p,p}\into E_1^{-p,p}$, whose image is characterized by the vanishing of all higher differentials.
    
    Summarizing, the defect degree filtration on $\tot Bar(C)$ defines $0 \leq D_0 \leq D_1 \leq \ldots \leq H^0_{Bar}(C_0)$, with $D_0\cong H^0_{Bar}(C)$ and associated graded
    \[
    \operatorname{gr}^D H^0_{Bar}(D_0) \cong \bigoplus_{p\geq 0} E_\infty^{-p,p} \into \bigoplus_{p\geq 0} E_1^{-p,p}
    \]
    as graded coalgebras. Furthermore, the $E_1$-coalgebra is filtered by the defect weight $W_\bullet E_1$, so that the associated graded with respect to $W_\bullet$ gives an inclusion of coalgebras
    \[
    \defect\maps \operatorname{gr}^W E_1 \into T(\delta \H^2(K) \oplus \H^1(C))
    \]
    respecting both defect weight and degree. We call the last bigraded coalgebra homomorphism the \emph{defect map}. The restriction of the defect weight filtration and the defect map to $E_{\infty} \leq E_{1}$ is the structure we set out to construct.
    
    Clearly, when $p=0$ the defect weight filtration is precisely the usual weight filtration on $\H^0_{Bar}(C)$ and the defect map is the inclusion recording the leading term as given by the $(n-1)$-fold reduced comultiplication on $\H^0_{Bar}(C)_{\leq n}$.
    
    Lastly, if $\alpha \in \H^1(C_0) \leq \H^0_{Bar}(C_0)$ has defect degree $1$, then it is not the image of a cocycle in $C$. But $p\maps C\onto C_0$ is surjective, so there does exist a lift $\tilde{\alpha}\in C^1$, which then necessarily has $0 \neq d\tilde{\alpha} \in K$. The sum
    \[
    \xymatrixcolsep{1ex}\xymatrixrowsep{1ex}\xymatrix{
    \delta\cdot d\tilde{\alpha} & \\
    + & \tilde{\alpha}
    }
    \]
    is a cocycle in $Bar(C,K)$ projecting to $\alpha$. Its image in $\operatorname{gr}^D_1\H^0_{Bar}(C_0)$ is therefore $[d\tilde{\alpha}]$, which is then also its defect. But this construction is the connecting homomorphism of the long exact sequence for $0\to K\to C\to C_0\to 0$, as claimed.
    \qedhere
\end{proof}
\begin{remark}\label{rmk:defect weight filtration}
    An explicit description of the defect weight filtration on $Bar(C,K)$ is that 
    \[
    W_{n,p} Bar(C,K) = \{\text{highest defect summand has $\leq p$ many factors in $\delta K$ and weight $\leq n$}\}.
    \]
    This differs of the weight filtration, that bounds the weights of all tensors in a cochain, and not just the top defect summand. The defect map returns the cohomology class of this leading term of defect $p$ and weight $n$.
\end{remark}

Shifting back to our topological setup of $C = C^*(Y_{|R};A)$ projecting onto $C_0 = C^*(Y;A)$, which satisfies Hypothesis \ref{hyp:defect setup} with $\H^2(K) \cong A^R$ giving functions on relations, the defect map is a bigraded coalgebra homomorphism
\begin{equation}
    \defect\maps \operatorname{gr}^W \operatorname{gr}^D \H^0_{Bar}(\Gamma_0;A) \longto T( \delta A^R \oplus \H^1(\Gamma;A) ).
\end{equation}
We claim that this defect map records nonvanishing letter-braiding invariants on relations $r\in R$.
\begin{example}[The universal relation]\label{ex:universal relation}
    Consider the case of the free group on one generator $F_{\{r\}}$ to which we introduce the relation $r$, making the group trivial. The associated topological construction is the inclusion $S^1(r) \into D^2(r)$.

    Since $\H^1(D^2) = 0$, the defect map takes values in the coalgebra $T( \delta\H^2(D^2,S^1) )$, in which weight and defect degree coincide. It then follows that the defect weight filtration is trivial in that it matches the defect degree. Since $\H^2(D^2,S^1)= A^{\{ r \}}$ is generated by $(\tau\maps r\mapsto 1)$ and the connecting homomorphism $\H^1(S^1) \to \H^2(D^2,S^1)$ sends the generator $t\in \H^1(S^1)$ to $\tau$, the defect filtration splits and is the evident isomorphism
    $\H^0_{Bar}(S^1;A) \cong A[t] \to A[\delta\tau]$, with $t^n$ having defect degree $n$.

    This gives an interpretation of defect: it records the letter-braiding invariants of the relation $r$. 
    The fact that every nontrivial braiding invariant is nonzero on $r$ reflects the fact that none are pulled back from the (trivial) quotient group.
\end{example}
This interpretation will be made precise in general. Indeed, suppressing $\delta$ from the notation, the defect map takes values in $T(A^R \oplus \H^1(\Gamma;A))$, but
\[
\bigoplus_{n\geq 0} ( A^R \oplus  \H^1(\Gamma;A) )^{\otimes n} \into \bigoplus_{n\geq 0} ( A^R \oplus  A^\Gamma )^{\otimes n} \into \bigoplus_{n\geq 0} A^{ (R \coprod \Gamma)^n  }
\]
defining $A$-valued $n$-ary functions on the relations $R$ and the quotient group $\Gamma$. We lift these to functions $\Gamma_0 \to A$ by setting for every $w\in \Gamma_0$
\begin{equation}\label{eq:extending to gamma0}
f\in A^\Gamma \rightsquigarrow \langle f,w\rangle = f(p(w)) \quad\text{and}\quad g\in A^R \rightsquigarrow \langle g,w\rangle = \begin{cases}
    g(r) & w=r\in R \\
    0 & \other.
\end{cases}    
\end{equation}
\begin{prop}[\textbf{Defect evaluation}] \label{prop:evaluation defect comparison}
    Consider the defect map of $\H^0_{Bar}(\Gamma_0;A)$ obtained from the surjection of dg algebras $C^*(Y_{|R};A) \onto  C^*(Y;A)$. If $T\in \ol{\H^0_{Bar}(\Gamma_0;A)}$ has defect degree $p$ and defect weight $n$, then for every sequence of group elements $w_1,\ldots ,w_m$ such that $q \geq p$ of them that are relations $r_i\in R$,
\begin{equation}\label{eq:evaluation on defect}
    \ell_{T} (w_1|\ldots|w_n) = \begin{cases}
        \langle \defect_{n,p}(T) \,,\, w_1 \otimes \ldots \otimes w_n \rangle & p=q \text{ and } n=m \\
        0 & p<q \text{ or } n< m,
    \end{cases}
\end{equation} 
where $\langle - , \bullet \rangle \maps (A^R \oplus \H^1(\Gamma;A))\otimes A[\Gamma_0] \to A$ extends the evaluations in \eqref{eq:extending to gamma0}. We allow $m<n$ but make no claims about its evaluation when $p=q$.
\end{prop}
\begin{proof}
    Lets first consider the case $p = 0$. In this case $T$ is pulled back from some $ \tilde{T} \in \H^0_{Bar}(\Gamma;A)$. By push-pull \eqref{eq:push-pull}
    \[
    \ell_{T}(w_1|\ldots |w_m) = (\ell_{\bullet}(w_1)\otimes\ldots\otimes \ell_{\bullet}(w_m))_{\ol\Delta^{m-1}(p^*\tilde{T})}
 = (\ell_{\bullet}(p(w_1))\otimes\ldots\otimes \ell_{\bullet}(p(w_m)))_{\ol\Delta^{m-1}(\tilde{T})}
    \]
    so since $\defect_{n,0}(T) = \ol{\Delta}^{n-1}(\tilde{T}) $, the two evaluations in \eqref{eq:evaluation on defect} coincide when $m=n$. In particular, when a positive number of $w_i$ belong to $R$, some of these evaluations vanish and thus so does their product. Also, when $m>n$ we have $\ol{\Delta}^{m-1}(T)=0$, as needed for \eqref{eq:evaluation on defect} to hold. This concludes the proof for all cases $(n,p) = (n,0)$, and we subsequently assume that $p\geq 1$.

    Now, for vanishing in the case $q>p$, assume by induction on $p$ that for all $p'<p$, whenever $T'$ has defect degree $\leq p'$ and $w_1,\ldots,w_{m'} \in \Gamma_0$ contains $q'>p'$ many relations $r_i\in R$ then
    \[
    \ell_{T'}(w_1|\ldots|w_{m'}) = 0.
    \]
    Let $T$ have defect degree $p$ and $w_1,\ldots,w_m\in \Gamma_0$ contain $q>p$ many relations. Say $w_i=r$ is the first relation, then Proposition \ref{prop:properties of multi-eval} shows
    \[
    \ell_{T}(w_1|\ldots|r|w_{i+1}|\ldots|w_m) = (\ell_{\bullet}(w_1|\ldots|r)\otimes \ell_{\bullet}(w_{i+1}|\ldots|w_m))_{\ol{\Delta}(T)}.
    \]
    By the first case considered, $\ell_{\bullet}(w_1|\ldots|r)$ vanishes on defect degree $0$, so the only nontrivial contributions to this evaluations come from summands of $\ol{\Delta}(T)$ in which the left-hand factor has defect degree $>0$. But since defect degree is a coalgebra filtration, all those summands must have right-hand factor of defect degree $< p$. Since there are $q-1\geq p$ many relations among $w_{i+1},\ldots,w_m$, the induction hypothesis implies that those evaluations vanish as well.
    
    We thus proved that the multi-evaluation vanishes whenever the number of relations in the evaluation is greater than the defect degree. In particular, if $w_1,\ldots,w_n$ contains $p$ many relations, then $\ell_{\bullet}(w_1|\ldots|w_n)$ vanishes on $D_{p-1}$ and factors through the associated graded 
    \[
    \ell_{\bullet}(w_1|\ldots|w_n)\maps \operatorname{gr}^D_p \H^0_{Bar}(\Gamma_0;A) \to A.
    \]
    Similarly, $(\ell_{\bullet}(w_1|\ldots|w_i)\otimes \ell_{\bullet}(w_{i+1}|\ldots|w_m)) \maps \H^0_{Bar}(\Gamma_0;A)\otimes \H^0_{Bar}(\Gamma_0;A) \to A$ vanishes on all tensors $D_j\otimes D_{k}$ with $j+k< p$. Thus the relation $(\ell_{\bullet}(w_1|\ldots|w_i)\otimes \ell_{\bullet}(w_{i+1}|\ldots|w_n))_{\ol{\Delta}(\bullet)} = \ell_{\bullet}(w_1|\ldots|w_n)$ descends to the associated graded quotients as well.

    We next prove that evaluations vanish when $m>n$. First, suppose $T$ has defect degree $p$ and defect weight $n=1$. Then $(\defect\otimes \defect)\circ \ol{\Delta}(T) = \ol{\Delta}(\defect(T)) = 0$, so by injectivity of the defect map, it must be that $\ol{\Delta}(T)$ has total defect degree $< p$. But then
    \[
    \ell_{T}(w_1|\ldots|w_m) = (\ell_{\bullet}(w_1)\otimes \ell_{\bullet}(w_2|\ldots|w_m))_{ \ol{\Delta}(T)} = 0
    \]
    from the previous discussion. The same argument now shows the case of general $m>n\geq 1$ by induction on $n$.

    We are left with the cases $m=n$ and $q=p$, where we want to show the defect map determines the multi-evaluation. Firstly, the vanishing of $\ell_{\bullet}(w_1|\ldots|w_n)$ on defect weight $W_{n-1}\operatorname{gr}^D_p \H^0_{Bar}(\Gamma_0;A)$ shows that these evaluations factor through the associated graded
    \[
    \ell_{\bullet}(w_1|\ldots|w_n)\maps \operatorname{gr}^W_n \operatorname{gr}^{D}_p \H^0_{Bar}(\Gamma_0)\to A
    \]
    on which the defect map is defined and injective. We prove \eqref{eq:evaluation on defect} by double induction on $(n,p)$. The cases with $(n,0)$ were already proved at the top of this proof.

    For $(n,p) = (1,1)$, we must calculate $\ell_T(r)$ for $r\in R$. The relation is witnessed topologically by a map of pairs $r\maps (D^2(r),S^1(r)) \to (Y_{|R},Y)$. This map induces a commutative square of cochain algebras
    \[
    \xymatrix{
    C^*(Y_{|R};A) \ar@{->>}[r] \ar[d]_{r^*} & C^*(Y;A) \ar[d]_{r^*} \\
    C^*(D^2(r);A) \ar@{->>}[r] & C^*(S^1(r);A)
    }
    \]
    showing that the pullback $r^*\maps \H^0_{Bar}(Y;A) \to \H^0_{Bar}(S^1(r);A)$ respects the respective defect structures. For the pair $(S^1\into D^2)$ we already computed the defect map in Example \ref{ex:universal relation}, so we are left with a commutative diagram
    \[
    \xymatrix{
\operatorname{gr}^W_1\operatorname{gr}^D\H^0_{Bar}(Y;A) \ar[r]^-{r^*} \ar[d]_{\defect_{1,1}} & \operatorname{gr}^W_1\operatorname{gr}^D\H^0_{Bar}(S^1;A) \ar[d]_{\defect_{1,1}} \ar[r]^-{\cong} & At \ar[d]_{\cong} \ar[r]^\cong & A \\ 
    \H^2(Y_{|R}, Y ; A) \ar[r]^-{r^*} \ar[d]_{\cong} &  \H^2(D^2,S^1;A) \ar[r]^-\cong \ar[d]_\cong & A\tau  & \\
    A^R \ar[r]^{\operatorname{ev}_r} & A & &
    }
    \]
    The composition the top row is the letter-braiding number $\ell_\bullet(r)$, while the bottom arrow identifies this number with the evaluation of functions $A^R$ at $r$. This proves the claim when $(n,p)=(1,1)$.

    Next, suppose by induction that \eqref{eq:evaluation on defect} holds for all $(n',p')$ such that $n'<n$ or $p'<p$. Then $\ell_{T}(w_1|\ldots|w_n) = (\ell_{\bullet}(w_1)\otimes\ell_{\bullet}(w_2|\ldots|w_n))_{\ol{\Delta}(T)}$ and $\ol{\Delta}(T)$ is a sum of tensors $T_1\otimes T_2$, each of which either have total defect degree $<p$, or it has total defect degree $p$ but respective defect weights $n_1$ and $n_2$ such that $n_1+n_2 = n$. In the former case the evaluation vanishes since there are $p$ relations in the list, while in the latter case the induction hypothesis gives
    \[
    \ell_{T_1}(w_1)\otimes \ell_{T_2}(w_2|\ldots|w_n) = \langle \defect(T_1) \otimes \defect(T_2) \,,\, w_1\otimes\ldots\otimes w_n \rangle
    \]
    Equality in \eqref{eq:evaluation on defect} thus holds since the defect map is a coalgebra homomorphism.
\end{proof}
With this we can complete the proof of our lifting criterion.
\begin{proof}[Proof of Theorem \ref{thm:lifting criterion} \eqref{case:vanishing multi-evaluations}$\implies$\eqref{case:lifting}]
    Suppose that $T\in \H^0_{Bar}(\Gamma_0;A)$ is not pulled back from $\Gamma$. Then by Theorem \ref{thm:defect filtration} it has defect degree $p>0$. Let $n$ denote the defect weight of $T$ modulo $D_{p-1}$. Then,
    \[
    0 \neq \defect_{n,p}(T) \in ( A^R \oplus \H^1(\Gamma;A) )^{\otimes n}
    \]
    a linear combination of tensors with exactly $p$ factors taken from $A^R$. As discussed in this subsection, the tensor power above injects into the space of $n$-ary functions on $R \coprod \Gamma$, and it does not vanish when evaluated on some sequence $w_1,\ldots,w_n$ with $p$ many entries in $R$. Now, the evaluation formula \eqref{eq:evaluation on defect} shows that 
    \[
    \ell_{T}(w_1|\ldots|w_n) = \langle \defect_{n,p}(T) \,,\, w_1\otimes\ldots\otimes w_n\rangle \neq 0. \qedhere
    \]
\end{proof}
\subsection{Consequences of the lifting criterion}\label{sec:consequences of lifting}
    There is a nice categorical restatement of the lifting criterion of Theorem \ref{thm:lifting criterion}.
    \begin{cor}[\textbf{Lifting criterion v2}] \label{cor:cartesian diagram}\label{cor:cartesian}
        Let $p\maps \Gamma_0\onto \Gamma$ be a surjective group homomorphism, then the commutative diagram of inclusions induced by the letter-braiding pairing
        \begin{equation*}\xymatrix{
            \H^0_{Bar}(\Gamma;A) \ar@{^(->}[r] \ar@{^(->}[d]^{p^*} & A[\Gamma]^\star \ar@{^(->}[d]^{p^*} \\
            \H^0_{Bar}(\Gamma_0;A) \ar@{^(->}[r] & A[\Gamma_0]^*}
        \end{equation*}
        is Cartesian. That is, a linear functional on $A[\Gamma]$ is realized by letter-braiding invariants on $\Gamma$ iff its pullback to $A[\Gamma_0]$ is realized by letter-braiding invariants on $\Gamma_0$.
    \end{cor}
    \begin{proof}
        This is immediate from the implication \eqref{case:lifting}$\iff$\eqref{case:group ring factors} of Theorem \ref{thm:lifting criterion}.
    \end{proof}
In the special case of a presentation $\Gamma = \langle S|R \rangle$, one can be more concrete.
\begin{cor}\label{cor:subcoalg of free groups}
    Let $\Gamma = \langle S|R\rangle$ be a group presentation (possibly infinite). Then the on bar cohomology of $\Gamma$ is the maximal sub-coalgebra
        \[
        \H^0_{Bar}(\Gamma;A) \leq T(A^S)
        \]
    of possibly nonhomogeneous tensors of forms in $A^S$, characterized by the property that its elements have vanishing braiding numbers of all relations $r\in R$.
\end{cor}
\begin{proof}
    The image of $\H^0_{Bar}(\Gamma;A) \into T(A^S)$ is clearly a sub-coalgebra and all its elements have vanishing letter-braiding numbers on $r\in R$. To see that it is the maximal such, if $H \leq T(A^S)$ is a coalgebra with this property, then every $T\in H$ satisfies $\ell_{T}(w_1|\ldots|w_n) = 0$ whenever $w_1,\ldots,w_n\in F_S$ is a sequence of elements at least one of which is an element in $R$. But the lifting criterion now shows that $T$ comes from $\H^0_{Bar}(\Gamma;A)$.
\end{proof}
With these facts establishes, we can prove the letter-braiding invariants are complete.
\begin{proof}[Proof of Theorem \ref{thm:completeness}]
    For a finitely generated free group $F_r$ generated by $\{x_1,\ldots,x_r\}$, we claim that the truncated group ring $A[F_r]/I^n$ is the free $A$-module generated by products
    \[
    (x_{i_1}-1)(x_{i_2}-1)\cdots (x_{i_d}-1)
    \]
    for sequences $(i_1,\ldots,i_d)$ of length $0\leq d < n$. Indeed, every generator $x_i$ can be written as $1+ (x_i-1)$, its inverse is $x_i^{-1} = 1-(x_i-1) + (x_i-1)^2 - \ldots +(-1)^{n-1}(x_i-1)^{n-1}+\ldots$, and every element in the free group can be written as a product of $x_i^{\pm 1}$ so expressible as a linear combination of shifted monomials. Freeness will follow from the following observation.
    
    Let $X_i\maps F_r \to A$ be the linear functional that counts instances of $x_i$ with sign, so that $\H^1(F_r;A) = \langle X_1,\ldots,X_r \rangle$. Then the letter-braiding invariants of $F_r$ are labeled by the tensor coalgebra $T(  X_1,\ldots,X_r  )$ of noncommutative polynomials in the $X_i$'s. Proposition \ref{prop:properties of multi-eval}
    \[
    \langle X_{i_1}|\ldots|X_{i_d} \,,\, (x_{j_1}-1)\cdots (x_{j_d}-1) \rangle = \ell_{x_{j_1}}(X_{i_1})\ldots \ell_{x_{j_d}}(X_{i_d}) = \delta_{(i_1,\ldots,i_d),(j_1,\ldots,j_d)}
    \]
    gives a dual basis to degree $d$ monomials, and vanishes on all higher degrees. By induction on degree, one can find appropriate linear combinations of weight $\leq n$ that define a dual basis to the monomials. This proves that letter-braiding $\H^0_{Bar}(F_r;A)_{\leq n}\into (A[F_r]/I^n)^*$ is surjective and thus an isomorphism (it also gives a recipe for constructing the Taylor-Maclaurin expansion of every group ring element, i.e. writing it as a linear combination of monomials).
    
    Now, let $F_n\onto \Gamma$ be a presentation of some group. Then the Cartesian square in Corollary \ref{cor:cartesian} shows that $\H^0_{Bar}(\Gamma;A)_{\leq n}\to (A[\Gamma]/I^n)^*$ is also an isomorphism.
\end{proof}
Lastly, we consider the case of infinitely generated groups.
\begin{proof}[Proof of Theorem \ref{intro-thm:infinite type universal pairing}]
    First, we show that the letter-braiding pairing is a universal conilpotent multiplicative pairing as defined in Definition \ref{def:universal pairing}. Suppose $(C,1_C)$ is a conilpotent coalgebra equipped with a multiplicative pairing
    \[
    \langle -,\bullet\rangle_C: C\otimes A[\Gamma]\to A.
    \]

    Pick any presentation $\pi:F_S\onto \Gamma$ so that $\pi$ defines an injection $\H^0_{Bar}(\Gamma;R)\into \H^0_{Bar}(F_S;A)\cong T(A^S)$ whose image is characterized by Corollary \ref{cor:subcoalg of free groups}. Define an $A$-linear map $C\to A^S$ by sending every $f\in C$ to the functional
    \[
    s\mapsto \langle f,\pi(s)\rangle_C.
    \]
    Since $T(A^S)$ is the cofree coalgebra on $A^S$ and $C$ is conilpotent, the map $p$ extends uniquely to a homomorphism of coaugmented coalgebras $h:C\to T(A^S)$.
    
    We claim that the two obvious pairings $C\otimes A[F_S]\to A$ coincide:
    \begin{equation}\label{eq:universal pullback}
    \langle f,\pi(w)\rangle_C = \langle h(f),w\rangle_{U}
    \end{equation}
    for all words $w\in F_S$ and $f\in C$, where $\langle-,\bullet\rangle_U:T(A^S)\otimes A[F_S]\to A$ is the letter-braiding pairing. Indeed, by the multiplicativity assumption on $\langle-,\bullet\rangle_C$, the element $1_C$ induces the same functional as $1\in T(A^S)$, so since $h$ is a coalgebra homomorphism both pairings are multiplicative. One can readily show that multiplicative pairings are determined by their values on a generating set, e.g. by induction on weight and word-length, so since the two pairings above agree on the generating set $S\subseteq F_S$ by construction, they agree on all words $w\in F_S$. These words form a basis for $A[F_S]$, so the two pairings coincide.

    By Corollary \ref{cor:subcoalg of free groups}, the coalgebra $\H^0_{Bar}(\Gamma;A)\leq T(A^S)$ is the maximal subcoalgebra whose letter-braiding pairing gives well-defined functions on $\Gamma$, so the homomorphism $h$ factors uniquely through $C\to \H^0_{Bar}(\Gamma;A)$. This homomorphism is uniquely determined by Equation \eqref{eq:universal pullback} since its restriction to generators $s\in S$ already determined $h$.
    
    It remains to characterize the conilpotent dual of $\Gamma$ as the maximal submodule $M\subseteq A[\Gamma]^*$ on which pullback along the convolution product $m:A[\Gamma]\otimes A[\Gamma]\to A[\Gamma]$ factors through a coproduct $M\to M\otimes M\subseteq A[\Gamma]^*\otimes A[\Gamma]^* \to (A[\Gamma]\otimes A[\Gamma])^*$. Observe that since $\sum (M_\alpha\otimes M_\alpha)\subseteq (\sum M_\alpha)\otimes (\sum M_\alpha)$, the sum of submodules with this property again has this property, thus there does exist a unique maximal submodule of this kind. Clearly, any multiplicative pairing between $C$ and $\Gamma$ defines a unique homomorphism $C\to M$ that pulls-back the pairing.

    Now, let $M_{conilp}\subseteq M$ denote the submodule of elements in $M$ annihilated by some iterate of the reduced coproduct. Then the pairing
    \[
    \langle -,\bullet \rangle: M_{conilp}\otimes A[\Gamma]\to A
    \]
    obtained by restricting the evaluation $A[\Gamma]^*\otimes A[\Gamma]\to A$ clearly defines a multiplicative conilpotent pairing, which is clearly universal. It follows that the letter-braiding pairing defines an isomorphism of coaugmented coalgebras $\H^0_{Bar}(\Gamma;A)\xrightarrow{\sim} M_{conilp}$.
\end{proof}

\subsection{Relative Massey products measure relations}\label{sec:relative Massey}
In Lemma \ref{lem:Bar spectral sequence} we saw that the obstructions to having a tensor $T(\H^1(\Gamma;A))$ define a letter-braiding invariant of $\Gamma$ is that all its cup and Massey products must vanish. In the last section we gave an alternative characterization as having vanishing defect, as measured by evaluations on relations. The relationship between Massey products and relations in a group has a long history, going back to the work of Dwyer \cite{Dwyer} and Fenn--Sjerve \cite{other-Massey-paper}), but they only cover a somewhat narrow set of examples: ones in which the Massey product of cohomology classes $\langle\alpha_1,\ldots,\alpha_n\rangle$ is defined, and it pairs with relations that are commutators.

However, even the most elementary examples we considered above, of surface groups and braid groups (Examples \ref{ex:surface grps},\ref{ex:braids}), exhibit interesting invariants in which one must consider a linear combination of tensors from $\H^1(\Gamma;A)$ to detect various group elements. In this section we expand the scope of Dwyer's analysis, giving a full picture of how Massey products coincide with letter-braiding numbers of relations, including non-commutators, via a \emph{relative Massey product}.

Let $\Gamma = \langle S|R\rangle$ be a group presentation. Applying the construction $Y\mapsto Y_{|R}$ from \S\ref{sec:introducing relations} to the wedge of circles $X_S$ from \S\ref{sec:free groups}, we obtain a presentation $2$-complex $X_{S|R}$ whose fundamental group is $\Gamma$. It contains a subcomplex $p\maps X_S\into X_{S|R}$, inducing the surjection $F_S \onto \Gamma$. Observe that the $1$-simplices of $X_S$ are in bijection with the set of generators $S\subseteq \Gamma$, while $X_{S|R}$ has additional edges coming from the relations. $1$-cochains are $A$-valued functions on those edges.
\begin{observation}
Restricting $1$-cochains of $X_{S|R}$ to the generators $S\subseteq \Gamma$ gives a surjection of coalgebras $p^* \maps Bar(X_{S|R};A) \onto Bar(X_S;A) = T(A^S)$, where the target is the coalgebra of letter-braiding invariants on $F_S$ from \S\ref{sec:free groups}. So every letter-braiding invariant of $F_S$ lifts to a bar cochain of $\Gamma$ by extending the functions from $X_S$ to the interiors of relation discs $D^2(r)\subseteq X_{S|R}$. But those extensions will typically not be $d_{Bar}$-closed: their coboundary involves a $2$-cochain of $X_{S|R}$, detecting some of the relations.
Let us quantify the obstruction to closed lifting.
\end{observation}
Consider then the short exact sequence
\begin{equation}\label{eq:bar short exact sequence}
0 \to (K)_{Bar} \to Bar(X_{S|R};A) \to T(A^S) \to 0    
\end{equation}
where $(K)_{Bar}$ is a coalgebraic coideal. In analogy with relative cochains $C^*(X_{S|R},X_S;A)$, denote this coideal by $Bar(X_{S|R},X_S;A)$. More generally,
\begin{definition}
    When $p\maps X\into Y$ is a subcomplex inclusion, denote
    \[
    Bar(Y,X;A) = \ker\left( p^*\maps Bar(Y;A)\onto Bar(X;A)  \right),
    \]
    the coideal cogenerated by $K=C^*(Y,X;A)$ in $T(\ol{C^*(Y;A)})$. Its cohomology will be denoted $\H^*_{Bar}(Y,X;A)$.
\end{definition}

    An alternative equivalent definition proceeds by recalling that a surjection $Bar(C)\onto Bar(C_0)$ is quasi-isomorphic to the injection $Bar(C)\into Bar(C,K)$. Then $Bar(Y,X;A)$ could be defined as $\coker(Bar(C)\into Bar(C,K))$ -- the quotient by the defect $0$ subcomplex -- up to a degree shift.

\begin{prop}\label{prop:relative Massey}
    An inclusion $X\into Y$ defines an exact sequence 
    \[
    \H^0_{Bar}(Y;A) \to \H^0_{Bar}(X;A) \xrightarrow{\mu} \H^1_{Bar}(Y,X;A),
    \]
    with $\mu(-)$ measuring the obstruction of lifting bar cocycles from $X$ to $Y$. We call $\mu(-)$ the \emph{relative Massey product}, for reasons explained below.
    
    In particular, a tensor $T\in T(A^S)$ extends to a closed cochain in $Bar^0(X_{S|R};A)$ if and only if its relative Massey product $\mu(T)\in \H_{Bar}^1(X_{S|R},X_S;A)$ vanishes.
\end{prop}
\begin{proof}
    The snake lemma applied to 
    \[
    0\to Bar(Y,X;A)\to Bar(Y;A) \to Bar(X;A) \to 0
    \]
    gives the claimed exact sequence. It follows that $T\in T(A^S)$ lifts to a class $\tilde{T}\in \H^0_{Bar}(X_{S|R};A)$ if and only if $\mu(T)=0$.
\end{proof}
    
For $T\in \H^0_{Bar}(X;A)$ its the relative Massey product is computed by picking any lift $\tilde{T}\in Bar^0(Y;A)$ and computing its Bar differential:
    \[
    \mu(T) = [d_{Bar}(\tilde{T})]\in \H^1_{Bar}(Y,X;A).
    \]
When $\mu(T)$ has weight $1$, it is a relative cohomology class $\mu(T) \in \H^2(Y,X;A)$ whose image in the absolute cohomology $\H^2(Y;A)$ is the usual Massey product of the lift $\tilde{T}$, as explained in Proposition \ref{prop:global Massey}. Thus $\mu(T)$ is a refinement of the Massey product that remembers that $T$ came from $X$.

To see how the relative Massey product detects relations in a group, return to the special case $(X_{S|R},X_S)$. Then $\mu(T)$ takes values in $\H^2(X_{S|R},X_S;A) \cong A^R$, and is therefore exactly a function on relations.
\begin{prop}[\textbf{Relative Massey = letter-braiding of relations}] \label{prop:massey=braiding}
    Suppose $T\in T(A^S)$ is a letter-braiding invariant of $F_S$ whose relative Massey product has weight $1$, that is $\mu(T)\in \H^2(X_{S|R},X_S;A) \cong A^R$. Then for every relation $r\in R$,
    \begin{equation}\label{eq:massey is braiding}
        \ell_T(r) = \langle\mu(T),[D^2(r)]\rangle,
    \end{equation}
    where $[D^2(r)]\in H_2(X_{S|R},X_S;A)$ is the relative homology class of the relation disc of $r$.
\end{prop}
\begin{proof}
    This is essentially a special case of Proposition \ref{prop:evaluation defect comparison}. With the interpretation of $\H^*_{Bar}(Y,X;A)$ using the defect filtration as in the paragraph preceding Proposition \ref{prop:relative Massey}, the condition on $\mu(T)$ is equivalent to saying that $T$ has defect $1$. The proof of \eqref{eq:evaluation on defect} here works, using the case of the universal relation $(D^2,S^1)$: the relative Massey product in the exact sequence
    \[
    \cancel{\H^0_{Bar}(D^2;A)} \to \H^0_{Bar}(S^1;A) \xrightarrow{\mu} \H^1_{Bar}(D^2,S^1;A) \to \cancel{\H^1_{Bar}(D^2;A)}
    \]
    is an isomorphism, identifying the letter-braiding invariant with the evaluation on the relative fundamental class $[D^2]\in \H^2(D^2,S^1;A)$.
\end{proof}
The absolute version of this proposition follows as a special case, and it appears e.g. in \cite[Main theorem]{other-Massey-paper}. They evaluate the absolute Massey product $\mu_n\maps \H^1(\Gamma;\Z)^{\otimes n}\to \H^2(\Gamma;\Z)$ on $H_2(\Gamma)$, which by Hopf's formula is generated by $R\cap [\Gamma,\Gamma]$ -- relations that happen to also be commutators. These correspond to disc classes $[D^2(r)]$ that come from $\H^2(X_{S|R};\Z)$ as they have vanishing boundary. Fenn--Sjerve show that the evaluation coincides with Magnus expansion coefficients, which by Theorem \ref{thm:magnus} is the same as $\ell_T(r)$. Proposition \ref{prop:massey=braiding} therefore subsumes their result and extend it in the following directions:
\begin{itemize}
    \item While older work considers Massey products $\langle \alpha_1,\ldots,\alpha_n\rangle$, which are $\mu_n(\alpha_1|\ldots|\alpha_n)$ on rank $1$ tensors, many groups require considering linear combinations of tensors, which our theory accommodates -- see Examples \ref{ex:surface grps},\ref{ex:braids}.
    \item Our generalization sees the entire set of relations $R$, and also considers products of cochains $C^1(\Gamma;A)^{\otimes n}$ that are no a priori closed. In fact, the pairing of $\alpha\in C^1(\Gamma;A)$ with a non-commutator $r \in R$ is exactly the obstruction to $\alpha$ being closed, with \eqref{eq:massey is braiding} specializing to
    \[
    \ell_\alpha(r) = d\alpha([D^2(r)])
    \]
    so $d\alpha=0$ iff $\ell_\alpha(r)=0$ for all $r\in R$.
    \item While we only discussed $\mu(T)$ of weight $1$ above, higher weights appear and are significant obstructions to lifting as well. This is necessary because $T(A^S)$ is a coalgebra, and one might attempt to lift braiding invariants $T$ that split into several factors with nontrivial relative Massey products. This situation is quantified by the defect map and Proposition \ref{prop:evaluation defect comparison}.
\end{itemize}

\section{Topological applications} \label{sec:applications}
So far, topology provided us with invariants of groups. In the other direction, we get a constructive proof the following fact.
\subsection{Cochains determine fundamental group ring over fields}
\begin{theorem}
    Let $X$ be any connected pointed CW-complex (or simplicial set) with a finitely generated fundamental group. Then with coefficients in any field $A=k$, the $2$-truncation of its singular (or simiplicial) cochain algebra $C^{*\leq 2}_{\cup}(X;k)$, along with its augmentation, determines the completed group ring $k[\pi_1(X,*)]_I^{\wedge}$ explicitly via letter-braiding invariants.
\end{theorem}
\begin{proof}
    Given a filtered coalgebra, its dual is naturally a completed algebra, and this operation is inverse to the continuous dual (see e.g. \cite[Subsection (2.5.3)]{sam-snowden}). In light of Theorem \ref{thm:completeness}, it follows that $k[\pi_1(X,*)]_I^{\wedge}$ is isomorphic to the complete algebra of $k$-linear functionals on $\H^0_{Bar}(X;k)$, where their pairing is the letter-braiding pairing. But the latter coalgebra is computed from $C^*_{\cup}(X;k)$ in degrees $\leq 2$ along with its augmentation coming form evaluation at the basepoint.
\end{proof}
    To get ones hands on the (completion of the) group $\pi_1(X,*)$ itself, additional structure is needed: the diagonal $X\to X\times X$, inducing the Hopf algebra structure on $A[\pi_1(X,*)]$. This involves the $E_\infty$-algebra structure on $C^*(X;A)$, and is not present in our discussion.

\subsection{Johnson homomorphisms} \label{sec:Johnson}
Let us construct the Johnson homomorphism mentioned in Theorem \ref{thm:intro-johnson}. Fix a group $\Gamma$ and let $A$ be a PID. Denote $\H = \H^1(\Gamma;A)$.
\begin{definition}
    The \emph{Johnson filtration modulo $A$} of $\Aut(\Gamma)$ consists of subgroups 
    \[\ldots \leq J^A(k) \leq J^A(k-1) \leq \ldots \leq J^A(0) = \Aut(\Gamma)
    \]
    where $J^A(k)$ is the pointwise stabilizer of the quotient ring $A[\Gamma]/I^{k+1}$.
\end{definition}
Since $\H^0_{Bar}(\Gamma;A)_{\leq k}$ injects into the functionals $\Hom_A(A[\Gamma]/I^{k+1}, A )$, it follows that any $\varphi\in J^A(k)$ defines a pullback $\varphi^*\maps \H^0_{Bar}(\Gamma;A)_{\leq k+1} \to \H^0_{Bar}(\Gamma;A)_{\leq k+1}$ fixing $\H^0_{Bar}(\Gamma;A)_{\leq k}$.
Therefore, for every $\varphi\in J^A(k)$ the assignment
\[
\tau_\varphi\maps  T \longmapsto T- \varphi^*(T)
\]
induces is an $A$-linear map 
\[
\tau_\varphi\maps \sfrac{\H^0_{Bar}(\Gamma;A)_{\leq k+1}}{\H^0_{Bar}(\Gamma;A)_{\leq k}} \to \H^0_{Bar}(\Gamma;A)_{\leq k+1},
\]
and by the weight spectral sequence (Lemma \ref{lem:Bar spectral sequence}), its domain is naturally isomorphic to $K_{k+1} \leq \H^1(\Gamma;A)^{\otimes k+1}$ -- the locus of definition and vanishing of the $(k+1)$-st Massey product.

\begin{prop}
    The image of $\tau_\varphi$ has weight $1$, and is therefore an element of $\H^1(\Gamma;A)$.
\end{prop}
\begin{proof}
    It suffices to prove that $T-\varphi^*(T)$ is killed by the reduced coproduct $\ol\Delta$.
    
    Suppose $\ol\Delta(T) = \sum T_1^i\otimes T_2^i$, then
    \[
    \ol\Delta(T-\varphi^*(T)) = \sum T_1^i\otimes T_2^i - \varphi^*(T_1^i)\otimes \varphi^*(T_2^i) = \sum (T_1^i-\varphi^*(T_1^i))\otimes T_2^i + \varphi^*(T_1^i)\otimes (T_2-\varphi^*(T_2^i)).
    \]
    But since both $T_1^i$ and $T_2^i$ have weight $\leq k$, they are fixed by $\varphi^*$ and the coproduct vanishes.
\end{proof}
\begin{prop}
    For $k\geq 1$, the assignment $\varphi \mapsto \tau_\varphi$ is a group homomorphism $J^A(k) \to \Hom_A(K_{k+1},H)$.
\end{prop}
\begin{proof}
    Let $\varphi,\psi\in J^A(k)$ be two automorphisms of $\Gamma$. Then
    \[
    \tau_{\varphi\psi} = T-(\varphi\psi)^*(T) = T-\psi^*\varphi^*(T) = T-\psi^*(T) + \psi^*( T - \varphi^*(T) ).
    \]
    But since $T-\varphi^*(T)$ has weight $1$, it is fixed by $\psi^*$ and we get $\tau_{\varphi\psi} = \tau_{\psi}+ \tau_{\varphi}$.
\end{proof}
\begin{prop}
    The Johnson homomorphism $J^A(k)\to \Hom_A(K_{k+1},\H)$ is $\Aut(\Gamma)$-equivariant, where the action on $J^A(k)$ is by conjugation.
\end{prop}
\begin{proof}
    For any $\chi\in \Aut(\Gamma)$,
    \[
    \tau_{\chi \varphi \chi^{-1}}(T) = T - (\chi^{-1})^*\varphi^* \chi^*(T) = (\chi^{-1})^*(\chi^*(T) - \varphi^*(\chi^*(T))) =(\chi^{-1})^* \tau_\varphi(\chi^*(T))
    \]
    which is exactly the action of $\chi$ on $\Hom_A(K_{k+1},\H)$.
\end{proof}
\begin{prop}
    The kernel of the Johnson homomorphism $J^A(k)\to \Hom_A(K_{k+1},\H)$ contains $J^A(k+1)$. When $\Gamma$ is finitely generated and $A[\Gamma]/I^{k+2}$ is $A$-torsion free (e.g. when $A$ is a field), this is the entire kernel.
\end{prop}
\begin{proof}
    Any $\varphi\in J^A(k+1)$ stabilizes all of $\H^0_{Bar}(\Gamma;A)_{\leq k+1}$ and thus $\tau_{\varphi}=0$. Conversely, if $\varphi\in J^A(k)$ has $\tau_\varphi = 0$ then $T = \varphi^*(T)$ for all $T\in \H^0_{Bar}(\Gamma;A)_{\leq k+1}$. We claim that, for $\Gamma$ finitely generated with $A[\Gamma]/I^{k+2}$ that is $A$-torsion-free, $\varphi$ fixes all of $A[\Gamma]/I^{k+2}$.

    Indeed, since $\Gamma$ is finitely generated, the $A$-module $A[\Gamma]/I^{k+2}$ is also finitely generated and we assume it to be torsion-free, thus it is in fact free. If $x\neq \varphi_*(x) \in A[\Gamma]/I^{k+2}$, then there exists a functional $f\maps A[\Gamma]/I^{k+2}\to A$ such that $f(x-\varphi(x))\neq 0$. But, by completeness of letter-braiding, there exists $T\in \H^0_{Bar}(\Gamma;A)_{\leq k+1}$ such that $\langle T,\bullet \rangle = f$. Now,
    \[
    f(x) = \langle T,x\rangle = \langle \varphi^*(T),x\rangle = \langle T, \varphi_*(x) \rangle = f(\varphi_*(x))
    \]
contradicting our choice of $f$.   It therefore must be that $\varphi_* = \Id$ on $A[\Gamma]/I^{k+2}$. 
\end{proof}
Together, these propositions prove Theorem \ref{thm:intro-johnson}. We apply this theory to finite $p$-groups.

\begin{lemma}\label{lem:p-johnson is elementary abelian}
    Let $p$ be a prime number and $\Gamma$ a finite $p$-group. Then for every $k\geq 1$, the quotient $\sfrac{J^{\F_p}(k)}{J^{\F_p}(k+1)}$ is an elementary abelian $p$-group.
\end{lemma}
\begin{proof}
    $\Gamma$ is finitely generated and $\F_p$ is a field, so the Johnson homomorphism gives an embedding
    \[
    \sfrac{J^{\F_p}(k)}{J^{\F_p}(k+1)} \into \Hom_{\F_p}(K_{k+1},H)
    \]
    whose domain is an $\F_p$-vector space. It follows that the quotient is elementary abelian.
\end{proof}
\begin{lemma}
    For $\Gamma$ a finite $p$-group, there exists $k\geq 1$ such that $J^{\F_p}(k+1) = \{1\} \leq \Aut(\Gamma)$. The same holds for any $p$-nilpotent group, i.e. one for which the lower-exponent $p$-central series terminates, with $k$ being the $p$-nilpotency class.
\end{lemma}
\begin{proof}
    Fix $x\in \Gamma$ and recall that $(x-1)\in I^k\lhd \F_p[\Gamma]$ if and only $x$ belongs to $\gamma_k(\Gamma;\F_p)$, the $k$-th term in the lower-exponent $p$-central series for $\Gamma$ (see \ref{sec:group ring prelims}).
    
    Suppose $\gamma_{k+1}(\Gamma;\F_p) = \{1\}$. Then if $\varphi\in \Aut(\Gamma)$ acts trivially on $\F_p[\Gamma]/I^{k+1}$, it follows that for all $x\in \Gamma$,
    \[
    (\varphi(x)-x) \in I^{k+1} \implies (\varphi(x)x^{-1} - 1) \in I^{k+1} \implies \varphi(x)x^{-1}\in \gamma_{k+1}(\Gamma;\F_p) = \{1\}
    \]
    so $\varphi(x) = x$ and thus $\varphi = 1\in \Aut(\Gamma)$.
\end{proof}
With this, we complete the proof of Theorem \ref{thm:aut of p-groups}.
\begin{proof}[Proof of Theorem \ref{thm:aut of p-groups}]
    The mod $p$ Johnson filtration
    \[
    \{1\} \leq J^{\F_p}(k) \leq \ldots \leq J^{\F_p}(1) \leq J^{\F_p}(0) = \Aut(\Gamma)
    \]
    has finite length and by Lemma \ref{lem:p-johnson is elementary abelian} all consecutive quotients below $J^{\F_p}(1)$ are elementary abelian $p$-groups. It follows that $J^{\F_p}(1)$ is itself a finite $p$-group.

    In particular, it follows that any $\varphi\in \Aut(\Gamma)$ that does not have $p$-power order cannot belong to $J^{\F_p}(1)$, implying that it acts nontrivially on $\F_p[\Gamma]/I^2$. But then, by the same argument as in the proof of the last lemma, there exists some $x\in \Gamma$ such that $\varphi(x)x^{-1}\notin \gamma_2(\Gamma;\F_p) = [\Gamma,\Gamma]\Gamma^p$, so $\varphi(x) \neq x$ in the quotient $\Gamma^{ab}\otimes \F_p$. 
\end{proof}

\subsection{Massey products of finite groups} \label{sec:Massey in finite groups}
For finite groups, our theory leads to the following fact about Massey products in $\H^1(-;\F_p)$.
\begin{theorem}
    Let $\Gamma$ be a finite group, and $p$ any prime number such that $\H^1(\Gamma;\F_p)\neq 0$. Then $\H^2(\Gamma;\F_p)\neq 0$ and nontivial Massey products $\mu_k\maps \H^1(\Gamma;\F_p)^{\otimes k} \to \H^2(\Gamma;\F_p)$ exist.

    Quantitatively, denote $b_1 = \dim_{\F_p}\H^1(\Gamma;\F_p)$, if $k\geq 2$ is the length shortest nontrivial Massey product
    \[
    0\neq \langle \alpha_1,\ldots,\alpha_k\rangle \in \H^2(\Gamma;\F_p),
    \]
    which is necessarily well-defined with no indeterminacy, then
    \begin{equation}
        k < \begin{cases}
            |\Gamma|+1 & b_1 = 1 \\
            \log_{b_1}|\Gamma| + 2 & b_1 > 1.
        \end{cases}
    \end{equation}
\end{theorem}
\begin{proof}
    Since $\F_p[\Gamma]$ is an $\F_p$-vector space of dimension $|\Gamma|$, and the letter-braiding invariants $\H^0_{Bar}(\Gamma;\F_p)$ inject into its linear dual, they must also have dimension $\leq |\Gamma|$. Now, suppose $k\geq 2$ is the smallest length of a nontrivial Massey product in $\H^1(\Gamma;\F_p)$. Then by the weight spectral sequence, 
    \[
    \gr \H^0_{Bar}(\Gamma;\F_p)_{\leq k-1} \cong \F_p \oplus \H^1 \oplus (\H^1)^{\otimes 2} \oplus \ldots \oplus (\H^1)^{\otimes k-1},
    \]
    so in particular $|\Gamma|\geq \dim_{\F_p} \H^0_{Bar}(\Gamma;\F_p) \geq 1+b_1+b_1^2+\ldots+b_1^{k-1}$. The claimed inequalities for $k$ follow.
\end{proof}
    \begin{example}[Cyclic $p$-groups]\label{ex:finite cyclic}
    Let $\Gamma = C_p$ by the cyclic group of prime order $p$ and take $A=\F_p$. Then $C^1(BC_p;\F_p) = \F_p^{C_p}$ is spanned by the power functions $m_k\maps C_p \to \F_p$,
    \[
    m_k(i) = i^k,
    \]
    for $0\leq k \leq p-1$. Only $m_1$ is a group homomorphism, and accordingly, $H^1(C_p;\F_p)= \F_p\cdot m_1$.

    One can check directly from the definitions that the Massey products of $m_1$ with itself are
    \[
    \langle \underbrace{m_1,\ldots, m_1}_{k\text{ times}}\rangle = \frac{1}{k!}\sum_{j=1}^{k-1}\frac{k!}{j!(k-j)!}m_j\smile m_{k-j} = d\left(\frac{1}{k!}m_k\right)
    \]
    so vanish in cohomology. They thus define letter-braiding invariants
    \[
    I_k := \frac{1}{k!}\sum_{\substack{i_1+\ldots+i_r = k \\ \forall j,\,i_j>0}} \frac{k!}{i_1!\cdots i_r!} (m_{i_1}|m_{i_2}|\ldots |m_{i_r}) = (\underbrace{m_1|\ldots|m_1}_{k \text{ times}}) + \text{(lower weight terms)}.
    \]
    This invariant pairs with $C_p$ as
    \[
    \langle I_k, i\rangle = \frac{i^k}{k!} \mod{p},
    \]
    which by the multinomial theorem is independent of a presentation of $i$ as $a_1+\ldots+a_n$.
    
    As for the $p$-th Massey power of $m_1$, observe that $m_p = (-)^p$ is already a group homomorphism\footnote{In fact $m_p=m_1$ by Fermat's little theorem}, i.e. $d m_p = 0$, allowing the $p$-th Massey power of $m_1$ to be nontrivial.

    Similar analysis for the group $C_{p^2}$ shows that now the $p$-th Massey power of $m_1$ is $d\gamma_1$ for the divided power operation
    \[
    \gamma_1(i) = \frac{i^p-i}{p} \mod p.
    \]
    It then follows that the weight $p$ invariant of $C_{p^2}$ is exactly the function $\gamma_1$, which detects $p$.
    \end{example}

\bibliographystyle{alpha}
\bibliography{bibliography}

\end{document}